\documentclass[]{amsart}

\usepackage{amsmath,amssymb}
\usepackage{graphicx,extarrows}
\usepackage{xcolor}
\usepackage{amsfonts, mathrsfs}
\usepackage{epstopdf,algorithm,algorithmic}

\newtheorem{Def}{Definition}[section]
\newtheorem{ASSUM}{Assumption}[section]

\newtheorem{Lma}{Lemma}[section]
\newtheorem{Thm}{Theorem}[section]
\newtheorem{Prop}{Proposition}[section]
\newtheorem{Cor}{Corollary}[section]
\theoremstyle{remark}
\newtheorem{Rmk}{Remark}[section]

\begin{document}
\title[Convergence of Hessian estimator]
{Convergence of Hessian estimator from random samples on a manifold with boundary}

\author{Chih-Wei Chen}
\address{Department of Applied Mathematics, National Sun Yat-sen University, Taiwan}
\email{BabbageTW@gmail.com, chencw@math.nsysu.edu.tw}
\author{Hau-Tieng Wu}
\address{Courant Institute of Mathematical Sciences, New York University, New York, 10012, USA}
\email{hauwu@cims.nyu.edu}
\date{}

\begin{abstract}
A common method for estimating the Hessian operator from random samples on a low-dimensional manifold involves locally fitting a quadratic polynomial. Although widely used, it is unclear if this estimator introduces bias, especially in complex manifolds with boundaries and nonuniform sampling. Rigorous theoretical guarantees of its asymptotic behavior have been lacking. We show that, under mild conditions, this estimator asymptotically converges to the Hessian operator, with nonuniform sampling and curvature effects proving negligible, even near boundaries. Our analysis framework simplifies the intensive computations required for direct analysis.
\end{abstract}

\maketitle

\section{Introduction}

In modern data analysis, researchers frequently encounter datasets characterized by high-dimensionality and nonlinearity. A common approach to modeling such datasets involves assuming that the point cloud is randomly sampled from a low-dimensional manifold or collected using a structured sampling scheme. Under this model, the objective of data analysis is to characterize the manifold structure for tasks such as inference and prediction. This paper focuses specifically on estimating the Hessian operator on manifolds with boundaries from random, potentially non-uniform samples. 

The Hessian, as a second-order differential operator on the manifold, has been widely considered  in various algorithms \cite{DG03,KSH09,SGWJ18,SJWG20}. On the theoretical side, the Hessian operator is intimately related to the Laplace-Beltrami operator and finds applications in studying distance functions on manifolds \cite{P16}, the geometry of manifolds with density or metric measure spaces  \cite{Mo05,S06a,S06b}, among others. In applications, it has been used in unsupervised manifold learning techniques like Hessian locally linear embedding (HLLE) \cite{DG03} and geodesic distance estimation \cite{CZXC20}, semi-supervised learning \cite{KSH09},  computer graphics \cite{HWAG09,SGWJ18,SJWG20}, scientific computation \cite{ITZ22}, and generally the optimization, where researchers have found it useful for capturing geometric features of data. 

There are different approaches to estimate the Hessian. We focus on the common one that involves fitting a quadratic polynomial to functions based on Taylor expansion principles and obtaining the second-order expansion as the Hessian estimator. This principle can be found in various fields of data analysis such as Savitzky-Golay filter \cite{SG64} or local polynomial regression \cite{FG96} in the Euclidean space setup, or its intuitive generalization to the manifold setup.
Despite empirical success, there is a lack of theoretical justification for how the Hessian estimator converges when random samples are independently and identically distributed (i.i.d.) according to a nonuniform density supported on a low-dimensional manifold with complex geometry and topology. Specifically, due to the intricate interplay among sampling density, curvature, and boundary effects, the behavior of the estimator remains unclear. A natural question arises: is such a Hessian estimator biased, and how does it asymptotically converge to the true Hessian operator as the number of sampled points increases, particularly in scenarios involving nontrivial curvature, nonuniform sampling, and non-empty boundaries?

It is worth noting that theoretical results on Hessian estimation in Euclidean settings can be found in \cite{FG96} within the framework of local polynomial regression, and in the flat manifold setup where the manifold $M$ is isometric to an open connected subset in $\mathbb R^d$ and sampling is uniform, as shown in \cite{DG03}. We shall also mention similar results that there is extensive theoretical support for estimating gradients, a first-order differential operator, from random samples under both manifold \cite{MWZ10,WGMM10,CW13} and Euclidean settings \cite{FG96} under the same Taylor expansion principle.

Our main contribution in this article is to demonstrate that the Hessian estimator, based on tangent space estimation using local principal component analysis (PCA) followed by polynomial fitting, asymptotically converges to the continuous Hessian on $M$. This convergence holds when the point cloud is i.i.d. sampled from a low-dimensional Riemannian manifold, with or without a boundary, and potentially nonuniformly.
Specifically, regardless of whether we estimate the Hessian near or away from the boundary, the impacts of nonuniformity and curvature are negligible as they appear in higher-order terms compared with the Hessian operator, albeit with a slower convergence rate near the boundary. These finds are detailed in Theorem \ref{main}.

The primary technical challenge addressed in Theorem \ref{main} is managing the intricate interplay among curvature, nonuniform density, and the function under analysis when fitting a quadratic function, where an inversion of a Gram matrix is involved. Specifically, for a manifold of dimension $d$, estimating the Hessian at point $x$ with samples within a local ball of radius $\varepsilon>0$ requires inverse the Gram matrix of size $(1+2d+\frac{d(d-1)}{2})\times (1+2d+\frac{d(d-1)}{2})$. This Gram matrix has the form $ZZ^T$, where $Z$ is the base matrix (or called design matrix in the statistics literature) that involves the constant, linear and quadratic coordinates. In general, the entries of $ZZ^T$  feature complex interaction terms of mixed orders of $\varepsilon$ as $\varepsilon\to 0$ asymptotically. See \eqref{Expansion of Z^TZ for the direct approach} and \eqref{Expansion of Z^Tf for the direct approach} for a depiction of this complex structure. Although a direct expansion is theoretically feasible, its execution is excessively intricate. The complexity heightens when estimating the Hessian near the boundary, where additional terms arise absent in boundary-free estimations. To handle this technical challenge, we propose a reduction trick in this work. In essence, we separate the function of interest to its dominant terms and higher order approximation first before applying fitting techniques, thereby significantly streamlining the analysis effort, even in scenarios near the boundary. See \eqref{reduction trick main equation} for details.

We shall mention that in quadratic kernel regression \cite{FG96}, the Hessian is typically employed to correct fitting errors and achieve more accurate function approximations, with primary focus on the function itself as the zeroth-order term through local approximation. This analysis thus prioritizes a few top-order terms, provided the zeroth-order term remains intact. When the focus shifts to the first-order term, the gradient \cite{MWZ10,WGMM10,CW13}, the analysis entails higher order terms. Luckily, in this case the interaction of higher order terms are simple. The left upper $2\times 2$ block matrix in \eqref{L0} is the Gram matrix used in the gradient estimate, where the dominant structure is diagonal, except for an additional off-diagonal term when boundary effects are factored in. However, when dealing the Hessian, which is the second-order term,  the complexity escalates significantly. Tracking higher-order terms becomes critical as they begin to interact in non-trivial ways. 
The more higher-order terms that require tracking, the greater the analytical complexity. 
Moreover, it is clear that the more higher order terms we track, the more complicated the computation will be.
It is noteworthy that the manifold model used in \cite{DG03} remains curvature-free without boundaries, and with uniform sampling, thus minimizing potential interactions between non-uniform sampling, curvature, and boundary effects. 
Our Theorem \ref{main} is thus an extension of the work in \cite{FG96,DG03,MWZ10,WGMM10,CW13}.

In summary, to our knowledge, our result is the first quantitative guarantee for the Hessian estimator from random samples, offering an explicit convergence rate that applied to non-flat data manifolds with boundaries and with nonuniform sampling.

The article is organized as follows. In Section 2, we illustrate the derivation of the discretized Hessian matrix through quadratic fitting in both the Euclidean and manifold settings. Section 3 offers a theoretical justification of our proposed algorithm. A comprehensive proof of the main theorem is provided in the Appendix. We provide two versions of proofs: one utilizing direct expansion and another employing the reduction trick to showcase its advantages.
The paper concludes with a discussion and conclusion section, highlighting the profound connection between HLLE and a complex fourth-order differential equation.

We will systematically use the following notations. When a function $h(t)$ satisfies $c_1\leq |h(t)|t^{-\lambda}\leq c_2$ for some $0<c_1\leq c_2$ and $\lambda\in \mathbb{R}$ as $t\to 0$, we denote $h=\Theta(t^{\lambda})$. When $0=c_1<c_2$, we denote $h=O(t^{\lambda})$.

\section{Quadratic fitting via Taylor's expansion}

Since the Hessian is characterized by quadratic approximation in Taylor's expansion, it is common practice to utilize Taylor's expansion to construct a Hessian estimator for functions defined on a data manifold $M$. To introduce this concept and establish notation, we begin with the Euclidean case.

\subsection{Euclidean case}
Considering the Euclidean case $M=\mathbb{R}^d$. Let $f:\mathbb{R}^d\to \mathbb{R}$ be a $C^2$-function around $x$. Without loss of generality, we assume $x=0$. Let ${\bf q}_1,\cdots,{\bf q}_k$ be $k$ points in $\mathbb{R}^d$, and we write ${\bf q}_j=[({\bf q}_j)_1, \cdots, ({\bf q}_j)_d]^T\in \mathbb{R}^d$ for all $j$. In terms of the standard basis, we denote
$$\left[\begin{array}{c}
{\bf q}_1^T\\
\vdots\\
{\bf q}_k^T
\end{array}\right]_{k\times d}
=\left[\begin{array}{ccc}
({\bf q}_1)_1& \cdots& ({\bf q}_1)_d\\
\vdots & \ddots& \vdots\\
({\bf q}_k)_1& \cdots& ({\bf q}_k)_d
\end{array}\right]
=
\left[\begin{array}{ccc}
({\bf y}_1)_1& \cdots& ({\bf y}_d)_1\\
\vdots & \ddots& \vdots\\
({\bf y}_1)_k& \cdots& ({\bf y}_d)_k
\end{array}\right]
=
\left[\begin{array}{ccc}
{\bf y}_1& \cdots& {\bf y}_d
\end{array}\right]_{k\times d},
$$
that is, ${\bf y}_i\in \mathbb{R}^k$ records the $i$-th coordinates of all $k$ vectors.
When ${\bf y}_j$ are all close to $0$, Taylor's expansion of $f$ at $0$ gives
\begin{align*}
\left[\begin{array}{ccc}
f({\bf q}_1)\\
\vdots\\
f({\bf q}_k)
\end{array}\right]
\approx &\
\left[\begin{array}{c}
 f(0) + (\nabla f(0))^T{\bf q}_1+\frac{1}{2}{\bf q}_1^T  H_f{\bf q}_1 \\
 \vdots\\
 f(0) + (\nabla f(0))^T{\bf q}_k+\frac{1}{2}{\bf q}_k^T  H_f {\bf q}_k 
\end{array}\right]\\
=&\
\left[\begin{array}{cccc}
{\bf 1}_{k\times 1}& 
({\bf y}_1\ \cdots \ {\bf y}_d )_{k\times d} &
({\bf y}_s \circ{\bf y}_s)_{1\leq s\leq d} &
({\bf y}_s \circ{\bf y}_t)_{1\leq s< t\leq d}  
\end{array}\right]_{k\times(1+d+d+\frac{d(d-1)}{2})}\\
&\ \cdot\left[\begin{array}{cccc}
 f(0)& \nabla f(0)_{1\times d}& \frac{1}{2}(h_{ss})_{1\leq s\leq d} &
(h_{st})_{1\leq s< t\leq d}
\end{array}\right]_{1\times(1+d+d+\frac{d(d-1)}{2})}^T,
\end{align*}
where $H_f=(h_{st})=\left(\frac{\partial^2f}{\partial x_s\partial x_t}\big|_0\right)\in \mathbb{R}^{d\times d}$ is the Hessian of $f$ at the origin, each ${\bf y}_s \circ{\bf y}_t$ is the Hadamard product of ${\bf y}_s$ and ${\bf y}_t$, namely the $k$-dimensional vector
$\left[({\bf y}_s)_1({\bf y}_t)_1\ \ ({\bf y}_s)_2({\bf y}_t)_2\ \ \cdots\ \ ({\bf y}_s)_k({\bf y}_t)_k\ \right]^T$, and $[({\bf y}_s \circ{\bf y}_s)_{1\leq s\leq d} \ \ \
({\bf y}_s \circ{\bf y}_t)_{1\leq s< t\leq d}]$ is a $k\times \frac{d(d+1)}{2}$ matrix. In particular, for any $j\in\{1,2,\dots,k\}$, we see that
\begin{align}
&\sum_{1\leq s\leq d}({\bf y}_s\circ{\bf y}_s)_j\frac{1}{2}h_{ss}+\sum_{1\leq s<t\leq d}({\bf y}_s\circ{\bf y}_t)_jh_{st}\\
=&\,\sum_{1\leq s\leq d}({\bf q}_j)_s({\bf q}_j)_s\frac{1}{2}h_{ss}+\sum_{1\leq s<t\leq d}({\bf q}_j)_s({\bf q}_j)_th_{st}={\bf q}_j^T  H_f {\bf q}_j\,.\nonumber
\end{align}
Throughout this paper, we introduce the following notations to simplify the discussion. We define the {\it Hessian vector} of $f$ at $0$ by 
\[
K_f:=\left[\frac{1}{2}h_{11}\ \cdots\ \frac{1}{2}h_{dd}\ h_{12}\ h_{13}\cdots\  h_{(d-1)(d)}\right]\,, 
\]
which is essentially a lineup of entries in the Hessian matrix $H_f$.
Moreover, we denote $Z \in \mathbb{R}^{k\times(1+d+d+\frac{d(d-1)}{2})}$ as
$$Z=\left[Z_A\ Z_B\ Z_C\ Z_D \right]=\left[\begin{array}{cccc}
{\bf 1}_{k\times 1
}& 
({\bf y}_1\ \cdots \ {\bf y}_d )_{k\times d} &
({\bf y}_s\circ {\bf y}_s)_{1\leq s \leq d} &
({\bf y}_s\circ {\bf y}_t)_{1\leq s<t\leq d}  
\end{array}\right]$$
as the {\it base matrix}, where
$$Z_A={\bf 1}_{k\times 1}, Z_B=\left[\begin{array}{ccc}
|&&|\\
{\bf y}_1&\cdots &{\bf y}_d\\
|&&|
\end{array}\right]
=\left[\begin{array}{ccc}
\mbox{---}&{\bf q}_1&\mbox{---}\\
&\vdots &\\
\mbox{---}&{\bf q}_{k}&\mbox{---}
\end{array}\right],$$
$$ Z_C=\left[\begin{array}{ccc}
|&&|\\
{\bf y}_1\circ{\bf y}_1&\cdots &{\bf y}_d\circ {\bf y}_d\\
|&&|
\end{array}\right],\mbox{ and }
Z_D=\left[\begin{array}{ccc}
|&&|\\
{\bf y}_1\circ{\bf y}_2&\cdots &{\bf y}_{d-1}\circ {\bf y}_d\\
|&&|
\end{array}\right].$$

Note that, by Taylor's approximation, $K_f$ is the vector which minimizes 
$$\left\| 
\left[\begin{array}{ccc}
f({\bf q}_1)& \cdots &f({\bf q}_k)
\end{array}\right]^T
-Z \left[\begin{array}{ccc}
f(0)& \nabla f(0)_{1\times d}& K
\end{array}\right]^T\right\|_2$$
among all $\frac{d(d+1)}{2}$-vectors $K$.
Hence, when given $f({\bf q}_1), \dots, f({\bf q}_k)$, and $Z$, one can solve the minimization problem $$\min_{G\in \mathbb{R}^{(1+d+\frac{d(d+1)}{2})\times 1} }\left\| \left[\begin{array}{ccc}
f({\bf q}_1)& \cdots &f({\bf q}_k)
\end{array}\right]^T
-ZG\right\|_2^2$$
and thus obtain an estimate of the Hessian vector $K_{f}$ from the solution $G_*$. Therefore, $K_f$ is our discrete Hessian estimator. 
It is well-known that the minimization can be solved with the solution as the normal equation $$G_*=(Z^TZ)^{-1}Z^T \left[\begin{array}{ccc}
f({\bf q}_1)& \cdots &f({\bf q}_k)
\end{array}\right]^T,$$
as long as $Z^TZ$ is invertible. The convergence of this estimator can be found in \cite{FG96,DG03}. 

\subsection{Manifold case}
From now on we assume that our data is stored in $\mathbb{R}^p$ and distributes in a $d$-dimensional Riemannian submanifold $M$ with $d<p$; that is, $M$ is isometrically embedded in $\mathbb{R}^p$ via an inclusion map $\iota$. 
Consider a point cloud $\{x_i\}_{i=1}^n\subset \iota(M)\subset  \mathbb{R}^p$ that is sampled i.i.d. from a random variable with the range $\iota(M)$. 
Consider a  point $z\in M$ and the Euclidean ball $B_{\varepsilon}(\iota(z))\subset\mathbb{R}^p$ centered at $\iota(z)$. 
Denote $B_{\varepsilon}(\iota(z))\cap \{x_i\}_{i=1}^n=\{x_{z,j}\}_{j=1}^{k_z}$, where $k_z$ is the number of sample points lying in $B_{\varepsilon}(\iota(z))$. We call $B_{\varepsilon}(\iota(z))\cap \iota(M)$ the $\varepsilon$-neighborhood and $x_{z,1},\dots,x_{z,k_z}$ the $\varepsilon$-neighbors of $z$. 
Let $T_{z}M$ be the tangent space of $M$ at $z$ and $\iota_*(T_{z}M)$ the embedded tangent space in $\mathbb{R}^p$. Denote the projection of $x_{z,j}$ on  $\iota_*(T_{z}M)\subset\mathbb{R}^p$ as $\widetilde{\bf q}_j$, where $j=1,\dots,k_z$. 
By choosing an orthonormal basis $\{e_1,\dots,e_d,e_{d+1},\dots,e_p\}$ of $\mathbb{R}^p$, where $e_1,\dots,e_d\in \iota_*(T_{z}M)$ and $e_{d+1},\dots,e_p\in (\iota_*(T_{z}M))^\perp$,
we can express 
$$
{\widetilde{\bf q}}_j=((\widetilde{\bf q}_j)_1,\dots,(\widetilde{\bf q}_j)_d,0,\dots,0)\mbox{ and } x_{z,j}=((\widetilde{\bf q}_j)_1,\dots,(\widetilde{\bf q}_j)_d, (x_{z,j})_{d+1}, \dots, (x_{z,j})_{p}).
$$
We shall clarify the notation a bit. While we shall use the embedding map $\iota$ to denote the submanifold as $\iota(M)\subset\mathbb{R}^p$, however, to alleviate the notational burden, we usually omit the notation $\iota$ and simply denote $\iota(M)$ by $M$ when there is no danger of confusion. Similarly, we use $z\in M$ and $\iota(z)\in \mathbb{R}^p$ interchangeably, and omit the notation $\iota_*$ when we discuss the tangent space because we always consider the tangent space as an affine subspace embedded in $\mathbb{R}^p$ and identify its origin to the point $z$.

Take a $C^2$-function $f:M\to\mathbb{R}$. In view of the fact that locally a manifold can be well approximated by an affine subspace, motivated by the Hessian estimate in the Euclidean case we have discussed above, we could estimate the Hessian at $z$ by the same way considered in the Euclidean setup via evaluating 
\begin{align}\label{G=ZTZ-1ZTf}
\widetilde G:=(\widetilde Z^T\widetilde Z)^{-1}\widetilde Z^T {\bf f},
\end{align}
where ${\bf f}:=\left[\begin{array}{ccc}
f(x_{z,1})& \cdots &f(x_{z,k_z})
\end{array}\right]^T\in \mathbb{R}^{k_z}$ is a discretization of $f$, 
and
$\widetilde Z=\left[\widetilde Z_A\ \ \widetilde Z_B\ \ \widetilde Z_C\ \ \widetilde Z_D \right]\in \mathbb{R}^{k_z\times(1+d+d+\frac{d(d-1)}{2})}$ is also called the {\em base matrix associated with $z$} with  

\begin{align}\label{definition Zb y p}
\widetilde Z_A={\bf 1}_{k_z\times 1},\ 
\widetilde Z_B=\left[\begin{array}{ccc}
|&&|\\
\widetilde{\bf y}_1&\cdots &\widetilde{\bf y}_d\\
|&&|
\end{array}\right]
=\left[\begin{array}{ccc}
\mbox{---}&\widetilde{\bf q}_1 &\mbox{---}\\
 &\vdots& \\
\mbox{---}&\widetilde{\bf q}_{k_z} &\mbox{---}
\end{array}\right], 
\end{align}

\begin{align*}
\widetilde Z_C=\left[\begin{array}{ccc}
|&&|\\
\widetilde{\bf y}_1\circ\widetilde{\bf y}_1&\cdots &\widetilde{\bf y}_d\circ \widetilde{\bf y}_d\\
|&&|
\end{array}\right],\mbox{ and }
\widetilde Z_D=\left[\begin{array}{ccc}
|&&|\\
\widetilde{\bf y}_1\circ\widetilde{\bf y}_2&\cdots &\widetilde{\bf y}_{d-1}\circ \widetilde{\bf y}_d\\
|&&|
\end{array}\right],
\end{align*}
where $\widetilde{\bf y}_i\in \mathbb{R}^{k_z}$ records the $i$-th coordinates of all $k_z$ neighboring points associated with the chosen basis.

However, in practice, we do not have an access to the manifold parametrization, nor the tangent space. Instead, we have only the point cloud $\{x_i\}_{i=1}^n$. Thus, to apply this idea, we need to estimate $T_{z}M$, particularly an orthonormal basis of $T_{z}M$. Since the full information of the manifold continuum is unknown, we can only estimate the tangent space $T_{z}M$ by using sample points near $z$. This can be done by the local principal component analysis (PCA) idea \cite{SW12,KM14,AL19}. Indeed, one may perform PCA on the $\varepsilon$-neighbors of $z$, $\{x_{z,j}\}_{j=1}^{k_z}$, and obtain an orthonormal basis by the first $d$ dominant principal directions, denoted as $\mathfrak{U}_z:=\{u_l\}_{l=1}^d$, that spans a $d$-dimensional vector space, denoted as $V_z$, as an estimate of $T_{z}M$. With $\mathfrak{U}_z$, we could project the samples within the $\varepsilon$-neighborhood of $z$ onto $V_z$, and obtain estimated coordinates associated with $\mathfrak{U}_z$.
Specifically, we denote the projection of $x_{z,j}$'s by ${\bf q}_j\in \mathbb{R}^d$. 
Note that the estimated coordinate ${\bf q}_j$ is expected to be close to $\widetilde{\bf q}_j$ by the measurement in $\mathbb{R}^p$ if the orthonormal basis of $T_{z}M$, denoted as $\{e_l\}_{l=1}^d$, is properly chosen. Indeed, it has been stated in Theorem B.1 in \cite{SW12} that there exists an orthonormal basis $\{e_l\}_{l=1}^d$ of $T_{z}M$ so that
\begin{align}\label{pq}
(\widetilde{\bf q}_j)_l=\langle x_{z,j}-z, e_l \rangle =\langle x_{z,j}-z,u_l \rangle +O(\varepsilon^{\tau+2}) 
=({\bf q}_j)_l+O(\varepsilon^{\tau+2})
\end{align}
when $M$ is a manifold with boundary, where $\tau=\frac{1}{2}$ when $z$ is close to the boundary and $\tau=1$ when $z$ is away from the boundary.

With the estimated coordinates of neighboring points, ${\bf q}_j\in V_z$, we consider an estimate of the base matrix at $z$, denoted as 
\begin{align}\label{formula: Zhat}
 Z=\left[ Z_A\ \  Z_B\ \  Z_C\ \  Z_D \right]\,, 
\end{align}
where
$$
 Z_A={\bf 1}_{k_z},\ 
 Z_B
=\left[\begin{array}{ccc}
\mbox{---}&{\bf q}_1&\mbox{---}\\
&\vdots &\\
\mbox{---}&{\bf q}_{k_z}&\mbox{---}
\end{array}\right]=\left[\begin{array}{ccc}
|&&|\\
{\bf y}_1&\cdots &{\bf y}_d\\
|&&|
\end{array}\right]\in \mathbb{R}^{k_z\times d},
$$
$$  Z_C=\left[\begin{array}{ccc}
|&&|\\
{\bf y}_1\circ{\bf y}_1&\cdots &{\bf y}_d\circ {\bf y}_d\\
|&&|
\end{array}\right]\in \mathbb{R}^{k_z\times d}
$$
and 
$$
 Z_D=\left[\begin{array}{ccc}
|&&|\\
{\bf y}_1\circ{\bf y}_2&\cdots &{\bf y}_{d-1}\circ {\bf y}_d\\
|&&|
\end{array}\right]\in \mathbb{R}^{k_z\times d(d-1)/2},
$$
where ${\bf y}_i\in \mathbb{R}^{k_z}$ records the estimated $i$-th coordinates of all $k_z$ neighboring points associated with the estimated  tangent space basis by local PCA.
Note that the notation here is abused to be coincide with the Euclidean case in Section 2.1. The only difference is that $\{{\bf q}_j\}_{j=1}^{k_z}$ in Section 2.1 could be any points around the center, while the specific $\{{\bf q}_j\}_{j=1}^{k_z}$ here are derived from local PCA.
With the base matrix $Z$, the commonly used estimator of Hessian of $f$ at $z$ is via evaluating
\[
(Z^TZ)^{-1}Z^T{\bf f}\,.   
\]
Note that we need to assume $k_z\geq 1+d+\frac{d(d+1)}{2}$ to avoid the invertibility issue of $Z^TZ$.  Specifically, we have the following definition of the Hessian estimator from random samples.

\begin{Def}
Suppose $X$ is a random vector with the support on $M\subset\mathbb{R}^p$ and $\{x_i\}_{i=1}^n\subset M$ are $n$ points i.i.d. sampled from $X$. 
Let $f:M\to\mathbb{R}$ be a $C^{2}$-function.
Fix $z\in M$. Denote all points in $z$'s $\varepsilon$-neighborhood $B_{\varepsilon}(z)\subset\mathbb{R}^p$ as $x_{z,j}$, where $j=1,\dots, k_z$ and $k_z\in \mathbb{N}$ is the number of $\varepsilon$-neighbors of $z$. 
Denote the base matrix generated by local PCA as $Z$. Rewrite
$$
(Z^TZ)^{-1}Z^T=\left[\begin{array}{c}
{\mathcal A}\\
{\bf grad}\\
{\mathcal H}
\end{array}\right],
$$
where ${\mathcal A}$, ${\bf grad}$, and ${\mathcal H}$ be the $1\times k_z$, $d\times k_z$, and $\frac{d(d+1)}{2}\times k_z$ matrices, respectively.
Then ${\bf grad}$ is an estimator of the {\rm gradient} at $z$ and ${\bf Hess}\in \mathbb{R}^{d\times d\times k_z}$, defined by
\[
{\bf Hess}_{i,j,m}:=\left\{
\begin{array}{ll}
2\mathcal{H}_{i,m} & \mbox{if }i=j\\
\mathcal{H}_{d+j-i+(i-1)(d-(i-1)),m} & \mbox{if }i<j
\end{array}\right.
\mbox{ and } {\bf Hess}_{i,j,m}={\bf Hess}_{j,i,m},
\] 
is an estimator of the {\rm Hessian} at $z$. Hence,
$$
(Z^TZ)^{-1}Z^T{\bf f}=\left[\begin{array}{c}
\mathcal{A}{\bf f}\\
{\bf grad\, f}\\
{\mathcal{H} {\bf f}}
\end{array}\right]\in \mathbb{R}^{(1+d+d(d+1)/2)\times 1}
$$
gives estimates of $f(z)$, $\nabla f(z)$ and $Hess f(z)$ via
$\mathcal{A}{\bf f}$, ${\bf grad\, f}$, and ${\bf Hess\, f}$ respectively,
where ${\bf Hess\, f}\in \mathbb{R}^{d\times d}$ is the symmetric matrix defined as
\[
{\bf Hess\, f}_{i,j}:=\left\{
\begin{array}{ll}
2({\mathcal{H}{\bf f}})_{i} & \mbox{if }i=j\\
({\mathcal{H}{\bf f}})_{d+j-i+(i-1)(d-(i-1))} & \mbox{if }i<j
\end{array}\right..
\] 
\end{Def}

The Hessian estimator algorithm is summarized in Algorithm \ref{alg:Hess ours}. Note that ${\bf grad}\in \mathbb{R}^{d\times k_z}$ and ${\bf Hess}\in \mathbb{R}^{d\times d\times k_z}$ at $z$ are associated with the estimated basis of the tangent space $T_zM$ determined by local PCA. In the next section, we will show that under proper conditions, ${\bf grad}\ {\bf f}$ converges to the vector consisting of partial derivatives $\frac{\partial}{\partial x^l} f$ with respect to the normal coordinates $\{x^l\}_{l=1}^d$, which fits $\{e_l\}_{l=1}^d$ at $z$, and ${\bf Hess}\, {\bf f}$ converges to the Hessian of $f$ at $z$ with respect to the same normal coordinates. This result justifies the nomination.

\begin{algorithm}
\caption{Hessian estimator}
\label{alg:Hess ours}
\begin{algorithmic}
\STATE{{\bf Input}\  $z\in M^d\subset\mathbb{R}^p$ and $\{x_i\}\subset\mathbb{R}^p$, where $i=1,\dots,n$.
}
\STATE{Step 1: Find $k$ $\varepsilon$-neighbors of $z$, $x_{z,1},\dots,x_{z,k}$.}
\STATE{Step 2: Run {\bf SVD}\ $X_{p\times k}=\left[x_{z,1}-z\ \cdots\ x_{z,k}-z\right]=U\Lambda V^T$.}
\STATE{Step 3: Construct $U_d = \left[u_1\ \cdots\ u_d\right]$, where $u_1, \dots, u_d$ are corresponding to the $d$ largest singular values.}
\STATE{Step 4: Set $X^TU_d$ as the first $d$ projected coordinates of every $x_{z,j}$.}
\STATE{Step 5: Set $Z$.}
\STATE{Step 6: Evaluate $(Z^TZ)^{-1}Z^T$.}
\STATE{{\bf Output}\ 
{\bf grad}, {\bf Hess}}
\end{algorithmic}
\end{algorithm}

Before closing this section, we shall mention the relationship between the considered Hessian estimator and other existing algorithms. For example, in the Hessian eigenmap \cite{DG03}, the calculation of $(Z^TZ)^{-1}$ in the definition of ${\bf Hess}$ is carried out by applying the Gram-Schmidt process to columns of the base matrix $Z$ so that the resulting matrix $\mathsf{Z}$ satisfies $\mathsf{Z}^T\mathsf{Z}=I$. Thus, ${\bf Hess}\, {\bf f}=(Z^TZ)^{-1}Z^T{\bf f}=\mathsf{Z}^T{\bf f}$ because $(Z^TZ)^{-1}Z^T$ is the projection map from $\mathbb{R}^{k_z}$ to the column space of $Z$ and the Gram-Schmidt process does not change the column space. 
In \cite{SGWJ18,SJWG20}, the authors estimate the Hessian by composing the gradient operator with a matrix of divergence operator.

\section{Convergence of the Hessian estimator with rates}

Before stating our result, we impose assumptions about the manifold and sampling scheme. 
Take a random vector $X:(\Omega,\mathcal{F},\mathbb{P})\rightarrow\mathbb{R}^p$ that we will sample from. 

\begin{ASSUM}\label{Assumption1}
Assume that the range of $X$ is supported on a $d$-dimensional compact smooth 
Riemannian manifold $(M, g)$ that is isometrically embedded in $\mathbb{R}^p$ via $\iota:M\hookrightarrow \mathbb{R}^p$. The manifold may have boundary and we denote $\overline{M}=M\cup\partial M$. When the boundary exists, we assume it is smooth.
\end{ASSUM}

Denote $M_\sigma$, where $\sigma>0$, to be the $\sigma$-neighborhood of $\partial M$ in $M$ defined as
\[
M_\sigma:=\{x\in M|d(x,\partial M)\leq \sigma\}\,.
\]
The random vector $X$ induces a probability measure supported on $\iota(M)$, denoted by $\widetilde{\mathbb{P}}_{X}$.

\begin{ASSUM}
Assume that $\widetilde{\mathbb{P}}_{X}$ is absolutely continuous with respect to the Riemannian measure on $\iota(M)$, denoted by $\iota_{\ast}d{\rm vol}(x)$, which by the Radon-Nikodym theorem leads to $d\widetilde{\mathbb{P}}_{X}(x)=\rho_{X}(\iota^{-1}(x))\iota_{\ast}d{\rm vol}(x)$ for a nonnegative function $\rho_X$ defined on $M$. 
\end{ASSUM}

We call $\rho_{X}$ defined above the probability density function (p.d.f.) associated with $X$. When $\rho_{X}=1$, we call $X$ uniform; otherwise  nonuniform.
Note that since $\widetilde{\mathbb{P}}_{X}$ is an induced probability measure, we immediately have ${\rm Vol}_{\rho_X}(M):=\int_M\rho_{X}(\iota^{-1}(x))\iota_{\ast}d{\rm vol}(x)=1$.

\begin{ASSUM}\label{Assumption3}
Assume that $\rho_{X}$ satisfies $\rho_{X} \in  C^2(M)$ and $\inf_{x\in M}\rho_{X}(x)>0$.
\end{ASSUM}

\begin{ASSUM}\label{Assumption4}
Assume that the observed data set  $\mathcal{X} = \{x_{i}\}_{i=1}^{n}\subset\mathbb{R}^{p}$ is i.i.d. sampled from $X$.
\end{ASSUM}

We adopt normal coordinates $\{x^j\}$ and the basis $\{e_j:=\frac{\partial}{\partial x^j}\}$ of $T_{\xi}M$ to compare ${\bf Hess}\, {\bf f}$ and $Hess(f)$. 
Note that $Hess(f)(\xi)$ depends on the Christoffel symbol of $M$, which can be set to 0 at $\xi$ when the normal coordinates centered at $\xi$ are used. See Section \ref{section:hessian review} for a quick summary. On the other hand, the computation of ${\bf Hess}\, {\bf f}$ involves numerical approximation, and the curvature is involved in the tangent space estimate. Hence ${\bf Hess}\, {\bf f}$ and $Hess(f)(\xi)$ might differ with the deviation term involving curvatures. The question is how large these deviations are, and whether they bias the Hessian estimator. Moreover, when the sampling distribution is nonuniform, it is expected to also play a role, and we need to know how much it impacts the estimator.

The main theorem in this paper quantifies how well ${\bf Hess}\, {\bf f}$ approximates $Hess(f)$ at a given point asymptotically when $f\in C^{2,\kappa}(M)$ for $\kappa\in(0,1]$ as $n\to \infty$.
Recall that $f\in C^{2,\kappa}(M)$ indicates that $f\in C^2$ and the second derivatives of $f$ are $\kappa$-H\"older continuous, i.e., 
$\displaystyle\sup_M \max_{|\alpha|=2}\frac{|D^{\alpha}f(x)-D^{\alpha}f(y)|}{dist(x,y)^\kappa}<\infty$.
Note that when the domain $M$ is compact with smooth boundary, $f\in C^3$ implies that $f\in C^{2,\kappa}$ for all $\kappa\in(0,1]$ (cf. \cite[Lemma 4.28]{AF03}). See \eqref{Taylor expansion of C2kappa function} for the associated Taylor expansion.

\begin{Thm}[Quadratic fitting theorem]\label{main}
Assume Assumptions \ref{Assumption1}-\ref{Assumption4} hold for the sample data $\{x_i\}_{i=1}^n$.
Let $f:M\to\mathbb{R}$ be a $C^{2,\kappa}$-function, where $\kappa\in (0,1]$.
Suppose $\varepsilon=\varepsilon(n)$ so that $\varepsilon\to 0$ and $\frac{n\varepsilon^d}{\log(n)}\to\infty$ when $n\to \infty$. 
Take $\sigma=\sqrt{\varepsilon}$ and $M_\sigma$ be the $\sigma$-neighborhood of $\partial M$ in $M$. 
Consider $z\in \overline{M}=M\cup\partial M$ and denote the sample points in the $\varepsilon$-neighborhood $B_{\varepsilon}(z)\subset\mathbb{R}^p$ as $x_{z,j}$, where $j=1,\dots, k_z$. 
Denote the base matrix generated by local PCA as $ Z$. Then, when $n$ is sufficiently large, with probability greater than $1-O(n^{-3})$, we have 
$$( Z^T Z)^{-1} Z^T\left[\begin{array}{c}
f(x_{z,1})\\
\vdots \\
f(x_{z,k_z})
\end{array}\right]
=\left[\begin{array}{c}
f(z)+O(\varepsilon^{2+\min(\tau,\kappa)})+O(\sqrt{\frac{\log(n)}{n\varepsilon^{d}}})\\
\\
\nabla_1 f|_{z}+O(\varepsilon^{1+\min(\tau,\kappa)})+O(\sqrt{\frac{\log(n)}{n \varepsilon^{d+2}}})\\ 
\vdots\\ 
\nabla_d f|_{z}+O(\varepsilon^{1+\min(\tau,\kappa)})+O(\sqrt{\frac{\log(n)}{n\varepsilon^{d+2}}})\\
\\
\frac{1}{2}\nabla_1\nabla_1 f|_{z}+O(\varepsilon^{\min(\tau,\kappa)})+O(\sqrt{\frac{\log(n)}{n \varepsilon^{d+4}}})\\
\vdots\\
\frac{1}{2}\nabla_d\nabla_d f|_{z}+O(\varepsilon^{\min(\tau,\kappa)})+O(\sqrt{\frac{\log(n)}{n \varepsilon^{d+4}}})\\
\\
\nabla_1\nabla_2 f|_{z}
+O(\varepsilon^{\min(\tau,\kappa)})+O(\sqrt{\frac{\log(n)}{n \varepsilon^{d+4}}})\\
\vdots\\
\nabla_{d-1}\nabla_d f|_{z}
+O(\varepsilon^{\min(\tau,\kappa)})+O(\sqrt{\frac{\log(n)}{n \varepsilon^{d+4}}})
\end{array}\right]
$$ 
for some normal coordinates $\{x^j\}, j=1,\dots,d,$ around $z$, where $\tau=\frac{1}{2}$ when $z\in M_{\sigma}$ and $\tau=1$ when $z\in M\setminus M_{\sigma}$, and the implied constants of the big O terms depend on $\|f\|_{C^{2,\kappa}}$, $\|\rho_X\|_{C^2}$, and the dimension $d$ and the curvature of $M$ at $z$.  
\end{Thm}

The proof of Theorem \ref{main} is deferred to Appendix \ref{section: proof of main theorem heuristic approach}.
This theorem establishes that the widely used Hessian estimator based on local PCA yields an accurate approximation of the true Hessian operator for $C^{2,\kappa}$ functions, even in the presence of nonuniform sampling, nontrivial curvature, and non-empty boundaries.
Moreover, the approximation error is uniform over $M$, as the implicit constants in the big-O notation can be uniformly bounded due to the compactness and smoothness assumptions. It is evident that the convergence rate improves as the point $z$ moves away from the boundary.
An immediate corollary is that as $n \to \infty$, we can recover the Hessian operator at all points ${x_i}_{i=1}^n$ almost surely by applying a direct union bound; i.e., Boole's inequality, along with the Borel–Cantelli lemma.

We now compare our result with existing related work. When $M$ is a flat manifold, our result recovers that of \cite{DG03}, and when $M$ is Euclidean space, it aligns with the classical analysis of local quadratic regression in \cite{FG96}.
Our estimator is closely related to the quadratic regression estimator with a 0–1 kernel. For instance, the error bound for the Hessian estimator under quadratic regression with an appropriately chosen kernel is uniform over $[a,b] \subset \mathbb{R}$ and given by $O_P\left(\varepsilon+\sqrt{\frac{\log(n)}{n \varepsilon^{5}}}\right)$, as shown in \cite[Theorem 6.5]{FY03}. This matches our result in the special case of $d=1$, $f \in C^{2,1}$, and when the evaluation point $z$ is sufficiently far from the boundary (so that $\tau = 1$). The primary discrepancy arises when $z$ is near the boundary. This discrepancy stems from the tangent space estimation, where $\tau = \frac{1}{2}$ reflects the boundary-induced error characterized in \cite[Theorem B.1]{SW12}. In essence, when estimating the tangent space near the boundary using local PCA, the asymmetry of the $\epsilon$-neighborhood prevents cancellation of the first-order term in the Taylor expansion, resulting in a larger bias compared to interior points. In other words, if we can design a tangent space estimator with faster convergence near the boundary (i.e., achieving $\tau \geq 1$), or if the tangent space is known a priori, the error near the boundary could be reduced to match the interior rate.
It is important to note, however, that even if the tangent space is exactly known, and thus eliminating errors from its estimation, the convergence rate cannot be further improved, since the dominant bias arises from the manifold curvature itself.

We give a sketch of proof here and leave the detailed proof in the Appendix \ref{section: proof of main theorem heuristic approach}. 
Recall that the first step of evaluating $( Z^T Z)^{-1} Z^T{\bf f}$ is finding the top principal vectors $\{u_l\}_{l=1}^d$ at $z$ by using the local PCA. Then we construct the base matrix $Z$ by using $\{u_l\}_{l=1}^d$. 
Intuitively, once we control the error in the tangent space estimation, we may {\em directly} evaluate the bias and variance terms between $\frac{1}{k_z} Z^T Z$ and its continuous counterpart, as well as those between $\frac{1}{k_z} Z^T{\bf f}$ and its continuous counterpart (see Lemma \ref{ZTf} for an example). Subsequently, we multiply $\left(Z^T Z\right)^{-1}$ and $Z^T{\bf f}$. However, this {\em direct expansion approach} faces challenges due to the heterogeneous asymptotic orders in the entries of $Z^T Z$ and $Z^T{\bf f}$ (see \eqref{Expansion of Z^TZ for the direct approach} for an expansion of $Z^TZ$ and \eqref{Expansion of Z^Tf for the direct approach} for an expansion of $Z^T{\bf f}$). Upon inversion, heterogeneous asymptotic orders in $Z^T Z$ give rise to several new terms (see Section \ref{subsection inversion of Z^TZ in the direct} for details) that complicate the multiplication step. To address all terms with heterogeneous asymptotic orders, we require higher-order local approximations of the function, density function, and curvature, further complicating the analysis.
In Section \ref{Appendix proof direct approach}, we demonstrate the limitation of this approach by showing that with extensive computation, completing the proof for interior points $z \in M \setminus M_{\sigma}$ is possible. For readers with interest in a glance of the challenge, see a detailed summary of this approach in the beginning of Section \ref{Appendix proof direct approach}. More technical explanations can be found in Remarks \ref{remark of comparison of two methods, 1}, \ref{remark of comparison of two methods, 2} and \ref{remark of comparison of two methods, 3}. The complexity of analyzing interior points highlights the significantly greater challenges of addressing boundary points.

To handle the above-mentioned challenge, in this work we propose to use a {\em reduction trick}. With this reduction trick, we can bypass the intensive computation and the difficulty one may encounter in the direct expansion approach involving taking inverse. The key observation underlying this reduction trick is that the leading order terms of 
$\frac{1}{k_z}  Z^T{\bf f}$ are the same as 
$\frac{1}{k_z}  Z^T  Z
\left[\begin{array}{c}
f(z)\\
\nabla_j f|_z\\
\nabla_i\nabla_jf|_z
\end{array}\right]$.  
With this observation, we have 
\begin{align}
\left(\frac{1}{k_z}  Z^T  Z\right)^{-1}\frac{1}{k_z}  Z^T{\bf f}
&\ =\left(\frac{1}{k_z}  Z^T  Z\right)^{-1}\left(
\frac{1}{k_z}  Z^T  Z\left[\begin{array}{c}
f(z)\\
\nabla_j f|_z\\
\nabla_i\nabla_jf|_z
\end{array}\right]
+V\right)\nonumber\\
&\ =\left[\begin{array}{c}
f(z)\\
\nabla_j f|_z\\
\nabla_i\nabla_jf|_z
\end{array}\right] + \left(\frac{1}{k_z}  Z^T  Z\right)^{-1}V,\label{reduction trick main equation}
\end{align}
where $V$ consists of higher order terms compared with $\frac{1}{k_z}  Z^T  Z
\left[\begin{array}{c}
f(z)\\
\nabla_j f|_z\\
\nabla_i\nabla_jf|_z
\end{array}\right]$.
This observation makes the computation much lighter for any $z$, either away from or near $\partial M$. Indeed, since $V$ contains higher order terms, we only need to track fewer higher order terms in $\left(\frac{1}{k_z}  Z^T  Z\right)^{-1}$ compared with the direct expansion approach. As a result, $\left(\frac{1}{k_z}  Z^T  Z\right)^{-1}V$ can be more easily calculated, and hence the bias and variability of our Hessian estimator compared to the Hessian operator. Indeed, since $V$ represents a higher-order term, we do not require as high-order an expansion of $\frac{1}{k_z} Z^T Z$ as needed in the direct expansion approach. This dramatically simplifies the proof.

It is noteworthy that the biases of the estimates for the function and its gradient are of a higher order compared to the bias of the Hessian estimate. This discrepancy is expected given that the Hessian represents a higher-order term. We hypothesize that this bias estimate can be improved with a more accurate tangent space estimation. In our theorem, the bias order of $\varepsilon$ is $\frac{1}{2}$ higher for interior points than for points near the boundary. This difference arises precisely due to the tangent space estimation derived in \cite{SW12}.

\section{Discussion of relevant topics}

With the result of estimating the Hessian on the manifold setup, one immediate interesting problem to study is the squared Hessian energy used in the HLLE algorithm \cite{DG03} and computer graphics \cite{SGWJ18,SJWG20,HWAG09}. 
Recall that the authors in HLLE \cite{DG03} define an embedding of the dataset by minimizing the squared Hessian energy 
\begin{equation}\label{equation mathcal H}
\mathscr{H}(f):=\int_M|Hess(f)|^2d{\rm vol}.
\end{equation}
They implement HLLE by using the fact that the kernel space of $\mathscr{H}$ consists of linear functions under the assumption that the manifold $M$ is isometric to a connected open subset in $\mathbb{R}^d$. However, in general, the kernel space of $\mathscr{H}$ is more complicated and the behavior of the Hessian energy $\mathscr{H}(f)$ when $M$ is a general manifold with non-trivial geometry and topology is unknown.
It is worth mentioning that a different interpretation about Hessian Eigenmap and its validity in view of locally linear embedding can be found in \cite{LC21}.

We shall elaborate the mathematical challenge when we study HLLE. In particular, for a given manifold, can we characterize the functions which minimize the Hessian energy $\mathscr{H}(f)$? 
On a closed $d$-dimensional smooth Riemannian manifold, we can calculate the first variation of Hessian energy $\mathscr{H}(f)$ and derive its Euler-Lagrange equation. We remind the reader that Euler-Lagrange equation is a necessary condition for $f$ to be a minimizer.

\begin{Lma}
Let $(M,g)$ be a closed $d$-dimensional smooth Riemannian manifold with Ricci curvature $Ric$ and scalar curvature $S$. Consider $\mathscr{H}(f):=\int_M|Hess(f)|^2d{\rm vol}$, where $f\in C^\infty(M)$. If $f$ is a minimizer of $\mathscr{H}$, then 
\begin{equation}\label{EL}
\Delta^2f+\langle Ric,Hess(f)\rangle+\frac{1}{2}\langle \nabla S,\nabla f\rangle=0.
\end{equation}
\end{Lma}

\begin{proof}
Consider a variation of $f$ via $f+th$, where $t\in(-\varepsilon, \varepsilon)$ and $h:M\to\mathbb{R}$ is a perturbation function. The variation of Hessian energy is
\begin{align*}
\frac{d}{dt}\bigg|_{t=0} \int_M|Hess(f+th)|^2\ d{\rm vol}
=&\ \sum_{i,j=1}^d\int_M  2 \langle \nabla_i\nabla_j h, \nabla_i\nabla_j f\rangle \ d{\rm vol}\\
=&\ \sum_{i,j=1}^d\int_M  2 \langle  h, \nabla_j\nabla_i\nabla_i\nabla_j f\rangle \ d{\rm vol}.
\end{align*}
By using the traced second Bianchi identity, we have
\begin{align*}
\sum_{i,j=1}^d\nabla_j\nabla_i\nabla_i\nabla_j f
=&\ \sum_{i,j=1}^d\nabla_j\left(\nabla_j\nabla_i\nabla_i f
+Ric(e_j,\nabla f)\right)\\
=&\ \Delta^2f
+ {\rm div}Ric(\nabla f,\cdot)\\
=&\ \Delta^2f+\langle Ric,Hess(f)\rangle+\frac{1}{2}\langle \nabla S,\nabla f\rangle,
\end{align*}
where $S$ is the scalar curvature and $Ric$ is the Ricci curvature of $M$.
So, if the minimizing function exists, it must satisfy the Euler-Lagrange equation (\ref{EL}).
\end{proof}

Therefore, while we do not have the spectral convergence result for the discretized $\mathscr{H}$ yet, we could conjecture that the dimension reduction achieved by the HLLE algorithm depends on the eigenstructure of the fourth order differential operator 
\[
\mathcal{L}:=\Delta^2+ Ric\cdot \nabla^2+\frac{1}{2} \nabla S\cdot \nabla. 
\]

Although a fourth order equation probably has many solutions, it is challenging to solve (\ref{EL}) explicitly due to the curvature-involving coefficients. 
Note that when $M$ is a Euclidean space or a flat manifold, $\mathcal{L}$ is simply the bi-Laplacian, and its kernel includes the span of the constant and linear functions. This partially explains how HLLE functions under the setup in \cite{DG03}. For a generic Riemannian manifold $M$, $\mathcal{L}$ still has good properties. However, we cannot expect $\mathcal{L}$ has a nontrivial kernel.

\begin{Thm} On a closed connected Riemannian manifold $(M,g)$ with Ricci curvature $Ric$ and scalar curvature $S$, the spectrum of \[
\mathcal{L}:=\Delta^2+ Ric\cdot \nabla^2+\frac{1}{2} \nabla S\cdot \nabla. 
\]
is discrete and diverges to infinity. Moreover, when the manifold is Einstein, that is, when $Ric=\Lambda g$, for some $\Lambda>0$, the zero eigenvalue of $\mathcal{L}$ is simple.
\end{Thm}

\begin{proof}
In general, although the order of $\mathcal{L}$ is high, it is  linear and self-adjoint. Note that
every eigenvalue $\lambda$ of $\mathcal{L}$ is nonnegative because 
$$\lambda\|\phi\|^2
=\int_M(\mathcal{L}\phi)\phi d{\rm vol}
= \sum_{i,j=1}^d\int_M (\nabla_j\nabla_i\nabla_i\nabla_j \phi)
\varphi d{\rm vol}
= \int_M\left|\nabla^2\phi\right|^2d{\rm vol}
\geq 0$$
provided that $\phi$ is an eigenfunction corresponding to $\lambda$.
Moreover, for any $\mu>0$,
\begin{align}\label{inner product}
\langle\!\langle u,v \rangle\!\rangle
:=\int_M \langle \nabla^2u,\nabla^2v\rangle +\mu uv\ d{\rm vol}
=\int_M \Delta u\Delta v - Ric(\nabla u,\nabla v)+\mu uv\  d{\rm vol}
\end{align}
defines an inner product on the Sobolev space $H^{2}(M)$. 
Given $\varphi\in L^2(M)$. Recall that a weak solution of the equation $(\mathcal{L}+\mu I)f=\varphi$ is a function $f\in H^2(M)$ that satisfies 
\[
\langle\!\langle f, v \rangle\!\rangle
= \int_M \langle \nabla^2 f,\nabla^2v\rangle +\mu f v\ d{\rm vol}
=\int_M \varphi v \ d{\rm vol}
\]
for all $v \in C_c^\infty(M)$. Define a functional $F_{\varphi}:H^2(M)\to \mathbb{R}$ by $F_{\varphi}(v):=\int_M \varphi v \ d{\rm vol}$. Since $F_{\varphi}$ is bounded, using Riesz representation theorem, we know that
there exists a unique $f\in H^2(M)$ such that $F_{\varphi}(v)=\langle\!\langle f,v\rangle\!\rangle$ for all $v\in H^2(M)$. Such $f$ is a weak solution to the equation $(\mathcal{L}+\mu I)f=\varphi$. Thus the inverse operator $(\mathcal{L}+\mu I)^{-1}$, which sends $\varphi$ to $f$, is well-defined. Analogue to the classical theory of Laplacian, the operator $(\mathcal{L}+\mu I)^{-1}$ is compact because the embedding $H^2(M) \hookrightarrow L^2(M)$ is compact.
The compactness implies that $\mathcal{L}$ has a discrete spectrum which diverges to infinity. 
So we obtain the first statement of this theorem.

When $Ric=\Lambda g$ with $\Lambda>0$, equation (\ref{EL}) becomes 
$\Delta (\Delta f+\Lambda f)=0$. Note that such manifolds must be compact and there are no non-constant harmonic functions on them. 
From Lichnerowicz theorem (cf. \cite[Theorem 5.1]{Li12}), we know that $-\Lambda$ is not an eigenvalue of $\Delta$, so the kernel of $\Delta +\Lambda$ is trivial. Hence, using Fredholm alternative, $\Delta f+\Lambda f=F$ has a unique solution $f\in H^{1}(M)$ for any given non-zero function $F\in L^2(M)$. Note that such solution $f$ must be smooth by the elliptic regularity. Thus $\Delta (\Delta f+\Lambda f)=\Delta F=0$ whenever $F$ is chosen to be harmonic. Since every harmonic function $F$ must be a constant, it is easy to see that the unique solution of $\Delta f+\Lambda f=F$ 
is the constant $f=\Lambda^{-1}F$. Therefore, $\mathcal{L}f=0$ only has constant solutions, i.e., the eigenvalue $0$ of $\mathcal{L}$ is simple.
\end{proof}

Note that this is an explicit example when the behavior of HLLE cannot be explained, since the behavior of eigenfunctions associated with the smallest eigenvalues is unclear at this moment. 
As a result, in general, the claimed behavior of HLLE may fail if we have a compact manifold without boundary, and we need more understanding of the general $\mathcal{L}$ before making a conclusion. 
For manifolds with $\Lambda<0$, ${\rm ker}(\Delta + \Lambda)$ is again trivial because $-\Delta$ is non-negatively definite. However, such manifold may be non-compact and any non-constant harmonic function $F$ gives us a non-trivial unique solution of equation (\ref{EL}).


We shall mention that if our data is modeled locally by a Ricci flat manifold, i.e., $\Lambda=0$, then the Euler-Lagrange equation becomes the well-known {\it biharmonic equation} $\Delta^2 f=0$. Harmonic maps and biharmonic maps arise naturally because they are critical points of the Dirichlet energy $E(u):=\int_M |du|^2d{\rm vol}$ and bienergy $E_2(u)=\int_M |\Delta u|^2d{\rm vol}$, respectively. They exhibit strong regularity; for instance, continuous biharmonic maps must be smooth (cf. \cite{CWY99a,CWY99b}). Related results about biharmonic maps can be found in \cite{GK90, H06, MO06}. Numerous findings about harmonic and biharmonic maps, among others due to our limited survey coverage, have been applied to different branches of mathematics \cite{T05, H06}, numerical computation \cite{D91, LTZ01}, theoretical physics \cite{N11,S87}, and geometric processing \cite{SGWJ18, FN19}. Note that the bi-Laplacian also appears in the traditional locally linear embedding (LLE) algorithm \cite{WW23}. To analyze the performance of these algorithms, it is critical to understand the behavior of the bi-Laplace operator and biharmonic maps, and more generally the $\mathcal{L}$ operator, bridging the gap between discrete operators and the continuous theory. For example, we need the associated spectral theory and spectral convergence results to fully understand LLE and HLLE. For fourth-order operators like $\mathcal{L}$, there have been some results for Paneitz-type operators, which are strongly related to conformal geometry and $Q$-curvature. See for instance \cite{CY97, DHL00}. Last but not least, higher-order differential operators have been studied for their usefulness in many practical problems; see, for example \cite{Mu94,Mu02}, among others. We look forward to discovering more geometric properties of the operator $\mathcal{L}$ and more general differential operator with their practical applications.

Last but not the least, while we do not explore the minimax rate of the Hessian operator estimate, it is an interesting future direction to explore it following the work in \cite{AL19}. 

\section{Conclusion}

We have presented a systematic analysis of the commonly used Hessian operator estimator in manifold settings. Our findings demonstrate that this estimator asymptotically converges to the target Hessian operator, even near boundaries, with the effects of nonuniform sampling and curvature proving asymptotically negligible. Our proposed analytical framework notably streamlines the rigorous theoretical calculations required to achieve this objective.

\section*{Acknowledgment}
The authors would like to acknowledge the anonymous reviewers for their constructive feedback and the suggestion to compare the minimax rate estimation. C.-W. Chen would like to thank National Center for Theoretical Sciences (NCTS) for its constant support. He is supported by Young Scholars Program and Short-Term Overseas Research Program of NSTC.

\setcounter{page}{1}
	\setcounter{equation}{0}
	\setcounter{figure}{0}
	\renewcommand{\thepage}{S.\arabic{page}}
	\renewcommand{\thesection}{S.\arabic{section}}
	\renewcommand{\theequation}{S.\arabic{equation}}
	\renewcommand{\thetable}{S.\arabic{table}}
	\renewcommand{\thefigure}{S.\arabic{figure}}

\appendix


\section{A summary of the Hessian operator}\label{section:hessian review}

We briefly recall the definition of gradient and Hessian under the manifold setup and the associated materials, and for readers with interest in a systematic treatment of the Riemannian geometry framework, we refer them to \cite{Lee12} for details.

Recall that a $d$-dimensional $C^k$-manifold $M$ is a collection of open domains with {\it local} charts, whose overlapping regions satisfy the $C^k$ transition condition, and each domain is $C^k$-diffeomorphic to an open ball in $\mathbb{R}^d$ via the chart function. The language of manifold is designed for doing local computation in each domain and integrating all quantities by gluing up all domains through overlapping regions.

Consider a point $\xi\in U\subset M$ and a chart (or coordinate) function $\varphi: U\to \varphi(U)\subset\mathbb{R}^d$. For a $C^2$-function $f:M\to\mathbb{R}$, the partial derivative of $f$ at $\xi$ is defined as the derivative of $f\circ\varphi$ at $\varphi(\xi)$. Specifically, if we denote $\varphi(\xi)=(x^1(\xi),\dots,x^d(\xi))$ as the coordinate function, we have $\frac{\partial f}{\partial x^k}(\xi):=\frac{\partial (f\circ \varphi)}{\partial x^k}(\varphi(\xi))$. Obviously these derivatives depend on the choice of coordinate functions. In order to define coordinate-independent derivatives, we need the concepts of {\it tensor} and {\it covariant derivative}.
Given a tangent vector $v=v^k\frac{\partial}{\partial x^k}$ in $T_{\xi}M$, the directional derivative of $f$ with respect to $v$ at $\xi$ is defined by $vf:=v^k\frac{\partial f}{\partial x^k}$. Here and below we adopt Einstein's convention which means that repeated indices are summed over even if there is no sum symbol. The differential of $f$ is defined as the 1-form $df$ which satisfies $df(v):=vf$ for all $v\in T_{\xi}M$. Moreover, the Hessian of $f$ is defined to be the covariant derivative of $df$, $Hess(f):=Ddf$, where $D$ is the covariant derivative based on the Levi-Civita connection associated with the given metric. $Hess(f)$ is an invariant 2-tensor which, using the local coordinate $\{x^j\}_{j=1}^d$, can be computed as 
$$Ddf=D(f_jdx^j) = f_{ij}dx^i\otimes dx^j-f_j\Gamma_{ik}^jdx^i\otimes dx^k =(f_{ij}-f_k\Gamma_{ij}^k)dx^i\otimes dx^j,$$ 
where $f_j:=\frac{\partial (f\circ\varphi)}{\partial x^j}$, $f_{ij}:=\frac{\partial^2 (f\circ\varphi)}{\partial x^i\partial x^j}$, and $\Gamma_{ij}^k$ are the Christoffel symbols of the Levi-Civita connection. 
To simplify the Hessian estimator on a manifold, it is desirable to get rid of the unknown $\Gamma_{ij}^k$. This can be achieved by using normal coordinates, i.e., $\varphi$ is chosen to be the inverse of the exponential map, denoted as $\exp^{-1}:U\to T_{\xi}M$. In this case,  $x^j$'s are called normal coordinates and all $\Gamma_{ij}^k$ vanish at $\xi$. Hence the Hessian operator is expressed as
\[
Hess(f)=f_{ij}dx^i\otimes dx^j=\frac{\partial^2 f}{\partial x^i\partial x^j}dx^i\otimes dx^j
\] 
at $\xi$ under the normal coordinates. In particular, when we compute the $Hess(f)$ at $\xi\in U\subset M$, the Christoffel symbol does not play a role under the normal coordinates centered at $\xi$, but it is possible that  $\Gamma_{ij}^k(\zeta)\neq 0$ for some $\zeta\in U$ that is different from $\xi$ under the same normal coordinates.
Note that to simplify the notation, usually researchers use $\frac{\partial^2 f}{\partial x^i\partial x^j}(\xi)$ as the shorthand for $\frac{\partial^2 (f\circ\exp)}{\partial x^i\partial x^j}$ at $\exp^{-1}(\xi)\in T_{\xi}M$. 
Recall that a tensor is invariant under coordinate change, namely, two tensors are identically the same if their coefficients in different coordinate systems obey the transformation law.

\section{Technical lemmas}
\subsection{Taylor expansion in local coordinates}
\label{coordinates}
Let $M$ be a smooth complete connected Riemannian manifold. In what follows, $M$ could be closed without boundary, or it has a smooth boundary $\partial M$ and $M\cup\partial M$ is compact. We embed the manifold into $\mathbb{R}^p$ via the inclusion map $\iota$. To analyze the local behavior of our Hessian operator at a point $z\in M$, it is necessary to differentiate between two distinct situations: 
$z$ is close to or away from the boundary. So we define
$M_\sigma:=\{x\in M\ |\ dist(x,\partial M)<\sigma \}$ as the $\sigma$-neighborhood of the boundary $\partial M$. Let $\sigma'$ be the distance between $M\setminus M_\sigma$ and $\partial M$ in $\mathbb{R}^p$, that is, $\sigma'=\inf\left\{ \|\iota(x)-\iota(y)\|_{\mathbb{R}^p}\ |\ x\in M\setminus M_\sigma, y\in\partial M \right\}$. 
Denote $\mathcal{B}:=\iota(M)\cap B_{\varepsilon}^p(\iota(z))$, where $B_{\varepsilon}^p(\iota(z))$ is a Euclidean ball in $\mathbb{R}^p$ with the center $\iota(z)$ and radius $\varepsilon$.

For any $z\in M\setminus M_\sigma$, choose $\varepsilon<\min\{\sigma', inj(M\setminus M_\sigma)\}$ to ensure that the whole set $\mathcal{B}$ lies in the injective region of the exponential map $\exp_z$, where $inj(M\setminus M_\sigma)$ is the infimum of the injectivity radius of points in $M\setminus M_\sigma$.
Then, for any $x\in \mathcal{B}$, $x=\exp_{z}(t\theta)$ for some $t\in\mathbb{R}$ and $\theta\in S^{d-1}\subset T_{z}M$. Any smooth function $f:M\to \mathbb{R}$ can be estimated by its Taylor polynomials
\begin{align}\label{taylor f}
f(x)=f(z)+f_1(\theta)t+f_2(\theta)t^2+f_3(\theta)t^3+O(t^4),
\end{align}
where $f_1(\theta)=\nabla f|_{z}\cdot\theta$, $f_2(\theta)=\frac{1}{2}\nabla^2f|_{z}(\theta,\theta)$, etc. 
When $f$ is merely a $C^{2,\kappa}$-function, we cannot differentiate $f$ three times and the above expansion fails. 
However, since $\displaystyle\sup_M \max_{|\alpha|=2}\frac{|D^{\alpha}f(x)-D^{\alpha}f(y)|}{dist(x,y)^\kappa}\leq C<\infty$, when $t$ is sufficiently small,
we have
\begin{align}
f(x)&=f(z)+f_1(\theta)t+f_2(\xi)t^2\mbox{ for some }\xi\mbox{ between } x\mbox{ and }z\label{Taylor expansion of C2kappa function}\\
&=f(z)+f_1(\theta)t+f_2(\theta)t^2+(f_2(\xi)-f_2(\theta))t^2\nonumber\\
&=f(z)+f_1(\theta)t+f_2(\theta)t^2+O(t^{2+\kappa})\,,\nonumber
\end{align} 
where in the last asymptotic control we use the fact that 
\[
|f_2(\xi)-f_2(\theta)|\leq C\cdot dist(\exp_z(\xi),\exp_z(\theta))^\kappa
\]
By using this estimate of $C^{2,\kappa}$-function, we can express the error term precisely as $O(t^{2+\kappa})$, instead of $o(t^2)$. Note that $f\in C^3(\overline\Omega)$ implies that $f\in C^{2,1}(\overline\Omega)$ only when $\overline\Omega$ is compact with regular boundary. See \cite[p. 53]{GT01} for an example that this property fails when $\partial\Omega$ has a cusp.

The position of $x$ in $\mathbb{R}^p$ can also be approximated by
\begin{align}\label{taylor x}
\iota \circ\exp_{z}(t\theta) -\iota(z) = K_1(\theta)t+K_2(\theta)t^2+K_3(\theta)t^3+K_4(\theta)t^4+O(t^5),
\end{align}
where $K_1(\theta)=\iota(\theta)$, $K_2(\theta)=\frac{1}{2}{\rm I\!I}_{z}(\theta,\theta)$, $K_3(\theta)=\frac{1}{6}\nabla_\theta {\rm I\!I}_{z}(\theta,\theta)$, etc.  Note that $\iota^{-1}(\mathcal{B})$ is not a geodesic ball in $M$ and the radial segment from $z$ to the boundary $\partial(\iota^{-1}(\mathcal{B}))\subset M$ can be approximated by
$$\tilde{\varepsilon}(\theta)=\varepsilon + H_1(\theta)\varepsilon^3 +H_2(\theta)\varepsilon^4+O(\varepsilon^{5}),$$
where $H_1(\theta)=\frac{1}{24}|{\rm I\!I}_{z}(\theta,\theta)|^2$, $H_2(\theta)=\frac{1}{24}\nabla_\theta {\rm I\!I}_{z}(\theta,\theta)\cdot {\rm I\!I}_{z}(\theta,\theta)$, etc. See \cite[Lemma B.3]{WW23} for a proof of this expansion. At last, the volume form on a local region can be approximated by  
$$d{\rm vol}= \left(t^{d-1}+R_1(\theta)t^{d+1}+R_2(\theta)t^{d+2}+O(t^{d+3})\right)dtdS_\theta,$$
where $R_1(\theta)=-\frac{1}{6}Ric_{z}(\theta,\theta)$, $R_2(\theta)=-\frac{1}{12}\nabla_\theta Ric_{z}(\theta,\theta)$, etc. (cf. Corollary 2.10 in \cite{G73}). 

For a point $z$ in $M_\sigma$, things could be done similarly but in a more nuanced way. The main issue is that $\mathcal{B}$ might not be covered by the range of $\exp_z$ due to the effect of boundary. There might be some point $x\in\mathcal{B}$ which cannot be 
approximated as in formula (\ref{taylor x}). This issue can be resolved by considering a Riemannian extension $N$ of $M$ (cf. \cite{PV20}). By extending the manifold a bit outwards from $\partial M$ and considering the exponential map of $N$ at $z$, $\exp_Z^N$, one sees that every $x$ in $\mathcal{B}$ can be expressed as $\exp_z^N(t\theta)$ and all the above formulas hold similarly. Moreover, we require $\varepsilon$ to be sufficiently small so that $\mathcal{B}$ is close to the truncated Euclidean ball $B_{\varepsilon,\delta}^d:=B_\varepsilon(0)\cap\{(x^1,\dots,x^d)\in\mathbb{R}^d\ |\ x^{d}>-(1-\delta)\varepsilon\}$ in $\mathbb{R}^d$ for some $\delta\in[0,1)$ in the following sense. 

\begin{Def}
Let $B_{\varepsilon,\delta}^d:=B_\varepsilon(0)\cap\{(x^1,\dots,x^d)\in\mathbb{R}^d\ |\ x^{d}>-(1-\delta)\varepsilon\}$ in $\mathbb{R}^d$ for some $\delta\in[0,1)$ and $e_k:=\frac{\partial}{\partial x^k}$. Consider the polar coordinates $t=(\sum_{k=1}^d(x^k)^2)^{\frac{1}{2}}$ and $\theta\in S^{d-1}$.
We say that $\mathcal{B}$ is close to $B_{\varepsilon,\delta}^d$ if

{\rm (i)} 
$(\exp_z^N)^{-1}(\mathcal{B})\subset T_zN$ can be parameterized by $(t,\theta)$ with $t\in[0,\widetilde{r}(\theta))$, and
$\exp_z^N(\widetilde{r}\theta)=\partial\mathcal{B}$; 

{\rm (ii)} $\partial\mathcal{B}$ consists of two parts, the truncated portion $\iota(\partial M)\cap B_{\varepsilon}^p(0)$ and the spherical portion $\iota(M)\cap \partial B_{\varepsilon}^p(0)$, and
$$\widetilde{r}(\theta)=\bigg\{
\begin{array}{ll}
\varepsilon+O(\varepsilon^3),& \forall  \theta \mbox{ s.t. }\exp_z^N(\widetilde{r}\theta)\in \iota(M)\cap \partial B_{\varepsilon}^p(0)\\
r+O(r^3),& \forall \theta \mbox{ s.t. }\exp_z^N(\widetilde{r}\theta)\in
 \iota(\partial M)\cap B_{\varepsilon}^p(0)
\end{array}, 
$$
where $r=r(\theta)=\varepsilon(1-\delta)\langle\theta,-e_d\rangle^{-1}$.
\end{Def}

\begin{figure}\vspace{-2cm}
    \includegraphics{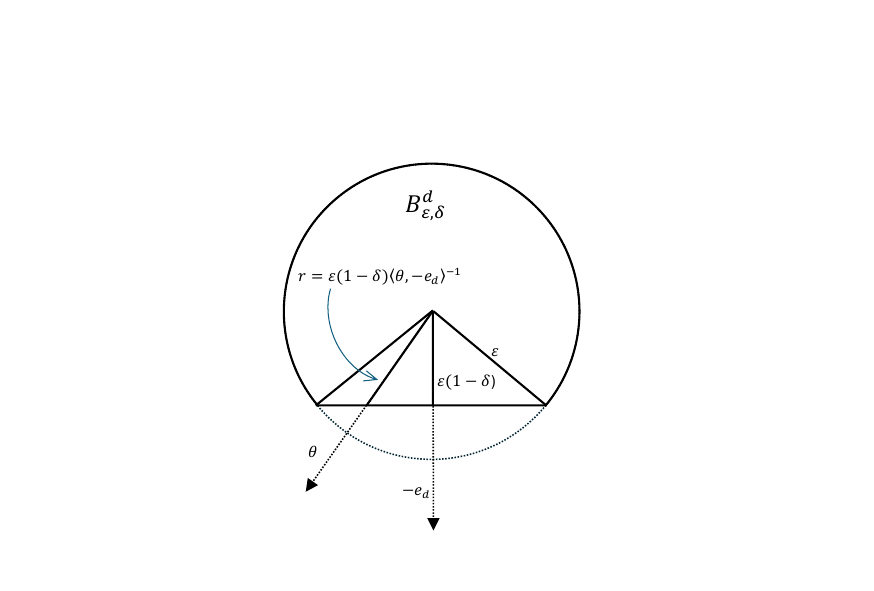}
    \vspace{-1cm}
    \caption{$B_{\varepsilon,\delta}^d$ and $r(\theta)=\varepsilon(1-\delta)\langle\theta,-e_d\rangle^{-1}$\hspace{2.5cm}\,}
    \label{fig:B}
\end{figure}

The choice of such $\varepsilon$ would depend on the magnitude of the second fundamental forms of $\partial M\subset N$ and of $\iota(M)\subset\mathbb{R}^p$, and thus could be chosen uniformly on $M\cup\partial M$.  

The smoothness assumption imposed in \cite{PV20} is the main reason we impose the smooth manifold assumption. If the manifold is boundary free, then a $C^3$ manifold assumption is sufficient to obtain our main theorem. If the result in \cite{PV20} can be 
generalized to manifolds with $C^3$ boundary,
our result can be extended as well. However, this is out of the scope of this paper.

\subsection{Integrals needed for the main proof}
In this subsection, we compute some integrals on domains in $T_zN \simeq\mathbb{R}^d$. Let the space be equipped with Cartesian coordinates $\{x^i\}$ and standard basis $\{e_i=\frac{\partial}{\partial x^i}\}$.
Moreover, we denote $B_{\varepsilon,\delta}$ as the truncated Euclidean ball in $T_zN$ and $\mathcal{B}$ as the preimage of $\mathcal{B}$ via $\exp_z^N$.
Let $(t,\theta)$ be the polar coordinates, $t\in\mathbb{R}_{\geq 0}$, $\theta\in S^{d-1}$. 
Denote $\bar\theta^i:=\langle \theta, e_i \rangle$. 
We use 
Einstein's convention by summing over all equally named indices and thus denote $\theta=\bar\theta^ie_i$.

For all non-negative integers $m$ and $k$, we denote
\begin{align*}
C_{m,2k}:=
&\ 
\varepsilon^{-(d+m+2k)}\int_{B_{\varepsilon,\delta}} t^{m+2k} (\bar\theta^{d})^m(\bar\theta^{d-1})^{2k} d{\rm vol}.
\end{align*}
Note that $C_{m,2k}$, being functions of $m,k,d,\delta$, are independent of $\varepsilon$.

\begin{Def} We introduce the following Greeks to simplify the expressions afterwards. 
$$\gamma_1=\varepsilon\frac{C_{1,0}}{C_{0,0}};\ \ \alpha_1=\varepsilon^2\frac{C_{2,0}}{C_{0,0}},\ \alpha_2=\varepsilon^2\frac{C_{0,2}}{C_{0,0}};$$

$$\mu_1=\varepsilon^3\frac{C_{3,0}}{C_{0,0}},\ \mu_2=\varepsilon^3\frac{C_{1,2}}{C_{0,0}};$$

$$\beta_1=\varepsilon^4\frac{C_{4,0}}{C_{0,0}},\ 
\beta_2=\varepsilon^4\frac{C_{2,2}}{C_{0,0}},\  
\beta_3=\varepsilon^4\frac{C_{0,4}}{C_{0,0}},\
\beta_4=\varepsilon^4\frac{C_{0,2,2}}{C_{0,0}}.$$
\end{Def}

It is not hard to check that, when $\delta=0$, these Greeks become $\gamma_i=\mu_i=0$ for all $i$; $\alpha_1=\alpha_2=\frac{1}{d+2}\varepsilon^2$; $\beta_1=\beta_3=\frac{3}{(d+2)(d+4)}\varepsilon^4$
; $\beta_2=\beta_4=\frac{1}{(d+2)(d+4)}\varepsilon^4$.

\begin{Lma} \label{volume}Suppose that there is a weighted function $\rho(x)\in C^2(M)$ on $M$, e.g. the probability density function (p.d.f.) associated with a random vector.
Let $\mathcal{B}:=\iota(M)\cap B_{\varepsilon}^p(\iota(z))$ be close to $B_{\varepsilon,\delta}$. Then
\begin{align*}
{\rm Vol_{\rho}}(\mathcal{B})
=&\ \varepsilon^d C_{0,0}\rho(z)+\varepsilon^{d+1}C_{1,0}\partial_d\rho(z)
+O(\varepsilon^{d+2})\\
=&\ \varepsilon^d C_{0,0}\rho(z)\left(1+\gamma_1\frac{\partial_d\rho(z)}{\rho(z)}+O(\varepsilon^{2})\right).
\end{align*}

\end{Lma}

\begin{proof}
We use geodesic polar coordinates to compute the local volume expansion. This computation is standard in Riemannian geometry (see for instance A. Gray's article \cite{G73} or his book \cite[Ch.9]{G04}).
{\allowdisplaybreaks
\begin{align*}
&\ {\rm Vol}_\rho(\mathcal{B})\\
=&\ \int_{S^{d-1}}\int_0^{\tilde{r}(\theta)}(\rho(z)+\rho_1(\theta)t+\rho_2(\theta)t^2+O(t^3))(t^{d-1}+R_1(\theta)t^{d+1}+O(t^{d+2}))dtdS_\theta\\
=&\ \int_{S^{d-1}}\int_0^{\tilde{r}(\theta)}\rho(z)t^{d-1}+\rho_1(\theta)t^{d}+(\rho_2(\theta)+\rho(z)R_1(\theta))t^{d+1}+O(t^{d+2}) dtdS_\theta\\
=&\ \int_{S^{d-1}}\int_0^{\tilde{r}(\theta)}\rho(z)t^{d-1}+\sum_{i=1}^d\partial_i\rho(z)\bar\theta^it^{d}dtdS_\theta\\
&\ +
\int_{S^{d-1}}\int_0^{\tilde{r}(\theta)}
\left(\frac{1}{2}\nabla_\theta^2\rho|_z-\frac{1}{6}\rho(z)Ric_z(\theta,\theta)\right)t^{d+1}+O(t^{d+2}) dtdS_\theta\\
=&\ \rho(z)\int_{S^{d-1}}\int_0^{\tilde{r}(\theta)}t^{d-1}
dtdS_\theta
+\sum_{i=1}^d\partial_i\rho(z)\int_{S^{d-1}}\int_0^{\tilde{r}(\theta)}\bar\theta^it^d dtdS_\theta\\
&\ +
\sum_{i,j}\left(\frac{1}{2}\partial_i\partial_j\rho|_z-\frac{1}{6}\rho(z)R_{ij}|_z\right)
\int_{S^{d-1}}\int_0^{\tilde{r}(\theta)} \bar\theta^i\bar\theta^j t^{d+1} dtdS_\theta+O(\varepsilon^{d+3})\\
=&\ \rho(z)C_{0,0}
\varepsilon^d+\frac{d}{d+2}\rho(z)\int_{S^{d-1}} \frac{1}{24}|{\rm I\!I}_z(\theta,\theta)|^2r(\theta)^{d+2}dS_{\theta}
+\partial_d\rho(z)\varepsilon^{d+1}C_{1,0}\\
&\ +
\left(\left(\frac{1}{2}\Delta\rho(z)-\frac{1}{6}\rho(z)S\right)
C_{0,2}+\left(\frac{1}{2}\partial_d^2\rho|_z-\frac{1}{6}\rho(z)R_{dd}|_z\right)(C_{2,0}-C_{0,2})\right)\varepsilon^{d+2}+O(\varepsilon^{d+3})\\
=&\ \rho(z)C_{0,0}\varepsilon^d+\partial_d\rho C_{1,0}\varepsilon^{d+1}+\left(\rho(z)\Lambda_{\delta} +\frac{1}{2}C_{0,2}\Delta\rho(z)+\frac{1}{2}(C_{2,0}-C_{0,2})\partial_d^2\rho|_z
\right)\varepsilon^{d+2}+O(\varepsilon^{d+3}),
\end{align*}
}
where $\Lambda_{\delta}$ is a combination of the second fundamental form and Ricci curvature at $z$.
Note that in the equation above we use $\rho\in C^{3}$ to conclude that the error is in $O(\varepsilon^{d+3})$. Since the statement of the lemma only concerns terms of order up to $\varepsilon^{d+2}$, $\rho\in C^2$ is sufficient to us. 
\end{proof}

\begin{Rmk}
By Gauss-Codazzi equation, a totally geodesic submanifold in $\mathbb{R}^p$ must be Ricci flat and thus $\Lambda_\delta=0$. On the other hand, if $M$ is merely intrinsically flat, then the second fundamental form may not be flat and $\Lambda$ may not be zero. For instance, if the data distributed on a 2-dimensional cylinder, then we have to consider the effect of extrinsic curvatures even though $M$ itself is intrinsically flat. 
\end{Rmk}

\begin{Lma}\label{lemma of local poly integration}
Let $\tau=\Big\{
\begin{array}{rl}
    \frac{1}{2}, &\mbox{ if }z\in M_\sigma \\
    1, &\mbox{ if }z\in M\setminus M_\sigma
\end{array}$.
Denote $\partial_k\rho=\langle \nabla\rho,e_k\rangle$ and $\{u_l\}_{l=1}^d$ as the basis obtained by local PCA. We have

{\rm (i)}
For all $k=1,\dots,d$, 
$$-\hspace{-3.8mm}\int_{\mathcal{B}}\langle x-z,u_k\rangle   d{\rm vol}_\rho
=\bigg\{
\begin{array}{rc}
 \gamma_1
+(\alpha_1-\gamma_1^2)\frac{\partial_d\rho}{\rho}\big|_{z}+O(\varepsilon^{\tau+2}),    & k=d\\ 
 \alpha_2\frac{\partial_k\rho}{\rho}\big|_{z}+O(\varepsilon^{\tau+3}),& k\neq d 
\end{array}.
$$

{\rm (ii)} For all $k=1,\dots,d$ and $j=1,\dots, d-1$, 
$$-\hspace{-3.8mm}\int_{\mathcal{B}}\langle x-z,u_k\rangle^2   d{\rm vol}_\rho
=\bigg\{
\begin{array}{rc}
 \alpha_1+(2\mu_1-\gamma_1\alpha_1)\frac{\partial_d\rho}{\rho}\big|_{z}+O(\varepsilon^{\tau+3}),    & k=d\\ 
 \alpha_2+(2\mu_2-\gamma_1\alpha_2)\frac{\partial_d\rho}{\rho}\big|_{z}+O(\varepsilon^{\tau+3}),& k\neq d 
\end{array},$$

$$-\hspace{-3.8mm}\int_{\mathcal{B}}\langle x-z,u_k\rangle\langle x-z,u_j\rangle   d{\rm vol}_\rho
=\bigg\{
\begin{array}{rc}
 \mu_2\frac{\partial_j\rho}{\rho}\big|_{z}+O(\varepsilon^{\tau+3}),    & k=d\\ 
 O(\varepsilon^{\tau+3}),& j\neq k\neq d 
\end{array}.$$

{\rm (iii)} 
For all $k=1,\dots,d$ and $j,l,m=1,\dots, d-1$, 
$$-\hspace{-3.8mm}\int_{\mathcal{B}}\langle x-z,u_k\rangle^2\langle x-z,u_d\rangle   d{\rm vol}_\rho
=\bigg\{
\begin{array}{rc}
 \mu_1+\beta_1\frac{\partial_d\rho}{\rho}\big|_{z}+O(\varepsilon^{\tau+4}),    & k=d\\ 
 \mu_2+\beta_2\frac{\partial_d\rho}{\rho}\big|_{z}+O(\varepsilon^{\tau+4}),& k\neq d, 
\end{array},$$
$$-\hspace{-3.8mm}\int_{\mathcal{B}}\langle x-z,u_k\rangle^2\langle x-z,u_j\rangle   d{\rm vol}_\rho
=\bigg\{
\begin{array}{rc}
 \beta_2\frac{\partial_j\rho}{\rho}\big|_{z}+O(\varepsilon^{\tau+4}),    & k=d\\ 
  \beta_4\frac{\partial_j\rho}{\rho}\big|_{z}+O(\varepsilon^{\tau+4}),& j\neq k\neq d
\end{array};$$
otherwise,
$$-\hspace{-3.8mm}\int_{\mathcal{B}}\langle x-z,u_j\rangle\langle x-z,u_l\rangle\langle x-z,u_m\rangle   d{\rm vol}_\rho
= O(\varepsilon^{\tau+4}).$$

{\rm (iv)} 
For all $j,l=1,\dots, d-1$ and $k,m,n=1,\dots,d$, 
$$-\hspace{-3.8mm}\int_{\mathcal{B}}\langle x-z,u_k\rangle^2\langle x-z,u_d\rangle^2   d{\rm vol}_\rho
=\bigg\{
\begin{array}{rc}
 \beta_1+O(\varepsilon^5),    & k=d\\ \beta_2+O(\varepsilon^5),& k\neq d 
\end{array},$$
$$-\hspace{-3.8mm}\int_{\mathcal{B}}\langle x-z,u_j\rangle^2\langle x-z,u_l\rangle^2   d{\rm vol}_\rho
=\bigg\{
\begin{array}{rc}
 \beta_3+O(\varepsilon^5),    & j=l\\ \beta_4+O(\varepsilon^5),& j\neq l 
\end{array};$$
otherwise,
$$-\hspace{-3.8mm}\int_{\mathcal{B}}\langle x-z,u_j\rangle\langle x-z,u_k\rangle\langle x-z,u_m\rangle\langle x-z,u_n\rangle   d{\rm vol}_\rho=O(\varepsilon^5).$$
\end{Lma}

\begin{proof}
\ \\
(i)
Recall that we denote $\theta=\bar\theta^je_j$.
As in equation (\ref{pq}), by Theorem B.1 in \cite{SW12}, there exists an orthonormal basis $\{e_l\}_{l=1}^d$ of $T_{z}M$ so that
$\langle x-z, e_l \rangle =\langle x-z,u_l \rangle +O(\varepsilon^{\tau+2})$. Hence we compute 
$\int_{\mathcal{B}}\langle x-z,e_k\rangle d{\rm vol}_\rho$ and use it to estimate $\int_{\mathcal{B}}\langle x-z,u_k\rangle d{\rm vol}_\rho$. Note that $\langle x-z,e_k\rangle$ is easier to compute than $\langle x-z,u_k\rangle$  because $e_k$ are the real tangent vector in $T_zM$ and thus $\langle {\rm I\!I}_z(\theta,\theta), e_k\rangle=0$. Therefore,

\begin{align*}
&\ \int_{\mathcal{B}}\langle x-z,e_k\rangle d{\rm vol}_\rho\\
=&\ \int_{S^{d-1}}\int_0^{\widetilde r(\theta)}
\left( t\bar \theta^k+O(t^3)\right)
\left(\rho+t\nabla_\theta \rho+O(t^2)\right)
\left( t^{d-1}+t^{d+1}R_1+O(t^{d+2})\right)dtdS_\theta\\
=&\ \int_{S^{d-1}}\int_0^{\widetilde r(\theta)}
\left( t^d\bar \theta^k \rho(z) +t^{d+1}\bar \theta^k \nabla_\theta \rho|_z+O(t^{d+2})\right)dtdS_\theta\\
=&\ \rho(z)\int_{S^{d-1}}\int_0^{\widetilde r(\theta)}
\bar \theta^k t^d dtdS_\theta
+
\sum_{j=1}^d\partial_j\rho\int_{S^{d-1}}\int_0^{\widetilde r(\theta)}\bar\theta^j\bar \theta^k t^{d+1} dtdS_\theta+O(\varepsilon^{d+3})\\
=&\ \bigg\{
\begin{array}{rc}
 \rho(z)C_{1,0}\varepsilon^{d+1}
+\partial_d\rho|_{z}C_{2,0}\varepsilon^{d+2}+O(\varepsilon^{d+3}),    & k=d\\ 
\partial_k\rho|_{z}C_{0,2}\varepsilon^{d+2}+O(\varepsilon^{d+3}),& k\neq d 
\end{array}.
\end{align*}
Recall that ${\rm Vol}_\rho(\mathcal{B})=\varepsilon^d C_{0,0}\rho(z)\left(1+\gamma_1\frac{\partial_d\rho}{\rho}\big|_z+O(\varepsilon^{2})\right)$,
so
\begin{align*}
&-\hspace{-3.8mm}\int_{\mathcal{B}}\langle x-z,u_k\rangle d{\rm vol}_\rho\\
=&\ (\varepsilon^d C_{0,0}\rho(z))^{-1}\left(1-\gamma_1\frac{\partial_d\rho}{\rho}\Big|_z+O(\varepsilon^{2})\right)\int_{\mathcal{B}}\langle x-z,e_k\rangle+O(\varepsilon^{\tau+2}) d{\rm vol}_\rho\\
=&\ \bigg\{
\begin{array}{rc}
 \gamma_1
+(\alpha_1-\gamma_1^2)\frac{\partial_d\rho}{\rho}\big|_{z}+O(\varepsilon^{\tau+2}),    & k=d\\ 
 \alpha_2\frac{\partial_k\rho}{\rho}\big|_{z}+O(\varepsilon^{\tau+2}),& k\neq d 
\end{array}.
\end{align*}

\ \\
{\rm (ii)}
\begin{align*}
&\ \int_{\mathcal{B}}\langle x-z,e_k\rangle^2 d{\rm vol}_\rho\\
=&\ \int_{S^{d-1}}\int_0^{\widetilde r(\theta)}
\left( t\bar \theta^k+O(t^3)\right)^2
\left(\rho(z)+t\nabla_\theta \rho|_z+O(t^2)\right)
\left( t^{d-1}+t^{d+1}R_1+O(t^{d+2})\right)dtdS_\theta\\
=&\ \int_{S^{d-1}}\int_0^{\widetilde r(\theta)}
\left( t^{d+1} (\bar\theta^k)^2 \rho(z) +t^{d+2} (\bar\theta^k)^2 \nabla_\theta \rho|_z+O(t^{d+3})\right)dtdS_\theta\\
=&\ \rho(z)\int_{S^{d-1}}\int_0^{\widetilde r(\theta)}
 (\bar\theta^k)^2 t^{d+1} dtdS_\theta
+
\sum_{j=1}^d\partial_j\rho\int_{S^{d-1}}\int_0^{\widetilde r(\theta)}\bar\theta^j(\bar \theta^k)^2 t^{d+2} dtdS_\theta+O(\varepsilon^{d+4})\\
=&\ \bigg\{
\begin{array}{rc}
 \rho(z)C_{2,0}\varepsilon^{d+2}
+\partial_d\rho|_{z}C_{3,0}\varepsilon^{d+3}+O(\varepsilon^{d+4}),    & k=d\\ 
 \rho(z)C_{0,2}\varepsilon^{d+2}
+\partial_d\rho|_{z}C_{1,2}\varepsilon^{d+3}+O(\varepsilon^{d+4}),& k\neq d 
\end{array}.
\end{align*}
As before, we have $\langle x-z, u_k \rangle^2 =\langle x-z,e_k \rangle^2+O(\varepsilon^{\tau+3})$.
Therefore,
$$-\hspace{-3.8mm}\int_{\mathcal{B}}\langle x-z,u_k\rangle^2   d{\rm vol}_\rho
=\bigg\{
\begin{array}{rc}
 \alpha_1+(\mu_1-\gamma_1\alpha_1)\frac{\partial_d\rho}{\rho}\big|_{z}+O(\varepsilon^{\tau+3}),    & k=d\\ 
 \alpha_2+(\mu_2-\gamma_1\alpha_2)\frac{\partial_d\rho}{\rho}\big|_{z}+O(\varepsilon^{\tau+3}),& k\neq d 
\end{array}.$$
Similarly, for $j=1,\dots, d-1$, one can derive
$$-\hspace{-3.8mm}\int_{\mathcal{B}}\langle x-z,u_k\rangle\langle x-z,u_j\rangle   d{\rm vol}_\rho
=\bigg\{
\begin{array}{rc}
 \mu_2\frac{\partial_j\rho}{\rho}\big|_{z}+O(\varepsilon^{\tau+3}),    & k=d\\ 
 O(\varepsilon^{\tau+3}),& k\neq d 
\end{array}.$$
\ \\
{\rm (iii)}
For all $k=1,\dots,d$, 
\begin{align*}
&\ \int_{\mathcal{B}}\langle x-z,e_k\rangle^2\langle x-z,e_d\rangle d{\rm vol}_\rho\\
=&\ \int_{S^{d-1}}\int_0^{\widetilde r(\theta)}
\left( t\bar \theta^k+O(t^3)\right)^2\left( t\bar \theta^d+O(t^3)\right)
\left(\rho(z)+t\nabla_\theta \rho|_z+O(t^2)\right)
\left( t^{d-1}+O(t^{d+1})\right)dtdS_\theta\\
=&\ \int_{S^{d-1}}\int_0^{\widetilde r(\theta)}
\left( t^{d+2} (\bar\theta^k)^2\bar\theta^d \rho(z) +t^{d+3} (\bar\theta^k)^2 (\bar\theta^d)^2\partial_d \rho|_z+O(t^{d+4})\right)dtdS_\theta\\
=&\ \bigg\{
\begin{array}{rc}
 \rho(z)C_{3,0}\varepsilon^{d+3}
+\partial_d\rho|_{z}C_{4,0}\varepsilon^{d+4}+O(\varepsilon^{d+5}),    & k=d\\ 
 \rho(z)C_{1,2}\varepsilon^{d+3}
+\partial_d\rho|_{z}C_{2,2}\varepsilon^{d+4}+O(\varepsilon^{d+5}),& k\neq d 
\end{array}.
\end{align*}
Moreover, for all $j=1,\dots, d-1$ and $k\neq j$, 
\begin{align*}
&\ \int_{\mathcal{B}}\langle x-z,e_k\rangle^2\langle x-z,e_j\rangle d{\rm vol}_\rho\\
=&\ \int_{S^{d-1}}\int_0^{\widetilde r(\theta)}
\left( t\bar \theta^k+O(t^3)\right)^2\left( t\bar \theta^j+O(t^3)\right)
\left(\rho(z)+t\nabla_\theta \rho|_z+O(t^2)\right)
\left( t^{d-1}+O(t^{d+1})\right)dtdS_\theta\\
=&\ \int_{S^{d-1}}\int_0^{\widetilde r(\theta)}
\left( t^{d+2} (\bar\theta^k)^2\bar\theta^j \rho(z) +t^{d+3} (\bar\theta^k)^2 (\bar\theta^j)^2\partial_j \rho|_z+O(t^{d+4})\right)dtdS_\theta\\
=&\ \bigg\{
\begin{array}{rc}
\partial_j\rho|_{z}C_{2,2}\varepsilon^{d+4}+O(\varepsilon^{d+5}),    & k=d\\ 
\partial_j\rho|_{z}C_{0,2,2}\varepsilon^{d+4}+O(\varepsilon^{d+5}),& k\neq d 
\end{array}.
\end{align*}

Therefore, 
$$-\hspace{-3.8mm}\int_{\mathcal{B}}\langle x-z,u_k\rangle^2\langle x-z,u_d\rangle   d{\rm vol}_\rho
=\bigg\{
\begin{array}{rc}
 \mu_1+\beta_1\frac{\partial_d\rho}{\rho}\big|_{z}+O(\varepsilon^{\tau+4}),    & k=d\\ 
 \mu_2+\beta_2\frac{\partial_d\rho}{\rho}\big|_{z}+O(\varepsilon^{\tau+4}),& k\neq d 
\end{array}.$$
$$-\hspace{-3.8mm}\int_{\mathcal{B}}\langle x-z,u_k\rangle^2\langle x-z,u_j\rangle   d{\rm vol}_\rho
=\bigg\{
\begin{array}{rc}
 \beta_2\frac{\partial_j\rho}{\rho}\big|_{z}+O(\varepsilon^{\tau+4}),    & k=d\neq j\\ 
  \beta_4\frac{\partial_j\rho}{\rho}\big|_{z}+O(\varepsilon^{\tau+4}),& k\neq d, k\neq j 
\end{array}.$$

At last, it is easy to see that the symmetry of domain gives 
$$-\hspace{-3.8mm}\int_{\mathcal{B}}\langle x-z,u_k\rangle^3d{\rm vol}_\rho= O(\varepsilon^{\tau+4}) \mbox{  and  }  -\hspace{-3.8mm}\int_{\mathcal{B}}\langle x-z,u_k\rangle\langle x-z,u_j\rangle\langle x-z,u_l\rangle   d{\rm vol}_\rho
= O(\varepsilon^{\tau+4})$$
when $k,j,l$ are distinct indices selected from $1,\dots,d-1$.
\ \\
{\rm (iv)} The proof is the same as the above. We leave the routine proof 
as an exercise for interested readers.
\end{proof}

\subsection{More integrals}

Below, we prepare more technical lemmas that are needed when we apply the direct expansion method shown in Section \ref{Appendix proof direct approach}.
Let us start by computing some integrals in the Euclidean space. 
\begin{Lma}\label{ints}
Let $B_r^d$ be the ball of radius $r$ in $\mathbb{R}^d$, $S^{d-1}$ be the unit sphere and 
$|B^d|$ denote the volume of the unit ball. Given an orthonormal basis $\{E_i\}\in T_0\mathbb{R}^d=\mathbb{R}^d$, for $x\in\mathbb{R}^d$ and $\theta\in S^{d-1}$ we denote $x^i:=\langle x,E_i\rangle$ and $\theta^i:=\langle \theta,E_i\rangle$. Then
{\allowdisplaybreaks
\begin{align*}
    I\! I_B(r)&:=\int_{B_r^d} (x^i)^4 dx=\frac{3}{(d+2)(d+4)}|B^d|r^{d+4},\\
    J\! I_B(r)&:=\int_{B_r^d} (x^j)^2(x^i)^2 dx=\frac{1}{(d+2)(d+4)}|B^d|r^{d+4},\\
    I\! I_S&:=\int_{S^{d-1}} (\theta^i)^4 dS_{\theta}=\frac{3}{d+2}|B^d|,\\
    J\! I_S&:=\int_{S^{d-1}} (\theta^j)^2(\theta^i)^2 dS_{\theta}=\frac{1}{d+2}|B^d|,\\
    I\! I\! I_S&:=\int_{S^{d-1}} (\theta^i)^6 dS_{\theta}=\frac{15}{(d+2)(d+4)}|B^d|,\\
    J\! J\! I_S&:=\int_{S^{d-1}} (\theta^j)^4(\theta^i)^2 dS_{\theta}=\frac{3}{(d+2)(d+4)}|B^d|,\\
    K\! J\! I_S&:=\int_{S^{d-1}} (\theta^k)^2(\theta^j)^2(\theta^i)^2 dS_{\theta}=\frac{1}{(d+2)(d+4)}|B^d|,\\
    J\! J\! I\! I_S&:=\int_{S^{d-1}} (\theta^j)^4(\theta^i)^4 dS_{\theta}=\frac{9}{(d+2)(d+4)(d+6)}|B^d|.
\end{align*}
}

\end{Lma}

\begin{proof}
Using spherical coordinates, one has

\begin{align*}
{\rm Vol}(B_r^d)
=&\ \int_0^{2\pi}\int_0^{\pi}\cdots\int_0^{\pi}\int_0^r\rho^{d-1}\sin^{d-2}\varphi_{d-1}\sin^{d-3}\varphi_{d-2}\cdots \sin\varphi_{2} d\rho d\varphi_{d-1}\cdots d\varphi_{2} d\varphi_{1}\\
=&\ \frac{r^d}{d}\cdot 2\pi\cdot\left(\int_0^\pi\sin^{d-2}\varphi d\varphi\right)\cdots\left(\int_0^\pi\sin\varphi d\varphi\right).
\end{align*}
By using formulas $C_d:=\int_0^{\pi/2}\sin^dtdt=\frac{\sqrt{\pi}\Gamma(\frac{d+1}{2})}{2\Gamma(\frac{d}{2}+1)}$ and $\int_0^\pi \cos^{2m}\varphi\sin^{d-2}\varphi d\varphi=\frac{\Gamma(\frac{2m+1}{2})\Gamma(\frac{d-1}{2})}{\Gamma(\frac{2m+d}{2})}$, one can derive 
$$\int_{B_r^d} (x^i)^{2m}dx=\frac{1\cdot3\cdot\cdots\cdot (2m-1)}{(d+2)(d+4)\cdots(d+2m)}|B^d|r^{d+2m}.$$
Note that $\int_{B_1^d} (x^i)^4dx=\int_{S^{d-1}}\int_0^1 t^4(\theta^i)^4t^{d-1}dtdS_{\theta}=\int_{S^{d-1}}(\theta^i)^4dS_{\theta}\cdot \frac{1}{d+4}$, so $I\! I_S=\int_{S^{d-1}} (\theta^i)^4 dS_{\theta}=\frac{3}{d+2}|B^d|$. All other integrals can be derived by similar computations.
\end{proof}

Note that all integrals with some odd powers of $\theta^j$ are zero because of the symmetries of $\theta$ and $S^{d-1}$. For instance, $\int_{S^{d-1}} (\theta^j)(\theta^i)^2 dS_{\theta}=0$ and $\int_{S^{d-1}} (\theta^k)(\theta^j)^2(\theta^i)^3 dS_{\theta}=0$. We need the following two technical lemmas which concern the integral of tensors.

\begin{Lma}\label{T}
For any symmetric contravariant continuous $2$-tensor $T$ defined around $x\in M^d$ and any orthonormal basis $\{e_j\}_{j=1}^d$ of $T_xM$, we have
$$\int_{S^{d-1}}T_x(\theta,\theta) dS_\theta
=|B^d|{\rm tr}T_x,\ \ \int_{S^{d-1}}T_x(\theta,\theta)\langle\theta,e_1\rangle^2dS_\theta
=\frac{|B^d|}{d+2}\left(2T_x(e_1,e_1)+{\rm tr}T_x\right),$$
and
$$\int_{S^{d-1}}T_x(\theta,\theta)\langle\theta,e_1\rangle\langle\theta,e_2\rangle dS_\theta
=\frac{2|B^d|}{d+2}T_x(e_1,e_2),$$
where $S^{d-1}\subset T_xM$ is the unit sphere centered at $0=\exp^{-1}(x)$. In particular, $$
\int_{S^{d-1}}\nabla_\theta^2 f(x) dS_\theta= |B^d|\Delta f(x)\ \ \mbox{ and }\ \  \int_{S^{d-1}}Ric_x(\theta,\theta) dS_\theta= |B^d|S(x),
$$
where $f\in C^2(M)$, $Ric_x$ is the Ricci tensor at $x$ and $S(x)$ is the scalar curvature at $x$.
\end{Lma}

\begin{proof}
We omit the subscript $x$ of $T_x$. Let $T$ be diagonalized at $x$ by an orthonormal basis which can be extended into a local frame $\{E_j\}_{j=1}^d$ around $x$. Let $e_j=\sum_{k=1}^d\eta_j^kE_k$ and $\theta=\sum_{k=1}^d\theta^kE_k$.  Note that $\sum_{k=1}^d(\eta_j^k)^2=1$ for any $j$. Then
\begin{align*}
\int_{S^{d-1}}T(\theta,\theta)\langle\theta,e_1\rangle^2dS_\theta
&\,=\int_{S^{d-1}} \sum_{j=1}^d T_{jj}(\theta^j)^2\sum_{k=1}^d(\theta^k)^2(\eta_1^k)^2dS_\theta\\
&\,= \sum_{j,k=1}^d T_{jj}(\eta_1^k)^2\int_{S^{d-1}}(\theta^j)^2(\theta^k)^2dS_\theta\,.
\end{align*}
Recall that terms with odd orders of $\theta^i$ vanish by symmetry. Lemma \ref{ints} says that $\int_{S^{d-1}}(\theta^j)^2(\theta^k)^2dS_\theta=I\!I_S$ when $k=j$ and $\int_{S^{d-1}}(\theta^j)^2(\theta^k)^2dS_\theta=J\!I_S$ when $k\neq j$. So  
\begin{align*}
\int_{S^{d-1}}T(\theta,\theta)\langle\theta,e_1\rangle^2dS_\theta
=&\ \sum_{j,k=1;k\neq j}^d T_{jj}(\eta_1^k)^2\cdot J\!I_S+\sum_{j=1}^d T_{jj}(\eta_1^j)^2\cdot I\!I_S\\
=&\ \sum_{k,j=1}^d T_{jj}(\eta_1^k)^2\cdot J\!I_S+\sum_{j=1}^d T_{jj}(\eta_1^j)^2\cdot\left(I\! I_S-J\! I_S\right)\\
=&\ {\rm tr}T\cdot J\! I_S  + T(e_1,e_1)\cdot \left(I\! I_S-J\! I_S\right),
\end{align*}
where the last equality comes from the symmetry of $T$. We get the claim by plugging Lemma \ref{ints}. Other statements can be similarly derived.
\end{proof}

\begin{Lma}\label{sff}
Denote the square norm of the second fundamental form ${\rm I\!I}_x$ by $|{\rm A}|^2(x)$ and the mean curvature vector by ${\bf H}_x:={\rm tr}{\rm I\!I}_x$. We have

$$\int_{S^{d-1}}|{\rm I\! I}_x(\theta,\theta)|^2 dS_\theta=\frac{|B^{d}|}{d+2}\left(2|{\rm A}|^2(x)+|{\bf H}_x|^2\right),$$
\begin{align*}
&\int_{S^{d-1}}|{\rm I\!I}_x(\theta,\theta)|
^2\langle\theta,e_1\rangle^2dS_\theta\\
=&\, \frac{|B^d|}{(d+2)(d+4)}\left(
12{\bf H}_x\cdot{\rm I\!I}_x(e_1,e_1)-8Ric_x(e_1,e_1)
+2|{\rm A}|^2(x) +|{\bf H}_x|^2 \right)
\end{align*}
and 
$$\int_{S^{d-1}}|{\rm I\!I}_x(\theta,\theta)|^2\langle\theta,e_1\rangle\langle\theta,e_2\rangle dS_\theta
=\frac{|B^d|}{(d+2)(d+4)}\left(12{\bf H}_x\cdot{\rm I\!I}_x(e_1,e_2)-8Ric_x(e_1,e_2)\right)$$
for any orthonormal basis $\{e_j\}_{j=1}^d$ of $T_xM$.
\end{Lma}

\begin{proof}
Since all the tensors are evaluated at the point $x$, we may omit the notation $x$ in $Ric_x$, ${\rm I\!I}_x$, ${\bf H}_x$ and $|A|^2(x)$ in this proof. Let ${\rm I\!I}$ be diagonalized at $x$ by an orthonormal basis $\{E_j\}_{j=1}^d\subset T_xM$. Extend $\{E_j\}$ to be a local frame around $x$. Denote $\theta=\sum_{j=1}^d\theta^jE_j$ and ${\rm I\!I}(E_j,E_j)=\sum_{l=d+1}^ph_{jj}^l e_l$, where $\{e_l\}_{j=d+1}^p$ is an orthonormal basis of $N_xM$, the normal space at $x$. In the following, we assume the codimension $p-d$ is $1$ and denote ${\rm I\!I}(E_j,E_j)$ simply by $h_{jj}$. The general case can be derived in exactly the same way and has the same conclusion.

Recall that $\displaystyle|{\rm A}|^2=\sum_{i=1}^dh_{jj}^2$ and $\displaystyle|{\bf H}|^2=\left(\sum_{i=1}^dh_{jj}\right)^2=\sum_{j,k=1}^d h_{jj}h_{kk}$, one can derive
{\allowdisplaybreaks
\begin{align*}
\int_{S^{d-1}}|{\rm I\! I}(\theta,\theta)|^2 dS_\theta
=&\ \sum_{j,k=1}^d h_{jj}h_{kk}\int_{S^{d-1}} (\theta^j)^2(\theta^k)^2 dS_\theta\\
=&\ \sum_{j=1}^d h_{jj}h_{jj}\int_{S^{d-1}} (\theta^j)^4 dS_\theta
+\sum_{j\neq k} h_{jj}h_{kk}\int_{S^{d-1}} (\theta^j)^2(\theta^k)^2 dS_\theta
\\
=&\ \sum_{j=1}^d h_{jj}h_{jj}\cdot I\!I_S
+\sum_{j\neq k} h_{jj}h_{kk}\cdot J\!I_S
\\
=&\ \sum_{j=1}^d h_{jj}h_{jj}\cdot I\!I_S
+\sum_{j,k=1}^d h_{jj}h_{kk}\cdot J\!I_S -\sum_{j=1}^d h_{jj}h_{jj}\cdot J\!I_S  
\\
=&\ (I\!I_S-J\!I_S  )\cdot \sum_{j=1}^d h_{jj}h_{jj}
+J\!I_S\cdot\sum_{j,k=1}^d h_{jj}h_{kk}
\\
=&\ \frac{1}{d+2}|B^{d}|\left(2|{\rm A}|^2+|{\bf H}|^2\right). 
\end{align*}
}
Let $e_j=\sum_{k=1}^d\eta_j^kE_k$. Recall that Gauss-Codazzi equation for a submanifold in $\mathbb{R}^p$ gives 
$$
\sum_{j=1}^d{\rm I\!I}(e_s,E_j){\rm I\!I}(e_t,E_j)={\bf H}\cdot{\rm I\!I}(e_s,e_t)-Ric(e_s,e_t)
$$ 
for all $s,t=1,\dots, d$. Hence, 
{\allowdisplaybreaks
\begin{align*}
&\int_{S^{d-1}}|{\rm I\!I}(\theta,\theta)|
^2\langle\theta,e_1\rangle^2dS_\theta\\
=&\ \sum_{j,k,l=1}^d h_{jj}h_{kk}(\eta_1^l)^2 \int_{S^{d-1}}(\theta^j)^2(\theta^k)^2(\theta^l)^2 dS_\theta\\
=&\ \sum_{j=1}^d h_{jj}h_{jj}(\eta_1^j)^2\cdot I\!I\!I_S
+  \sum_{l\neq j;j,l=1}^d h_{jj}h_{jj}(\eta_1^l)^2\cdot J\!J\!I_S \\
&\ + 2\sum_{k\neq j;j,k=1}^d h_{jj}h_{kk}(\eta_1^j)^2\cdot J\!J\!I_S+ \sum_{l\neq k\neq j;j,k,l=1}^d h_{jj}h_{kk}(\eta_1^l)^2\cdot K\!J\!I_S \\
=&\ \left({\bf H}\cdot{\rm I\!I}(e_1,e_1)-Ric(e_1,e_1)\right)\cdot \left(I\!I\!I_S -3 J\!J\!I_S+ 2K\!J\!I_S\right)\\
&\ + (2 {\bf H}\cdot{\rm I\!I}(e_1,e_1)+|{\rm A}|^2)\cdot\left( J\!J\!I_S - K\!J\!I_S\right)
+ |{\bf H}|^2\cdot K\!J\!I_S \\
=&\ \frac{|B^d|}{(d+2)(d+4)}\left(
12{\bf H}\cdot{\rm I\!I}(e_1,e_1)-8Ric(e_1,e_1)
+2|{\rm A}|^2 +|{\bf H}|^2 \right).
\end{align*}
}The last equation can be similarly derived.
\end{proof}

\begin{Rmk}
Taking trace of the Gauss-Codazzi equation, one can see that 
\begin{align*}
\sum_{s=1}^d\sum_{j=1}^d{\rm I\!I}(e_s,E_j){\rm I\!I}(e_s,E_j)&\,=\sum_{s=1}^d\sum_{j=1}^d{\rm I\!I}(\sum_{k=1}^d\eta_s^kE_k,E_j){\rm I\!I}(\sum_{k=1}^d\eta_s^kE_k,E_j)\\
&\,=\sum_{j=1}^d \sum_{s=1}^d(\eta_s^jh_{jj})^2=|{\rm A}|^2
\end{align*}
equals to
$\sum_{s=1}^d\left({\bf H}\cdot{\rm I\!I}(e_s,e_s)-Ric(e_s,e_s)\right)=|{\bf H}|^2-S$. Hence the scalar curvature $S=|{\bf H}|^2-|{\rm A}|^2$.
\end{Rmk}

At last, we need the following lemma which can be derived by using the same method as previous ones.

\begin{Lma}\label{Aij}
For distinct indices $s,l\in\{1,\dots,d\}$ and any orthonormal basis $\{e_j\}_{j=1}^d$ of $T_xM$, we have
$$\int_{S^{d-1}}\langle {\rm I\!I}_x(\theta,\theta),{\rm I\!I}_x(\theta,e_s) \rangle\langle\theta,e_l\rangle dS_\theta
=\frac{|B^d|}{d+2}\left(3{\bf H}_x\cdot{\rm I\!I}_x(e_s,e_l)-2Ric_x(e_s,e_l)\right).$$
\end{Lma}

\begin{proof}
We omit the lower subscript $x$ for all tensors. Using frames $\{E_j\}$ and $\{e_j\}$ as in the previous lemmas,
we have 
\begin{align*}
&\int_{S^{d-1}}\langle {\rm I\!I}(\theta,\theta),{\rm I\!I}(\theta,e_s) \rangle\langle\theta,e_l\rangle dS_\theta\\
=&\ \sum_{j,k,m=1}^dh_{jj}h_{kk}\int_{S^{d-1}}(\theta^j)^2\theta^k\eta_s^k\theta^m\eta_l^m dS_{\theta}\\
=&\ \sum_{j,k=1}^dh_{jj}h_{kk}\eta_s^k\eta_l^k\int_{S^{d-1}}(\theta^j)^2(\theta^k)^2 dS_{\theta}\\
=&\ \sum_{j=1}^d{\rm I\!I}(E_j,e_s){\rm I\!I}(E_j,e_l)\frac{2|B^d|}{d+2}
+{\bf H}\cdot{\rm I\!I}(e_s,e_l)\frac{|B^d|}{d+2}\\
=&\ \frac{|B^d|}{d+2}\left(3{\bf H}\cdot{\rm I\!I}(e_s,e_l)-2Ric(e_s,e_l)\right).
\end{align*}
\end{proof}

To continue, recall from the above that we have
\begin{align*}
f(x)=&\ f(z)+f_1(\theta)t+f_2(\theta)t^2+f_3(\theta)t^3+O(t^4);\\
\iota(x) -\iota(z) 
=&\  K_1(\theta)t+K_2(\theta)t^2+K_3(\theta)t^3+K_4(\theta)t^4+O(t^5);\\
\tilde{\varepsilon}(\theta)
=&\ \varepsilon + H_1(\theta)\varepsilon^3 +H_2(\theta)\varepsilon^4+O(\varepsilon^{5});\\
d{\rm vol}=&\ \left(t^{d-1}+R_1(\theta)t^{d+1}+R_2(\theta)t^{d+2}+O(t^{d+3})\right)dtdS_\theta.
\end{align*}
Moreover, deriving from Lemma \ref{T} and Lemma \ref{sff} that
\begin{align*}
\int_{S^{d-1}} R_1(\theta)+(d+2)H_1(\theta) dS_\theta&\,=|B^d|\left(-\frac{1}{6}S(z)+\frac{1}{12}|A|^2(z)+\frac{1}{24}|{\bf H}_{z}|^2\right)\\
&\,=:-|B^d|(d+2)\Lambda,
\end{align*}
where 
$$
\Lambda=\frac{1}{8(d+2)}\left(|{\bf H}_{z}|^2-2|A|^2(z)\right),
$$ 
one obtains (cf. Theorem 3.1 in \cite{G73})
\begin{align*}
{\rm Vol}(\iota(M)\cap B_\varepsilon^p(z))=&\ \int_{S^{d-1}}\int_0^{\tilde\varepsilon(\theta)}t^{d-1}+R_1(\theta)t^{d+1}+R_2(\theta)t^{d+2}+O(t^{d+3})dtdS_\theta\\
=&\ \varepsilon^d|B^d|\left(1-\Lambda\varepsilon^2+O(\varepsilon^3)\right).
\end{align*}

Suppose that there is a weighted function $\rho(x)\in C^4(M)$ on $M$, e.g. the probability density function (p.d.f.) associated with a random vector, we can similarly compute the weighted volume \begin{align}\label{taylor volume}
\begin{split}
&\ {\rm Vol}_\rho(\iota(M)\cap B_\varepsilon^p(z))\\
=&\ \int_{S^{d-1}}\int_0^{\tilde\varepsilon(\theta)}(\rho(z)+\rho_1(\theta)t+\rho_2(\theta)t^2+O(t^3))(t^{d-1}+R_1(\theta)t^{d+1}+O(t^{d+2}))dtdS_\theta\\
=&\ \varepsilon^d|B^d|\left(\rho(z)+\varepsilon^2\left(-\Lambda\rho(z) +\frac{1}{2(d+2)}\Delta\rho(z)\right)+O(\varepsilon^4)\right).
\end{split}
\end{align}
Hence,
$$\frac{1}{{\rm Vol}_\rho(\iota(M)\cap B_\varepsilon^p(z))}
= \frac{1}{\varepsilon^d|B^d|}\left(\rho^{-1}+\varepsilon^2\rho^{-2}\left(\Lambda\rho-\frac{1}{2(d+2)}\Delta\rho\right)+O(\varepsilon^4)\right).$$

\begin{Thm} \label{intf}
Let $\alpha=\frac{\varepsilon^2}{d+2}$, $\beta=\frac{\varepsilon^4}{(d+2)(d+4)}$, and $B_\varepsilon^p(z)\subset\mathbb{R}^p$ be the $\varepsilon$-ball centered at $z$. Let $\{e_j\}_{j=1}^d$ be an orthonormal basis of $T_{z}M$ and $\widetilde{B}:=\iota(M)\cap B_\varepsilon^p(z)$. Suppose that the sample distribution on $M$ is govern by the density function $\rho(x)$ and denote $\rho(x)d{\rm vol}$ by $d{\rm vol}_{\rho}$. Then for $f\in C^4(M)$, we have\\
{\rm (i)}
$$-\hspace{-3.8mm}\int_{\widetilde B}f(x)d{\rm vol}_{\rho}
= f(z)+\frac{1}{2}\alpha\left(\Delta f(z)+2\rho^{-1}(z)\langle \nabla f,\nabla\rho\rangle(z)\right)+ O(\varepsilon^4).
$$\\
{\rm (ii)}
For all $j=1,\dots,d$,
$$-\hspace{-3.8mm}\int_{\widetilde B}\langle x-z,e_j\rangle f(x)d{\rm vol}_{\rho}
=\alpha\left(\nabla_jf|_{z}+\rho^{-1}(z)f(z)\nabla_j\rho|_{z}\right)+O(\varepsilon^{4}).
$$
{\rm (iii)} For all $j=1,\dots,d$,
\begin{align*}
&\int_{\widetilde B}\langle x-z,e_j\rangle^2f(x){\rm dvol}_{\rho}\\
=&\ \alpha f(z)+\beta\rho^{-1}(z)\left(\nabla_j^2(f\rho)|_{z}+\frac{1}{2}\Delta (f\rho)(z)-\frac{d
+4}{2(d+2)}f(z)\Delta\rho(z)\right)\\
&\ +\beta\left(2 \Lambda-\frac{1}{2}{\bf H}_{z}\cdot{\rm I\!I}_{z}(e_j,e_j)\right)f(z)+O(\varepsilon^6)\\
=&\ \alpha f(z)+\beta\left(\nabla_j^2f|_{z}+\frac{1}{2}\Delta f(z)\right)\\
&\ +\beta \rho^{-1}(z)\left( 2\nabla_jf\nabla_j\rho|_{z}+ \langle \nabla f,\nabla\rho\rangle(z)+f(z)\nabla_j^2\rho|_{z}-\frac{1}{d+2}f(z)\Delta\rho(z)\right)\\
&\ +\beta \left(2 \Lambda-\frac{1}{2}{\bf H}_{z}\cdot{\rm I\!I}_{z}(e_j,e_j)\right)f(z)+O(\varepsilon^6).  
\end{align*}
\\
{\rm (iv)} 
For all $s,l=1,\dots,d$ and $s\neq l$,
\begin{align*}
&-\hspace{-3.8mm}\int_{\widetilde B}\langle x-z,e_s\rangle\langle x-z,e_l\rangle f(x)d{\rm vol}_{\rho}\\
=&\  \beta\rho^{-1}(z)\nabla_s\nabla_l(f\rho)|_{z} +\beta\left(-\frac{1}{2}{\bf H}_{z}\cdot{\rm I\!I}_{z}(e_s,e_l)\right)f(z)
+O(\varepsilon^6)\\
=&\  \beta\nabla_s\nabla_l f|_{z}
+\beta\rho^{-1}(z)\left(\nabla_sf\nabla_l\rho+\nabla_lf\nabla_s\rho+f\nabla_s\nabla_l\rho\right)|_{z}\\
& +\beta\left(-\frac{1}{2}{\bf H}_{z}\cdot{\rm I\!I}_{z}(e_s,e_l)\right)f(z)+O(\varepsilon^6).
\end{align*}
\\
{\rm (v)} 
$$-\hspace{-3.6mm}\int_{\widetilde B}\langle x-z,e_s \rangle\langle x-z,e_l \rangle^2 d{\rm vol}_\rho=
\Big\{
\begin{array}{rc}
 \beta\rho^{-1}(z)\nabla_s\rho|_{z}+O(\varepsilon^6),& s\neq l \\
 3\beta\rho^{-1}(z)\nabla_s\rho|_{z}+O(\varepsilon^6),    & s=l
\end{array}.$$
 
$$-\hspace{-3.6mm}\int_{\widetilde B}\langle x-z,e_j \rangle\langle x-z,e_s \rangle\langle x-z,e_l \rangle d{\rm vol}_\rho=O(\varepsilon^6)\ {\mbox when }\ j,s,l \mbox{ are mutually distinct}.$$
\\
{\rm (vi)} 
For all $s,l=1,\dots,d$,
$$
-\hspace{-3.8mm}\int_{\widetilde B}\langle x-z,e_s\rangle^2\langle x-z,e_l\rangle^2d{\rm vol}_{\rho}
=
\Big\{
\begin{array}{rc}
 \beta+O(\varepsilon^6),& s\neq l \\
 3\beta+O(\varepsilon^6),    & s=l
\end{array}.
$$
When $s,t,l,m$ are mutually distinct, 
$-\hspace{-3.7mm}\int_{\widetilde B}\langle x-z,e_s\rangle\langle x-z,e_l\rangle^3d{\rm vol}_{\rho}$,
$-\hspace{-3.7mm}\int_{\widetilde B}\langle x-z,e_s\rangle\langle x-z,e_t\rangle\langle x-z,e_l\rangle^2d{\rm vol}_{\rho}$
and 
$-\hspace{-3.7mm}\int_{\widetilde B}\langle x-z,e_s\rangle\langle x-z,e_t\rangle\langle x-z,e_l\rangle\langle x-z,e_m\rangle d{\rm vol}_{\rho}$ are all of $O(\varepsilon^6)$.

\end{Thm}

\begin{proof}
{\allowdisplaybreaks
For the sake of simplicity, we sometimes omit the subscript $z$ of tensors evaluated at $z$, namely, we denote $Ric_{z}$, ${\bf H}_{z}$, ${\rm I\! I}_{z}$ and $|A|^2(z)$ simply by $Ric$, ${\bf H}$, ${\rm I\! I}$ and $|A|^2$. 
We also omit the inclusion map $\iota$ and do not distinguish $M$ and $\iota(M)$. 
Note that $K_1,K_3$ and $R_2$ are odd functions of $\theta$, $R_1$ is an even function of $\theta$, and $K_2$ is the second fundamental form and thus is perpendicular to any tangent vector $e_l$.

(i) Since $\int_{S^{d-1}}(f\rho)_1dS_\theta=\int_{S^{d-1}}\langle\nabla (f\rho),\theta\rangle dS_\theta=0$ and $R_2$ is an odd function of $\theta$, by the symmetry of sphere, one has
\begin{align*}
&\ -\hspace{-3.8mm}\int_{\widetilde B}f(x)d{\rm vol}_{\rho}\\
=&\ 
\frac{1}{{\rm Vol}_{\rho}(\widetilde B)}\int_{S^{d-1}}\int_0^{\tilde\varepsilon} \left(f(z)\rho(z)+(f\rho)_1t+(f\rho)_2t^2+O(t^3)\right)
\left(t^{d-1}+R_1t^{d+1}+O(t^{d+2})\right)dtdS_\theta\\
=&\ 
\frac{1}{{\rm Vol}_{\rho}(\widetilde B)}\int_{S^{d-1}}\int_0^{\tilde\varepsilon} f(z)\rho(z)t^{d-1}+t^{d+1}\left((f\rho)R_1+(f\rho)_2\right)+O(t^{d+3})dtdS_\theta\\
=&\ 
\frac{1}{{\rm Vol}_{\rho}(\widetilde B)}\int_{S^{d-1}}
\frac{\varepsilon^d}{d}f(z)\rho(z)
+\frac{\varepsilon^{d+4}}{d+2}\left[f(z)\rho(z)(R_1+(d+2)H_1)+(f\rho)_2\right]+O(\varepsilon^{d+3})dS_\theta.
\end{align*}
Recall that $\int_{S^{d-1}} \nabla^2(f\rho)|_{z}dS_{\theta}=|B^d|\Delta(f\rho)(z)$ (cf. Lemma \ref{T}) and $$\frac{\varepsilon^{d+2}}{d+2}\int_{S^{d-1}} R_1+(d+2)H_1 dS_\theta=\frac{\varepsilon^{d+2}}{d+2}|B^d|\left(-\frac{1}{6}S+\frac{1}{12}|A|^2+\frac{1}{24}|{\bf H}|^2\right)=-\varepsilon^{d+2}|B^d|\Lambda,$$
one has 
\begin{align*}&\
-\hspace{-3.8mm}\int_{\widetilde B}f(x)d{\rm vol}_{\rho}\\
=&\ 
\left(\rho^{-1}(z)+\varepsilon^2\rho^{-2}(z)\left(\Lambda\rho(z)-\frac{1}{2(d+2)}\Delta\rho \right)+O(\varepsilon^4)\right)f(z)\rho(z)\\
&\ +\frac{1}{\varepsilon^d|B^d|}\left(\rho^{-1}(z)+O(\varepsilon^2)\right)
\left[-\varepsilon^{d+2}|B^d|\Lambda(f(z)\rho(z)+\frac{\varepsilon^{d+2}}{d+2}\int_{S^{d-1}}\frac{1}{2}\nabla^2(f\rho)|_{z} dS_\theta+ O(\varepsilon^{d+4})\right]\\
=&\ f(z)+\varepsilon^2\left[f(z)\Lambda -\rho^{-1}\frac{1}{2(d+2)}\Delta\rho-\Lambda f(z) +\rho^{-1} \frac{1}{2(d+2)}\Delta (f\rho)\right]+ O(\varepsilon^4)\\
=&\ f(z)+\frac{\varepsilon^2}{2(d+2)}\left(\Delta f(z)+2\rho^{-1}\langle\nabla f,\nabla\rho \rangle(z)\right)+ O(\varepsilon^4).
\end{align*}
\\
(ii) For simplicity we omit the notation $z$ and denote $f(z)$ and $\rho(z)$ by $f$ and $\rho$.
\begin{align*}
&\  -\hspace{-3.8mm}\int_{\widetilde B} \langle x- z,e_j\rangle f(x)d{\rm vol}_{\rho}\\
=&\ \frac{1}{{\rm Vol}_{\rho}(\widetilde B)} 
\int_{S^{d-1}}\int_0^{\tilde\varepsilon} \langle K_1t+K_3t^3+O(t^4),e_j\rangle\cdot\left(f\rho+(f\rho)_1t+(f\rho)_2t^2+(f\rho)_3t^3+O(t^4)\right)\\
&\hspace{5cm}\cdot\left(t^{d-1}+R_1t^{d+1}+R_2t^{d+2}+O(t^{d+3})\right)dtdS_\theta\\
=&\ 
\frac{1}{\varepsilon^d|B^d|}\left(\rho^{-1}+\varepsilon^2\rho^{-2}\left(\Lambda\rho-\frac{1}{2(d+2)}\Delta\rho\right)+O(\varepsilon^4)\right)\frac{\varepsilon^{d+1}}{d+1} \int_{S^{d-1}} \langle K_1,e_j\rangle f\rho dS_\theta\\
&\ 
+\frac{1}{\varepsilon^d|B^d|}
\left(\rho^{-1}+\varepsilon^2\rho^{-2}\left(\Lambda\rho-\frac{1}{2(d+2)}\Delta\rho\right)\right)\frac{\varepsilon^{d+2}}{d+2}\int_{S^{d-1}} \langle K_1,e_j\rangle (f\rho)_1 dS_\theta\\
&\ +\frac{1}{\varepsilon^d|B^d|}\frac{\varepsilon^{d+3}}{d+3}\rho^{-1}\int_{S^{d-1}} \langle K_1,e_j\rangle (f\rho)_2+\langle K_3,e_j\rangle f\rho + \langle K_1,e_j\rangle f\rho (d+3)H_1\\
&\hspace{5cm}+\langle K_1,e_j\rangle f\rho R_1 dS_\theta\\
&\ +\frac{1}{\varepsilon^d|B^d|}\frac{\varepsilon^{d+4}}{d+4}\rho^{-1}\int_{S^{d-1}} \langle K_1,e_j\rangle (f\rho)_3+\langle K_3,e_j\rangle (f\rho)_1 +\langle K_4,e_j\rangle f+\langle K_1,e_j\rangle f\rho (d+4)H_2\\
&\hspace{4cm}+\langle K_1,e_j\rangle (f\rho)_1(d+4)H_1+\langle K_1,e_j\rangle (f\rho)_1 R_1+\langle K_1,e_j\rangle f\rho R_2 dS_\theta\\
&\ +O(\varepsilon^{5})\\
=&\ 
\frac{1}{\varepsilon^d|B^d|}
\left(\rho^{-1}+\varepsilon^2\rho^{-2}\left(\Lambda\rho-\frac{1}{2(d+2)}\Delta\rho\right)\right)\frac{\varepsilon^{d+2}}{d+2}\int_{S^{d-1}} \langle K_1,e_j\rangle (f\rho)_1 dS_\theta\\
&\ +\frac{1}{\varepsilon^d|B^d|}\frac{\varepsilon^{d+4}}{d+4}\rho^{-1}\int_{S^{d-1}} \langle K_1,e_j\rangle (f\rho)_3+\langle K_3,e_j\rangle (f\rho)_1 +\langle K_4,e_j\rangle f+\langle K_1,e_j\rangle f\rho (d+4)H_2\\
&\hspace{4cm}+\langle K_1,e_j\rangle (f\rho)_1(d+4)H_1+\langle K_1,e_j\rangle (f\rho)_1 R_1+\langle K_1,e_j\rangle f\rho R_2 dS_\theta\\
&\ +O(\varepsilon^{5}).
\end{align*}
Applying the divergence theorem, we have  
\begin{align*}
&\ \int_{S^{d-1}} \langle K_1,e_j\rangle (f\rho)_1 dS_\theta
=\int_{S^{d-1}} \langle \langle\theta,e_j\rangle \nabla (f\rho)|_{z},\theta\rangle dS_\theta\\
=&\ \int_{B^d} {\rm div}( \langle x,e_j\rangle \nabla (f\rho)|_{z}) dx
= \int_{B^d}\nabla_j(f\rho)|_{z}dx
= |B^d|\cdot\nabla_j(f\rho)|_{z}\,.
\end{align*}
Therefore,
$$ -\hspace{-3.8mm}\int_{\widetilde B} \langle x-z,e_j\rangle f(x)d{\rm vol}_{\rho}
=\frac{\varepsilon^2}{d+2}\left(\nabla_jf|_{z}+\rho^{-1}f\nabla_j\rho|_{z}\right) +O(\varepsilon^{4}).$$
\\
(iii) As before, we denote $f(z)$ and $\rho(z)$ by $f$ and $\rho$. Using the fact that the second fundamental form is a normal vector, i.e., $\langle{\rm I\! I}(\theta,\theta), e_j\rangle=0$, we can derive 
$$
\langle K_3,e_j\rangle
= \frac{1}{6}\langle\nabla_\theta{\rm I\! I}(\theta,\theta), e_j\rangle
= -\frac{1}{6}\langle{\rm I\! I}(\theta,\theta), \nabla_\theta e_j\rangle
= -\frac{1}{6}\langle{\rm I\! I}(\theta,\theta), {\rm I\! I}(\theta,e_j)\rangle
$$
and thus, by Lemma \ref{Aij}, $\int_{S^{d-1}}\langle K_1,e_i\rangle\langle K_3,e_j\rangle dS_\theta=-\frac{|B^d|}{6(d+2)}(3{\bf H}\cdot{\rm I\!I}(e_j,e_j)-2Ric(e_j,e_j))$. Therefore,
\begin{align*}
&\ -\hspace{-3.8mm}\int_{\widetilde B}\langle x- z,e_j\rangle^2f(x)d{\rm vol}_{\rho}\\
=&\ \frac{1}{{\rm Vol}_{\rho}(\widetilde B)} \int_{S^{d-1}}\int_0^{\tilde\varepsilon} \langle K_1t+K_3t^3+O(t^4),e_j\rangle^2\cdot\left(f\rho+(f\rho)_1t+(f\rho)_2t^2+(f\rho)_3t^3+O(t^4)\right)\\
&\hspace{3cm}\cdot\left(t^{d-1}+R_1t^{d+1}+R_2t^{d+2}+R_3t^{d+3}+O(t^{d+4})\right)dtdS_\theta\\
=&\ \frac{\varepsilon^{2}}{(d+2)|B^d|}
\left(\rho^{-1}+\varepsilon^2\rho^{-2}\left(\Lambda\rho-\frac{1}{2(d+2)}\Delta\rho\right)+O(\varepsilon^4)\right) \int_{S^{d-1}} \langle K_1,e_j\rangle^2f\rho dS_\theta\\
&\ +\frac{\varepsilon^{4}}{(d+4)|B^d|}\rho^{-1}\int_{S^{d-1}} 2\langle K_1,e_j\rangle\langle K_3,e_j\rangle f\rho + \langle K_1,e_j\rangle^2\left((f\rho)_2+ f\rho(d+4)H_1+f\rho R_1\right) dS_\theta\\
&\ +O(\varepsilon^6)\\
=&\ \frac{\varepsilon^2}{d+2}f+\varepsilon^4\rho^{-1}\left[\frac{1}{d+2}f\left(\Lambda\rho-\frac{1}{2(d+2)}\Delta\rho\right)\right.\\
&\hspace{3cm} -\frac{1}{3(d+2)(d+4)}(3{\bf H}\cdot{\rm I\!I}(e_j,e_j)-2Ric(e_j,e_j))f\rho \\
&\hspace{3cm} \left.+\frac{1}{(d+4)|B^d|}\int_{S^{d-1}} \langle \theta,e_j\rangle^2\left((f\rho)_2+ f\rho(d+4)H_1+f\rho R_1\right) dS_\theta \right]+O(\varepsilon^6).
\end{align*}

Recall that $(f\rho)_2=\frac{1}{2}\nabla_{\theta}^2(f\rho)|_{z}, H_1=\frac{1}{24}|{\rm I\!I}|^2$, and $R_1=-\frac{1}{6}Ric(\theta,\theta)$. 
So, using Lemma \ref{T} and Lemma \ref{sff}, all $Ric$ are cancelled and one obtains 
\begin{align*}
&\ -\hspace{-3.8mm}\int_{\widetilde B}\langle x-z,e_j\rangle^2f(x){\rm dvol}_{\rho}\\
=&\ \frac{\varepsilon^2}{d+2}f
+\frac{\varepsilon^4}{d+2}f\rho^{-1}\left(\Lambda\rho-\frac{1}{2(d+2)}\Delta\rho\right)\\
&\ +\frac{\varepsilon^4}{(d+2)(d+4)}\left[f\left(-\frac{1}{2}{\bf H}\cdot{\rm I\!I}(e_j,e_j)-(d+2)\Lambda\right)
+\rho^{-1}\left(\nabla_j^2(f\rho)+\frac{1}{2}\Delta(f\rho)\right)\right] +O(\varepsilon^6)\\
=&\ \frac{\varepsilon^2}{d+2}f+\frac{\varepsilon^4}{(d+2)(d+4)}\left(\nabla_j^2f|+\frac{1}{2}\Delta f\right)\\
& +\frac{\varepsilon^4}{(d+2)(d+4)} \rho^{-1}\left( 2\nabla_jf\nabla_j\rho+ \langle \nabla f,\nabla\rho\rangle+f\nabla_j^2\rho-\frac{1}{d+2}f\Delta\rho\right)\\
& +\frac{\varepsilon^4}{(d+2)(d+4)} \left(2 \Lambda-\frac{1}{2}{\bf H}\cdot{\rm I\!I}(e_j,e_j)\right)f+O(\varepsilon^6).   
\end{align*}
\\
(iv) Similarly, we have 
\begin{align*} 
&\ -\hspace{-3.8mm}\int_{\widetilde B}\langle x-z,e_s\rangle\langle x-z,e_l\rangle f(x)d{\rm vol}_{\rho}\\
=&\ \frac{1}{{\rm Vol}_{\rho}(\widetilde B)} \int_{S^{d-1}}\int_0^{\tilde\varepsilon(\theta)} \langle K_1t+K_3t^3+O(t^4),e_s\rangle\langle K_1t+K_3t^3+O(t^4),e_l\rangle\\
&\hspace{6cm}\cdot\left(f\rho+(f\rho)_1t+(f\rho)_2t^2+(f\rho)_3t^3+O(t^4)\right)\\
&\hspace{6cm}\cdot\left(t^{d-1}+R_1t^{d+1}+R_2t^{d+2}+O(t^{d+3})\right)dtdS_\theta\\
=&\ \frac{\varepsilon^{2}}{(d+2)|B^d|}
\left(\rho^{-1}+\varepsilon^2\rho^{-2}\left(\Lambda\rho-\frac{1}{2(d+2)}\Delta\rho\right)+O(\varepsilon^4)\right) \int_{S^{d-1}} \langle K_1,e_s\rangle\langle K_1,e_l\rangle f\rho dS_\theta\\
&\ +\frac{\varepsilon^{3}}{(d+3)|B^d|}
\left(\rho^{-1}+\varepsilon^2\rho^{-2}\left(\Lambda\rho-\frac{1}{2(d+2)}\Delta\rho\right)+O(\varepsilon^4)\right) \int_{S^{d-1}} \langle K_1,e_s\rangle\langle K_1,e_l\rangle (f\rho)_1dS_\theta\\
&\ +\frac{\varepsilon^{4}}{(d+4)|B^d|}\rho^{-1}\int_{S^{d-1}} \langle K_1,e_s\rangle\langle K_3,e_l\rangle f\rho+\langle K_1,e_l\rangle\langle K_3,e_s\rangle f\rho\\ 
&\hspace{3cm} +\langle K_1,e_s\rangle\langle K_1,e_l\rangle\left((f\rho)_2+ f\rho(d+4)H_1+f\rho R_1\right) dS_\theta\\
&\ +O(\varepsilon^5).
\end{align*}
Note that all the terms of odd order are zero. Moreover, the first term is also zero because 
$\int_{S^{d-1}}\langle \theta,e_s\rangle\langle \theta,e_l\rangle dS_\theta
= \sum_{m,j=1}^d\eta_s^m\eta_l^j\int_{S^{d-1}} \theta^m\theta^jdS_\theta
=\sum_{m=1}^d\eta_s^m\eta_l^m\int_{S^{d-1}} (\theta^m)^2dS_\theta
=\langle e_s,e_l\rangle |B^d|=0$. Indeed, it can be seen easily by using the symmetry of $S^{d-1}$. On the other hand, by using Lemma \ref{Aij}, we see that $$\int_{S^{d-1}}\langle K_1,e_s\rangle\langle K_3,e_l\rangle+\langle K_1,e_l\rangle\langle K_3,e_s\rangle dS_\theta
=\frac{|B^d|}{d+2}\left(\frac{2}{3} Ric(e_s,e_l)- {\bf H}\cdot{\rm I\!I}(e_s,e_l)\right),$$
hence
\begin{align*}
&\ -\hspace{-3.8mm}\int_{\widetilde B}\langle x-z,e_s\rangle\langle x- z,e_l\rangle f(x)d{\rm vol}_{\rho}\\
=&\ \frac{\varepsilon^{4}}{(d+4)(d+2)}\rho^{-1}\left(\frac{2}{3}Ric(e_s,e_l)- {\bf H}\cdot{\rm I\!I}(e_s,e_l)\right)f\rho\\
&\ +\frac{\varepsilon^{4}}{(d+4)|B^d|}\rho^{-1}\int_{S^{d-1}} \langle\theta,e_s\rangle\langle \theta,e_l\rangle\left(\frac{1}{2}\nabla_\theta^2 (f\rho)|_{z}+ \left(\frac{d+4}{24}|{\rm I\!I}(\theta,\theta)|^2-\frac{1}{6} Ric(\theta,\theta)\right)f\rho\right) dS_\theta+O(\varepsilon^6)\\
=&\  \frac{\varepsilon^{4}}{(d+4)(d+2)}\left(\frac{2}{3}Ric(e_s,e_l)- {\bf H}\cdot{\rm I\!I}(e_s,e_l)\right)f\\
&\ +\frac{\varepsilon^{4}}{(d+2)(d+4)}\left(\rho^{-1}\nabla_s\nabla_l (f\rho)|_{z}+ \left({\bf H}\cdot\frac{1}{2}{\rm I\!I}(e_s,e_l)-\frac{2}{3} Ric(e_s,e_l)\right)f\right)+O(\varepsilon^6)\\
=&\ \frac{\varepsilon^{4}}{(d+2)(d+4)}\left(\rho^{-1}\nabla_s\nabla_l (f\rho)|_{z}-\frac{1}{2}{\bf H}\cdot{\rm I\!I}(e_s,e_l)f\right)+O(\varepsilon^6).
\end{align*}
\\
(v)
When $\rho$ is a constant, it is easier to compute and one can obtain 
\begin{align*}
&\ -\hspace{-3.8mm}\int_{\widetilde B}\langle x-z,e_j \rangle\langle x-z,e_k \rangle^2 d{\rm vol}\\
=&\ \frac{1}{{\rm Vol}(\widetilde{B})}\int_{S^{d-1}}\int_0^{\tilde\varepsilon(\theta)}  \langle K_1t+K_3t^3+O(t^4),e_j\rangle\langle K_1t+K_3t^3+O(t^4),e_k\rangle^2\\
&\hspace{8cm} \cdot(t^{d-1}+t^{d+1}R_1+O(t^{d+2}))dtdS_\theta\\
=&\ \frac{1}{{\rm Vol}(\widetilde{B})}\int_{S^{d-1}}\frac{\varepsilon^{d+3}}{d+3}\langle \theta,e_j \rangle\langle \theta,e_k \rangle^2 + \frac{{\varepsilon}^{d+5}}{24}|{\rm I\!I}|^2\langle \theta,e_j \rangle\langle \theta,e_k \rangle^2 \\
&\hspace{14mm} +\frac{{\varepsilon}^{d+5}}{d+5}\left(\langle\frac{1}{6}\nabla_\theta{\rm I\!I},e_j\rangle\langle \theta,e_k \rangle^2
+2\langle \theta,e_j \rangle\langle \theta,e_k \rangle\langle\frac{1}{6}\nabla_\theta{\rm I\!I},e_k\rangle +\langle \theta,e_j \rangle\langle \theta,e_k \rangle^2 \frac{-1}{6}Ric(\theta,\theta)
\right)\\
&\hspace{14mm} +O(\varepsilon^{d+6})\ dS_{\theta}\\
=&\ O(\varepsilon^6),
\end{align*}
because $\int_{S^{d-1}}\langle \theta,e_j \rangle\langle \theta,e_l \rangle^2dS_\theta$, $\int_{S^{d-1}}|{\rm I\!I}|^2\langle \theta,e_j \rangle\langle \theta,e_l \rangle^2dS_\theta$, and other explicit terms written in the above integrals are odd integrals and thus vanish.

In general, when the density is  nonuniform, one has
\begin{align*}
&\ -\hspace{-3.8mm}\int_{\widetilde B}\langle x-z,e_j \rangle\langle x-z,e_k \rangle^2 d{\rm vol}_{\rho}\\
=&\ \frac{1}{{\rm Vol}_{\rho}(\widetilde{B})}\int_{S^{d-1}}\int_0^{\tilde\varepsilon(\theta)}  \langle K_1t+K_3t^3+O(t^4),e_j\rangle\langle K_1t+K_3t^3+O(t^4),e_k\rangle^2\\
&\hspace{5cm} \cdot(t^{d-1}+t^{d+1}R_1+O(t^{d+2}))(\rho(z)+\rho_1t+O(t^2))dtdS_\theta
\end{align*}
and now the term $\rho_1t=t\nabla_\theta\rho$ will produce an $O(\varepsilon^4)$-term 
\begin{align*}
&\ \frac{1}{{\rm Vol}_{\rho}(\widetilde{B})}\int_{S^{d-1}}\int_0^{\tilde\varepsilon(\theta)}  \langle K_1,e_j\rangle\langle K_1,e_k\rangle^2
t^{d+3}\nabla_\theta\rho\ dtdS_\theta\\
=&\ \frac{1}{{\rm Vol}_{\rho}(\widetilde{B})}\int_{S^{d-1}} \frac{\varepsilon^{d+4}}{d+4}\langle \theta,e_j\rangle\langle \theta,e_k\rangle^2
\nabla_\theta\rho\ dS_\theta\\
=&\ \frac{1}{{\rm Vol}_{\rho}(\widetilde{B})}\int_{S^{d-1}} \frac{\varepsilon^{d+4}}{d+4}\langle \theta,e_j\rangle\langle \theta,e_k\rangle^2
\langle\nabla\rho,\theta\rangle\ dS_\theta.
\end{align*}
As in (ii), using the divergence theorem and $\int_{B^d}x_k^2dx=\frac{1}{d+2}|B^d|$, one can show that this $O(\varepsilon^4)$-term is 
\[
\frac{1}{{\rm Vol}_{\rho}(\widetilde{B})}\frac{\varepsilon^{d+4}}{d+4}\int_{B^d}\left(\langle x,e_k\rangle^2+\langle x,e_j\rangle\langle x,e_k\rangle\right)\nabla_j\rho|_{z} dx=\beta\rho^{-1}\nabla_j\rho|_{z}
\] 
when $j\neq k$ and is 
\[
\frac{1}{{\rm Vol}_{\rho}(\widetilde{B})}\frac{\varepsilon^{d+4}}{d+4}\int_{B^d}3\langle x,e_j\rangle^2\nabla_j\rho|_{z} dx=3\beta\rho^{-1}\nabla_j\rho|_{z}
\]
when $j=k$. 
The second integral can be shown easily, so we skip the proof.

(vi) The proof is similar to previous cases. For the first integral, it is easy to see that the leading term is 
$$ \frac{1}{{\rm Vol}_{\rho}(\widetilde{B})}\int_{S^{d-1}}\int_0^{\tilde\varepsilon(\theta)}  \langle K_1,e_s\rangle^2\langle K_1,e_l\rangle^2 t^{d+3}\rho\ dtdS_\theta=
\Big\{
\begin{array}{rc}
 \beta+O(\varepsilon^6),& s\neq l \\
 3\beta+O(\varepsilon^6),    & s=l
\end{array}.
$$
Other integrals can be similarly derived.
}
\end{proof}

Since $\langle x-z,e_j\rangle=\langle x-z,u_j\rangle+O(\varepsilon^3)$, we can derive the following estimates.

\begin{Cor} \label{intfu}
Let $\{u_j\}_{j=1}^d$ be the orthonormal basis given by local PCA. Following Theorem \ref{intf}, we have
\\
{\rm (i)} 
For all $j=1,\dots,d$,
$$-\hspace{-3.8mm}\int_{\widetilde B}\langle x-z,u_j\rangle f(x)d{\rm vol}_{\rho}
=\alpha\left(\nabla_jf|_{z}+\rho^{-1}(z)f(z)\nabla_j\rho|_{z}\right)+O(\varepsilon^{3});
$$\\
{\rm (ii)} 
\begin{align*}
&\ -\hspace{-3.8mm}\int_{\widetilde B}\langle x-z,u_j\rangle^2f(x){\rm dvol}_{\rho}\\
=&\ \alpha f(z)+\beta\left(\nabla_j^2f|_{z}+\frac{1}{2}\Delta f(z)\right)\\
&\ +\beta \rho^{-1}(z)\left( 2\nabla_jf\nabla_j\rho|_{z}+ \langle \nabla f,\nabla\rho\rangle(z)+f(z)\nabla_j^2\rho|_{z}-\frac{1}{d+2}f(z)\Delta\rho(z)\right)\\
&\ +\beta \mathcal{U}_{jj}f(z)+O(\varepsilon^5), 
\end{align*}
where $\mathcal{U}_{jj}=2 \Lambda-\frac{1}{2}{\bf H}_{z}\cdot{\rm I\!I}_{z}(e_j,e_j)+D_{jj}$
and $\beta D_{jj}f(z)$ is the leading term of $-\hspace{-3.4mm}\int_{\widetilde B}\langle x-z,u_j\rangle^2f(x){\rm dvol}_{\rho}
--\hspace{-3.4mm}\int_{\widetilde B}\langle x-z,e_j\rangle^2f(x){\rm dvol}_{\rho}$;
\\
{\rm (iii)} 
For all $s,l=1,\dots,d$ and $s\neq l$,
\begin{align*}
&-\hspace{-3.8mm}\int_{\widetilde B}\langle x-z,u_s\rangle\langle x-z,u_l\rangle f(x)d{\rm vol}_{\rho}\\
=&\  \beta\rho^{-1}(z)\nabla_s\nabla_l(f\rho)|_{z} +\beta\mathcal{V}_{s,l}f(z)
+O(\varepsilon^6)\\
=&\  \beta\nabla_s\nabla_l f|_{z}
+\beta\rho^{-1}(z)\left(\nabla_sf\nabla_l\rho+\nabla_lf\nabla_s\rho+f\nabla_s\nabla_l\rho\right)|_{z}\\
& +\beta\mathcal{V}_{s,l}f(z)+O(\varepsilon^5),
\end{align*}
where $\mathcal{V}_{s,l}=-\frac{1}{2}{\bf H}_{z}\cdot{\rm I\!I}_{z}(e_s,e_l)+D_{s,l}
$
and
$\beta D_{s,l}f(z)$ is the leading term of $-\hspace{-3.2mm}\int_{\widetilde B}\langle x-z,u_s\rangle\langle x-z,u_l\rangle f(x)d{\rm vol}_{\rho}
--\hspace{-3.2mm}\int_{\widetilde B}\langle x-z,e_s\rangle\langle x-z,e_l\rangle f(x)d{\rm vol}_{\rho}$.
\\
{\rm (iv)} 
$$-\hspace{-3.8mm}\int_{\widetilde B}\langle x-z,u_s \rangle\langle x-z,u_l \rangle^2 d{\rm vol}_\rho=
\Big\{
\begin{array}{rc}
 \beta\rho^{-1}(z)\nabla_s\rho|_{z}+O(\varepsilon^5),& s\neq l \\
 3\beta\rho^{-1}(z)\nabla_s\rho|_{z}+O(\varepsilon^5),    & s=l
\end{array}.$$
 
$$-\hspace{-3.8mm}\int_{\widetilde B}\langle x-z,u_j \rangle\langle x-z,u_s \rangle\langle x-z,u_l \rangle d{\rm vol}_\rho=O(\varepsilon^5)\ {\mbox when }\ j,s,l \mbox{ are mutually distinct}.$$
\\
{\rm (v)} 
For all $s,l=1,\dots,d$,
$$
-\hspace{-3.8mm}\int_{\widetilde B}\langle x-z,u_s\rangle^2\langle x-z,u_l\rangle^2d{\rm vol}_{\rho}
=
\Big\{
\begin{array}{rc}
 \beta+O(\varepsilon^6),& s\neq l \\
 3\beta+O(\varepsilon^6),    & s=l
\end{array}.
$$
When $s,t,l,m$ are mutually distinct, 
$-\hspace{-3.2mm}\int_{\widetilde B}\langle x-z,e_s\rangle\langle x-z,e_l\rangle^3d{\rm vol}_{\rho}$,
$-\hspace{-3.2mm}\int_{\widetilde B}\langle x-z,e_s\rangle\langle x-z,e_t\rangle\langle x-z,e_l\rangle^2d{\rm vol}_{\rho}$
and 
$-\hspace{-3.2mm}\int_{\widetilde B}\langle x-z,e_s\rangle\langle x-z,e_t\rangle\langle x-z,e_l\rangle\langle x-z,e_m\rangle d{\rm vol}_{\rho}$ are all of $O(\varepsilon^6)$.

\end{Cor}

\section{Proof of the main theorem}\label{section: proof of main theorem heuristic approach}
The proof of the main theorem is divided into several steps, which we detail below. In the following, we use dashes to indicate the lower triangular portion of a symmetric matrix.

\subsection{Step 1:}

Recall that $\{x_{z,j}\}_{j=1}^{k_z}$ are $k_z$ neighbors of $z$ and ${\bf q}_j$ is the projection of $x_{z,j}$ onto the subspace $V_z={\rm span}\{u_s\}_{s=1}^d$ that is derived from the local PCA. Thus, $({\bf q}_j)_s= \langle \iota(x_{z,j})-\iota(z),u_s\rangle$. Denote $\mathcal{B}=\iota(M)\cap B_{\varepsilon}^p(z)$ and $d{\rm vol}_\rho=\rho(x)d{\rm vol}$ as before. Also, denote $J$ to be the 1-1 map that indexing the upper triangular matrix of size $d\times d$ by $1,\dots,\frac{d(d-1)}{2}$.

To study the asymptotic behavior of $\frac{1}{k_z}Z^TZ$, we divide it into two steps. The first step is evaluating its bias from the desired asymptotic quantity, and the second step is evaluating its large deviation from the desired quantity.
First of all, under the manifold assumption and the law of large number, asymptotically when $n\to\infty$, we expect that $\frac{1}{k_z}Z^TZ$ converges to
$$
L=\left[\begin{array}{cccc}
1&L_{AB}& L_{AC}& L_{AD}\\
- &L_{BB}&L_{BC}&L_{BD}\\ 
- & - & L_{CC}&L_{CD}\\
- & - &- &L_{DD}
\end{array}\right]\in \mathbb{R}^{(1+d+d+\frac{d(d-1)}{2})\times (1+d+d+\frac{d(d-1)}{2})},
$$
where $L_{AB}$ is a $d$-dim vector
$$\left[
\begin{array}{ccc}
     -\hspace{-3mm}\int_{\mathcal{B}} \langle \iota(x)-\iota(z),u_1\rangle d{\rm vol}_\rho& \cdots& -\hspace{-3mm}\int_{\mathcal{B}} \langle \iota(x)-\iota(z),u_d\rangle d{\rm vol}_\rho
\end{array}\right]\,,
$$
$L_{AC}$ is a $d$-dim vector
$$
\left[
\begin{array}{ccc}
     -\hspace{-3mm}\int_{\mathcal{B}} \langle \iota(x)-\iota(z),u_1\rangle^2d{\rm vol}_\rho& \cdots& -\hspace{-3mm}\int_{\mathcal{B}} \langle \iota(x)-\iota(z),u_d\rangle^2d{\rm vol}_\rho
\end{array}\right]\,,
$$
$L_{AD}$ is a $\frac{d(d-1)}{2}$-vector 
$$\left[
\begin{array}{ccc}
-\hspace{-3mm}\int_{\mathcal{B}} \langle\iota(x)-\iota(z),u_1\rangle\langle\iota(x)-\iota(z),u_2\rangle d{\rm vol}_\rho 
&\cdots
&-\hspace{-3mm}\int_{\mathcal{B}} \langle\iota(x)-\iota(z),u_{d-1}\rangle\langle\iota(x)-\iota(z),u_d\rangle d{\rm vol}_\rho \end{array}\right],
$$
$L_{BB}$ is a $d\times d$ diagonal matrix 
$$\left[
\begin{array}{ccc}
-\hspace{-3mm}\int_{\mathcal{B}} \langle\iota(x)-\iota(z),u_1\rangle^2d{\rm vol}_\rho&&0\\
&\ddots&\\
0&&-\hspace{-3mm}\int_{\mathcal{B}} \langle\iota(x)-\iota(z),u_d\rangle^2d{\rm vol}_\rho
\end{array}\right],
$$
$L_{BC}$ is a $d\times d$ matrix with the $(s,t)$-th entry, $s,t=1,\ldots,d$,
$$
-\hspace{-3.6mm}\int_{\mathcal{B}} \langle\iota(x)-\iota(z),u_s\rangle \langle\iota(x)-\iota(z),u_t\rangle^2d{\rm vol}_\rho \,,
$$
$L_{BD}$ is a $d\times \frac{d(d-1)}{2}$ matrix with the $(s,J(k,l))$-th entry, $s=1,\ldots,d$, $k,l=1\ldots,d$ and $k\neq l$,
$$
-\hspace{-3.6mm}\int_{\mathcal{B}} \langle\iota(x)-\iota(z),u_s\rangle\langle\iota(x)-\iota(z),u_k\rangle\langle\iota(x)-\iota(z),u_l\rangle d{\rm vol}_\rho\,,
$$
$L_{CC}$ is the $d\times d$ matrix with the $(s,t)$-th entry, $s,t=1,\ldots,d$,
$$
-\hspace{-3.6mm}\int_{\mathcal{B}} \langle\iota(x)-\iota(z),u_s\rangle^2 \langle\iota(x)-\iota(z),u_t\rangle^2d{\rm vol}_\rho \,,
$$
$L_{CD}$ is the $d\times \frac{d(d-1)}{2}$ matrix with the $(s,J(k,l))$-th entry, $s=1,\ldots,d$, $k,l=1\ldots,d$ and $k\neq l$,
$$
-\hspace{-3.6mm}\int_{\mathcal{B}} \langle\iota(x)-\iota(z),u_s\rangle^2\langle\iota(x)-\iota(z),u_k\rangle\langle\iota(x)-\iota(z),u_l\rangle d{\rm vol}_\rho\,,
$$ and
$L_{DD}$ is the $\frac{d(d-1)}{2}\times \frac{d(d-1)}{2}$ matrix with the $(J(s,t),J(k,l))$-th entry, $s,t=1,\ldots,d$, $s\neq t$, $k,l=1\ldots,d$ and $k\neq l$,
$$
-\hspace{-3.6mm}\int_{\mathcal{B}} \langle\iota(x)-\iota(z),u_s\rangle\langle\iota(x)-\iota(z),u_t\rangle\langle\iota(x)-\iota(z),u_k\rangle\langle\iota(x)-\iota(z),u_l\rangle d{\rm vol}_\rho\,.
$$ 

Note that the difference between $\frac{1}{k_z}Z^TZ$ and $L$ is a large deviation term and we denote it by $\tt LD$. We will compute $\tt LD$ precisely in the later sections. Now we deal with $L$ first.
When $\mathcal{B}$ is close to $B_{\varepsilon,\delta}^d$, we can show in the following lemma that this seeming complicated matrix $L$ can be well approximated by the  leading matrix

\begin{align}\label{L0}
L^0
=\left[\begin{array}{c|cccc|ccccc|c}
1 & 0& \cdots& 0 &\gamma_1  &\alpha_2 &\cdots& \cdots&\alpha_2& \alpha_1 &0\\ \hline
- &\alpha_2 &0&\cdots&0 & 0&\cdots &\cdots&\cdots&0 &\\ 
- & - &\ddots &\ddots&\vdots &\vdots&\ddots&\ddots &\ddots &\vdots&\\ 
- & - &- &\alpha_2& 0&0& \cdots&\cdots&\cdots& 0& L^0_{BD}\\ 
- & - &- &- &\alpha_1  &\mu_2&\cdots&\cdots&\mu_2&\mu_1& \\ \hline 
- & - &-&-&-& \beta_3&\beta_4&\cdots&\beta_4&\beta_2&0 \\
- & - &-&-&-& -    &\ddots&\ddots&\vdots&\vdots&\vdots\\
- & - &-&-&-& -    &-     &\ddots&\beta_4&\vdots&0\\
- & - &-&-&-& -    &-     &-&\beta_3&\beta_2&0\\
- & - &-&-&-& -    &-     &-     &\cdots&\beta_1&0\\ \hline 
- & - &-&-&-& -    &-     &-     &0&0&L^0_{DD} \end{array}\right],
\end{align}
where 
$L^0_{BD}$
is
$\left[
\begin{array}{cccc|cccc|c|ccc|cc|c}
0&\cdots&0&\mu_2 &   0&\cdots& 0& 0    &  \cdots  &  0&0&0&      0&0      &   0 \\   
0&\cdots&0& 0    &   0&\cdots& 0&\mu_2 &  \cdots  &  0&0&0&      0&0      &   0 \\   
 &      & &      &    &      &  &      &  \ddots  &   & & &       &       &    \\   
0&\cdots&0& 0    &   0&\cdots& 0& 0    &  \cdots  &  0&0&\mu_2&  0&0      &   0 \\
0&\cdots&0& 0    &   0&\cdots& 0& 0    &  \cdots  &  0&0&0&      0&\mu_2  &   0 \\
0&\cdots&0& 0    &   0&\cdots& 0& 0    &  \cdots  &  0&0&0&      0&0      &\mu_2\\
0&\cdots&0& 0    &   0&\cdots& 0& 0    &  \cdots  &  0&0&0&      0&0      &   0 \\   
\end{array}
\right]_{d\times \frac{d(d-1)}{2}}
$
and 
$L^0_{DD}$
is the diagonal matrix
$\left[
\begin{array}{cccc|c|ccc|cc|c}
\beta_4 & &  &  &&    & & &     &    & \\
  & \ddots&  &  &&    & & &     &    & \\
  & & \beta_4&  &&    & & &     &    & \\
  & &     &\beta_2 &&    & & &     &            & \\ \hline
  & &     &&\ddots&  & & &     &            & \\ \hline
  & &     &     &&  \beta_4  & & &     &        & \\
  & &     &     &&    &\beta_4& &     &          & \\
  & &     &     &&    &      &\beta_2&     &    & \\ \hline
  & &     &     &&    &      &       &\beta_4&         &\\
  & &     &     &&    &      &       &       & \beta_2& \\ \hline
  & &     &     &&    &      &       &   &        & \beta_2
\end{array}
\right]_{\frac{d(d-1)}{2}\times \frac{d(d-1)}{2}}
$.

In the following, we usually omit the notation $\iota$ and thus denote $\iota(M)$, $\iota(x)$ and $\iota(z)$ by $M,x$ and $z$, respectively.

\begin{Lma}\label{C.W bound of L and L0}
Denote $\frac{1}{k_z}Z^TZ=L+{\tt LD}=L^0+W+{\tt LD}$. 
Then
$$W:=L-L^0=
\left[\begin{array}{cccc}
0 & O(\varepsilon^2) & O(\varepsilon^3) &O(\varepsilon^3) \\
- & O(\varepsilon^3) & O(\varepsilon^4) &O(\varepsilon^4) \\
    - & - & O(\varepsilon^5) &O(\varepsilon^5) \\
    - & - & - &O(\varepsilon^5)
\end{array}\right]
.$$ 
\end{Lma}

\begin{proof} This Lemma comes immediately by applying Lemma \ref{lemma of local poly integration}.
\end{proof}

Now we have 
$$\frac{1}{k_z}Z^TZ=L+{\tt LD}=L^0+W+{\tt LD},$$
where $L^0$ has an explicit expression and the order of $W$ is also known by the above lemma. To obtain the Hessian estimator, 
one can compute $(\frac{1}{k_z}Z^TZ)^{-1}(\frac{1}{k_z}Z^T{\bf f})$ directly. However, it would be painful and time-consuming. The task can be accomplished in a much better way as follows. 

Let ${\bf f}_{Taylor}$ be the vector 
\begin{equation}\label{definition: fTaylor}
{\bf f}_{Taylor}:=\left[
\begin{array}{cccc}
  f(z)
& \nabla f(z)_{1\times d}
& \frac{1}{2}(h_{ss})_{1\leq s\leq d} 
& (h_{st})_{1\leq s< t\leq d}
\end{array}\right]_{1\times(1+d+d+\frac{d(d-1)}{2})}^T,
\end{equation}
where $(h_{st})=\left(\frac{\partial^2f}{\partial x_s\partial x_t}\big|_z\right)\in \mathbb{R}^{d\times d}$ is the Hessian of $f$ at $z$.
Observing that $\frac{1}{k_z}Z^T{\bf f}$ is actually very similar to $\frac{1}{k_z}Z^TZ {\bf f}_{Taylor}$, we can 
denote $$\frac{1}{k_z}Z^T{\bf f}=\frac{1}{k_z}Z^TZ {\bf f}_{Taylor}-V$$ and derive
$$\left(\frac{1}{k_z}Z^TZ\right)^{-1}\left(\frac{1}{k_z}Z^T{\bf f}\right)
=\left(\frac{1}{k_z}Z^TZ\right)^{-1}\left(\frac{1}{k_z}Z^TZ {\bf f}_{Taylor}-V\right)={\bf f}_{Taylor}-\left(\frac{1}{k_z}Z^TZ\right)^{-1}V.$$

Therefore we reduce the computation from $\left(\frac{1}{k_z}Z^TZ\right)^{-1}\frac{1}{k_z}Z^T{\bf f}$ to $\left(\frac{1}{k_z}Z^TZ\right)^{-1}V$ and makes the higher order estimation doable. Indeed, if we compute $\left(\frac{1}{k_z}Z^TZ\right)^{-1}$ and $\frac{1}{k_z}Z^T{\bf f}$ directly and multiply them together, then we can also prove our result for interior points of $M$, as shown in the next section, Appendix D. For points near the boundary, direct approach might work if someone can complete the tremendous computation which we do not think is realistic. Another advantage of considering $V$ is that we can get rid of the demanding regularity assumption of $\rho$ and $f$
(cf. Theorem \ref{ThmInterior}) and explain why the non-uniform distribution only affects higher order errors.

\begin{Thm}\label{V} For any $C^{2,\kappa}$-function $f$ with $\kappa\in(0,1]$,
$$V:=\frac{1}{k_z}Z^TZ {\bf f}_{Taylor}-\frac{1}{k_z}Z^T{\bf f}
=\left[
\begin{array}{c}
O(\varepsilon^{2+\min(\tau,\kappa)})  \\
O(\varepsilon^{3+\min(\tau,\kappa)}) \\
O(\varepsilon^{4+\min(\tau,\kappa)}) \\
O(\varepsilon^{4+\min(\tau,\kappa)}) 
\end{array}\right]+{\tt LD}.$$
\end{Thm}

\begin{proof}
Let's compare the first entries in 
$\frac{1}{k_z}Z^TZ {\bf f}_{Taylor}$ and $\frac{1}{k_z}Z^T{\bf f}$.
Recall that $\frac{1}{k_z}Z^TZ=L+{\tt LD}$, 
where $\tt LD$ contains the deviation terms depending on the relationship between $n$ and $\varepsilon$. Note that ${\tt LD}$
might be different in each equation. As before, we omit the notation $\iota$ in the following computations.
Let 
{\allowdisplaybreaks
$$
L{\bf f}_{Taylor}=\left[\begin{array}{cccc}
1&L_{AB}& L_{AC}& L_{AD}\\
- &L_{BB}&L_{BC}&L_{BD}\\ 
- & - & L_{CC}&L_{CD}\\
- & - &- &L_{DD}
\end{array}\right]
\left[\begin{array}{c}
  f(z)\\
 \nabla f(z)_{1\times d}\\
 \frac{1}{2}(h_{ss})_{1\leq s\leq d}\\ 
 (h_{st})_{1\leq s< t\leq d}
\end{array}\right].
$$
The first entry of $L{\bf f}_{Taylor}$ is  
\begin{align*}
&\ f(z)
+\left[-\hspace{-3.8mm}\int_{\mathcal B} \langle x-z,u_j\rangle d{\rm vol}_\rho\right]_j \cdot \nabla f\big|_z
 +\left[-\hspace{-3.8mm}\int_{\mathcal B} \langle x-z,u_j\rangle\langle x-z,u_k\rangle d{\rm vol}_\rho\right]_{j,k} \cdot \nabla^2 f\big|_z\\
=&\ f(z)
+\sum_{j=1}^d\left(-\hspace{-3.8mm}\int_{\mathcal B} t(\langle K_1,e_j\rangle+O(\varepsilon^{\tau+1}))
+t^3(\langle K_3,e_j\rangle+O(\varepsilon^{\tau+1}))+O(t^4) d{\rm vol}_\rho\right) \partial_j f\big|_z\\
&\ +\sum_{j,k=1}^d\left(-\hspace{-3.8mm}\int_{\mathcal B} t^2(\langle K_1,e_j\rangle \langle K_1,e_k\rangle+O(\varepsilon^{\tau+1})) +t^4 (\langle K_3,e_j\rangle \langle K_1,e_k\rangle+O(\varepsilon^{\tau+1})) +O(t^5)  d{\rm vol}_\rho\right)\frac{1}{2} \partial_j\partial_k f\\
=&\ f(z)
+\left(\sum_{j=1}^d\left(-\hspace{-3.8mm}\int_{\mathcal B} t\bar\theta^j
+t^3\langle K_3,e_j\rangle d{\rm vol}_\rho\right) 
+O(\varepsilon^{\tau+2})\right)\partial_j f\big|_z\\
&\ +\left(\sum_{j,k=1}^d\left(-\hspace{-3.8mm}\int_{\mathcal B} t^2\bar\theta^j\bar\theta^k +t^4 \langle K_3,e_j\rangle \langle K_1,e_k\rangle  d{\rm vol}_\rho\right)+O(\varepsilon^{\tau+3})\right)\frac{1}{2} \partial_j\partial_k f.
\end{align*}

On the other hand, the first entry of $\frac{1}{k_z}Z^T{\bf f}$ is
\begin{align*}
&\
-\hspace{-3.8mm}\int_{\mathcal B} \left(f(z)+t\nabla_\theta f\big|_z+\frac{1}{2}t^2\nabla_\theta^2f\big|_z +O(t^{2+\kappa})\right)  d{\rm vol}_\rho +{\tt LD}\\
=&\
f(z)+-\hspace{-3.8mm}\int_{\mathcal B}t\nabla_\theta f d{\rm vol}_\rho 
+-\hspace{-3.8mm}\int_{\mathcal B}\frac{1}{2}t^2\nabla_\theta^2f  d{\rm vol}_\rho  
+O(\varepsilon^{2+\kappa}) +{\tt LD}\\
=&\
f(z)+\sum_{j=1}^d\left(-\hspace{-3.8mm}\int_{\mathcal B}t\bar\theta^j d{\rm vol}_\rho\right) \partial_j f\big|_z
+\sum_{j,k}\left(-\hspace{-3.8mm}\int_{\mathcal B}t^2 \bar\theta^j\bar\theta^k d{\rm vol}_\rho \right)\frac{1}{2}\partial_j\partial_k f\big|_z +O(\varepsilon^{2+\kappa}) +{\tt LD}.
\end{align*}

It is easy to see that the first entries of $\frac{1}{k_z}Z^TZ {\bf f}_{Taylor}$ and $\frac{1}{k_z}Z^T{\bf f}$ coincide with each other in the leading terms and the difference is 
\begin{align*}
V_A
=&\ \sum_{j}\left(-\hspace{-3.8mm}\int_{\mathcal B} 
t^3\langle K_3,e_j\rangle d{\rm vol}_\rho+O(\varepsilon^{\tau+2})\right) \partial_j f\big|_z
+O(\varepsilon^{2+\kappa}) -{\tt LD}.
\end{align*}

Similarly, one can derive
\begin{align*}
[V_B]_j
=&\ 
\left(-\hspace{-3.8mm}\int_{\mathcal B} \langle x-z,u_j\rangle d{\rm vol}_\rho\right) f(z)
+\sum_{k=1}^d\left(-\hspace{-3.8mm}\int_{\mathcal B} \langle x-z,u_j\rangle \langle x-z,u_k\rangle d{\rm vol}_\rho\right) \partial_k f\big|_z\\
&\ +\sum_{k,l}\left(-\hspace{-3.8mm}\int_{\mathcal B} \langle x-z,u_j\rangle \langle x-z,u_k\rangle \langle x-z,u_l\rangle d{\rm vol}_\rho\right) \frac{1}{2} \partial_k\partial_l f\big|_z\\
&\ - -\hspace{-3.8mm}\int_{\mathcal B} \langle x-z,u_j\rangle \left(f(z)+t\nabla_\theta f\big|_z+\frac{1}{2}t^2\nabla_\theta^2 f\big|_z+O(t^{2+\kappa})\right)  d{\rm vol}_\rho -{\tt LD}
\\
=&\ 
\sum_{k=1}^d\left(-\hspace{-3.8mm}\int_{\mathcal B} t^2(\bar\theta^j +O(\varepsilon^{\tau+1}))(\bar\theta^k +O(\varepsilon^{\tau+1}))+O(t^4) d{\rm vol}_\rho\right)  \partial_k f\big|_z\\
&\ +\sum_{k,l}\left(-\hspace{-3.8mm}\int_{\mathcal B} t^3(\bar\theta^j +O(\varepsilon^{\tau+1}))(\bar\theta^k +O(\varepsilon^{\tau+1}))(\bar\theta^l +O(\varepsilon^{\tau+1}))
d{\rm vol}_\rho\right) \frac{1}{2} \partial_k\partial_l f\big|_z\\
&\ - \sum_{k=1}^d\left(-\hspace{-3.8mm}\int_{\mathcal B} \left( t^2(\bar\theta^j +O(\varepsilon^{\tau+1}))\bar\theta^k + t^4(\langle K_3,e_j\rangle+O(\varepsilon^{\tau+1}) )\bar\theta^k \right) d{\rm vol}_\rho\right) \partial_kf\big|_z\\
&\ - \sum_{k,l}\left(-\hspace{-3.8mm}\int_{\mathcal B}t^3(\bar\theta^j +O(\varepsilon^{\tau+1}))\bar\theta^k\bar\theta^l d{\rm vol}_\rho\right) \frac{1}{2}\partial_k\partial_l f\big|_z  +O(\varepsilon^{3+\kappa})-{\tt LD}\\
=&\ 
\sum_{k=1}^d\left(-\hspace{-3.8mm}\int_{\mathcal B} t^4\bar\theta^j\langle K_3,u_k\rangle d{\rm vol}_\rho +O(\varepsilon^{\tau+3})\right) \partial_k f\big|_z\\
&\ +\sum_{k,l}O(\varepsilon^{\tau+4})\frac{1}{2} \partial_k\partial_l f\big|_z
+O(\varepsilon^{3+\kappa})-{\tt LD},
\end{align*}

\begin{align*}
[V_C]_j
=&\ 
\left(-\hspace{-3.8mm}\int_{\mathcal B} \langle x-z,u_j\rangle^2 d{\rm vol}_\rho\right) f(z)
+\sum_{k=1}^d\left(-\hspace{-3.8mm}\int_{\mathcal B} \langle x-z,u_j\rangle^2 \langle x-z,u_k\rangle d{\rm vol}_\rho\right) \partial_k f\big|_z\\
&\ +\sum_{k,l}\left(-\hspace{-3.8mm}\int_{\mathcal B} \langle x-z,u_j\rangle^2 \langle x-z,u_k\rangle \langle x-z,u_l\rangle d{\rm vol}_\rho\right)\frac{1}{2} \partial_k\partial_l f\big|_z\\
&\ - -\hspace{-3.8mm}\int_{\mathcal B} \langle x-z,u_j\rangle^2 \left(f(z)+t\nabla_\theta f\big|_z+\frac{1}{2}t^2\nabla_\theta^2 f\big|_z+O(t^{2+\kappa})\right)  d{\rm vol}_\rho -{\tt LD}\\
=&\ 
\sum_{k=1}^d\left(-\hspace{-3.8mm}\int_{\mathcal B}t^5(\bar\theta^j)^2\langle K_3,u_k\rangle d{\rm vol}_\rho +O(\varepsilon^{\tau+4})\right)\partial_k f\big|_z 
+O(\varepsilon^{4+\kappa})-{\tt LD},
\end{align*}
}
and
{\allowdisplaybreaks
\begin{align*}
[V_D]_{ij}  
=&\ 
\left(-\hspace{-3.8mm}\int_{\mathcal B} \langle x-z,u_i\rangle\langle x-z,u_j\rangle d{\rm vol}_\rho\right) f(z)\\
&\ 
+\sum_{k=1}^d\left(-\hspace{-3.8mm}\int_{\mathcal B} \langle x-z,u_i\rangle\langle x-z,u_j\rangle \langle x-z,u_k\rangle d{\rm vol}_\rho\right) \partial_k f\big|_z\\
&\ +\sum_{k,l}\left(-\hspace{-3.8mm}\int_{\mathcal B}\langle x-z,u_i\rangle\langle x-z,u_j\rangle \langle x-z,u_k\rangle \langle x-z,u_l\rangle d{\rm vol}_\rho\right)\frac{1}{2} \partial_k\partial_l f\big|_z\\
&\ - -\hspace{-3.8mm}\int_{\mathcal B} \langle x-z,u_i\rangle\langle x-z,u_j\rangle \left(f(z)+t\nabla_\theta f\big|_z+\frac{1}{2}t^2\nabla_\theta^2 f\big|_z+O(t^{2+\kappa})\right)  d{\rm vol}_\rho -{\tt LD}\\
=&\ \sum_{k=1}^d\left(-\hspace{-3.8mm}\int_{\mathcal B} 
t^5\bar\theta^i\bar\theta^j\langle K_3,u_k\rangle d{\rm vol}_\rho  +O(\varepsilon^{\tau+4})\right) \partial_k f\big|_z
+O(\varepsilon^{4+\kappa}) -{\tt LD}.
\end{align*}
}
Therefore, 
$$V:
=\left[
\begin{array}{c}
V_A  \\
V_B \\
V_C \\
V_D
\end{array}\right]
=\left[
\begin{array}{c}
O(\varepsilon^{2+\min(\tau,\kappa)})  \\
O(\varepsilon^{3+\min(\tau,\kappa)}) \\
O(\varepsilon^{4+\min(\tau,\kappa)}) \\
O(\varepsilon^{4+\min(\tau,\kappa)}) 
\end{array}\right]+{\tt LD}.$$

\end{proof}

\subsection{Step 2: The large deviation estimate for $\frac{1}{k_z}Z^TZ$ and $V$} 
Next, we study the large deviation of $\frac{1}{k_z}Z^TZ$ caused by the random samples. Write as before
\begin{align}\label{approximation of ZtZ/t}
\frac{1}{k_z}Z^TZ=L+{\tt LD}=L^0+W+{\tt LD},
\end{align}
where ${\tt LD}$ indicates the large deviation.
We will divide ${\tt LD}$ into blocks 
\begin{align}\label{decomposition block LD}
{\tt LD}=\left[\begin{array}{cccc}
{\tt LD}_{AA}&{\tt LD}_{AB}& {\tt LD}_{AC}& {\tt LD}_{AD}\\
- &{\tt LD}_{BB}&{\tt LD}_{BC}& {\tt LD}_{BD}\\ 
- & - & {\tt LD}_{CC}&{\tt LD}_{CD}\\
- & - &- &{\tt LD}_{DD}
\end{array}\right]\in \mathbb{R}^{(1+d+d+\frac{d(d-1)}{2})\times (1+d+d+\frac{d(d-1)}{2})}\,,
\end{align}
where ${\tt LD}_{AA}\in \mathbb{R}$, ${\tt LD}_{AB},{\tt LD}_{AC}\in \mathbb{R}^{1\times d}$, ${\tt LD}_{AD}\in \mathbb{R}^{1\times d(d-1)/2}$, ${\tt LD}_{BB},{\tt LD}_{BC},{\tt LD}_{CC}\in \mathbb{R}^{d\times d}$, ${\tt LD}_{BD},{\tt LD}_{CD}\in \mathbb{R}^{d\times d(d-1)/2}$, and ${\tt LD}_{DD}\in \mathbb{R}^{d(d-1)/2\times d(d-1)/2}$,
and bound each term separately with high probability using the Bernstein's inequality. 

We first provide detailed computation of ${\tt LD}_{AD}\in \mathbb{R}^{1\times d(d-1)/2}$. Recall the notations $Z$, ${\bf y}$ and ${\bf q}$ from \eqref{formula: Zhat}. Fix $s,l=1,\ldots,d$ and $s\neq l$, and denote 
\begin{align}\label{Definition Ymsl}
Y_m^{(s,l)}:= ({\bf y}_s)_m({\bf y}_l)_m=({\bf q}_m)_s({\bf q}_m)_l\,,
\end{align}
for $m=1,\ldots,k_z$, where 
$\{{\bf y}_s\}_{s=1}^d$ are the column vectors given by local PCA and $\{{\bf q}_m\}_{m=1}^{k_z}$ are row vectors of the $(k_z\times d)$-matrix $\left[{\bf y}_1\ \cdots\ {\bf y}_d\right]$. 
Denote 
\begin{align}\label{Definition Qmsl}
Q^{(s,l)}:=\langle X-z,u_s\rangle \langle X-z,u_l\rangle \ \
\mbox{and}\ \ Q_i^{(s,l)}=\langle x_i-z,u_s\rangle \langle x_i-z,u_l\rangle=({\bf q}_i)_s({\bf q}_i)_l\,.
\end{align}
When $n\to\infty$, by the law of large number we would expect that the following term is close to its mean 
\begin{align*}
\frac{1}{k_z}\sum_{m=1}^{k_z}Y_m^{(s,l)}
=\frac{\frac{1}{n}\sum_{m=1}^{k_z}Y_m^{(s,l)}}{\frac{1}{n}\sum_{m=1}^{k_z}1}
=\frac{\frac{1}{n}\sum_{i=1}^{n}Q_i^{(s,l)}\varepsilon^{-d}\chi_{ B_\varepsilon(z)}(x_i)}{\frac{1}{n}\sum_{i=1}^{n}\varepsilon^{-d}\chi_{ B_\varepsilon(z)}(x_i)}\,. 
\end{align*}
Note that $\varepsilon^{-d}$ in front of $\chi_{ B_\varepsilon(z)}(X_i)$ is a normalization so that $\varepsilon^{-d}\chi_{ B_\varepsilon(z)}(X_i)$ is a proper kernel function for the approximation of identity purpose. Also note that   $Q_i^{(s,l)}\varepsilon^{-d}\chi_{ B_\varepsilon(z)}(X_i)$ is an i.i.d. random vector.
We now apply Bernstein's inequality in the denominator and numerator separately to achieve our goal.

\begin{Prop}[Bernstein's inequality] Suppose $X_1,\ldots,X_n$ are i.i.d. sampled from a random variable $X$, where $E(X)=0$, $|X|\leq a$ and ${\tt Var}(X)<\infty$. Then for $b>0$ we have 
$$
\mathbb{P}\left\{\left|\frac{1}{n}\sum_{i=1}^{n} X_i\right|> b\right\}\leq 2\exp\left( -\frac{n b^2}{2 {\tt Var}(X)+\frac{2}{3}a b}\right)\,.
$$
\end{Prop}

Now, denote 
\begin{align}\label{Definition Xmsl}
X_i^{(s,l)}:=Q_i^{(s,l)}\varepsilon^{-d}\chi_{ B_\varepsilon(z)}(X_i)-E\left(Q_i^{(s,l)}\varepsilon^{-d}\chi_{ B_\varepsilon(z)}(X_i)\right)\,, 
\end{align}
where $i=1,\ldots,n$ and $s\neq l$.
Note that $E(X_i^{(s,l)})=0$, $|X_i^{(s,l)}|\leq \varepsilon^{-d}\|{\bf q}\|^2\leq a=\Theta(\varepsilon^{-d+2})$ and 
\begin{align*}
&{\tt Var}(X_i^{(s,l)})\\
=&\ \varepsilon^{-2d}\int_{\mathcal{B}} {\bf q}_s^2{\bf q}_l^2  d{\rm vol}_\rho
-\varepsilon^{-2d}\left(\int_{\mathcal{B}} {\bf q}_s{\bf q}_l  d{\rm vol}_\rho\right)^2\\
=&\ \varepsilon^{-2d}\int_{\mathcal{B}} \langle x-z,u_s\rangle^2\langle x-z,u_l\rangle^2  d{\rm vol}_\rho
-\varepsilon^{-2d}\left(\int_{\mathcal{B}} \langle x-z,u_s\rangle\langle x-z,u_l\rangle d{\rm vol}_\rho\right)^2\\
=&\ \varepsilon^{-2d}\int_{\mathcal{B}} \langle x-z,e_s\rangle^2\langle x-z,e_l\rangle^2+O(\varepsilon^{\tau+5})  d{\rm vol}_\rho\\
&-\varepsilon^{-2d}\left(\int_{\mathcal{B}} \langle x-z,e_s\rangle\langle x-z,e_l\rangle +O(\varepsilon^{\tau+3})d{\rm vol}_\rho\right)^2\\
=&\ \varepsilon^{-d}\rho(z)C_{0,2,2}\varepsilon^4 +O(\varepsilon^{-d+5})=\Theta(\varepsilon^{-d+4}).
\end{align*}
Clearly, $a$ is much larger than the variance, so we would apply Bernstein's inequality to obtain a sharp bound.
Since we are interested in the Hessian, which is of order $\varepsilon^2$, we focus on the region $b\ll \varepsilon^2$ when we apply Bernstein's inequality; that is, the deviation caused by the random sampling should be much smaller than $\varepsilon^2$, and we also hope that this bound happens with a high probability. To this goal, note that we have
\begin{align*}
\frac{n b^2}{2 {\tt Var}(X_i^{(s,l)})+\frac{2}{3}a b}\geq \frac{n b^2}{4 {\tt Var}(X_i^{(s,l)})}\,,
\end{align*}
and we set
\begin{align*}
\frac{n b^2}{4 {\tt Var}(X_i^{(s,l)})}= 3\log(n)
\end{align*}
since we hope the probability that the large deviation happens to be controlled from above by $n^{-3}$. That is, we take
\begin{align}
b= \Theta\left(\sqrt{\frac{\log(n)}{n\varepsilon^{d-4}}}\right)\,,
\end{align}
and achieve the bound
\[
\mathbb{P}\left(\left|\frac{1}{n}\sum_{i=1}^{n}Q_i^{(s,l)}\varepsilon^{-d}\chi_{ B_\varepsilon(z)}(x_i)-
\int_M Q^{(s,l)} \varepsilon^{-d}\chi_{B_{\varepsilon}(z)}(x) d{\rm vol}_\rho
\right|>\sqrt{\frac{\log(n)}{n\varepsilon^{d-4}}}\right)\leq Cn^{-3}
\]
for some $C>0$. For a sanity check, note that since we have assumed $\frac{n\varepsilon^d}{\log(n)}\to \infty$, $\sqrt{\frac{\log(n)}{n\varepsilon^{d-4}}}$ is dominated by $\varepsilon^2$ when $n$ is sufficiently large. 
Similarly, for the denominator we have
\[
\mathbb{P}\left(\left|\frac{1}{n}\sum_{i=1}^{n}\varepsilon^{-d}\chi_{ B_\varepsilon(z)}(x_i)-\int_M \varepsilon^{-d} \chi_{B_{\varepsilon}(z)}(x)d{\rm vol}_\rho\right|>\sqrt{\frac{\log(n)}{n\varepsilon^d}}\right)\leq Cn^{-3}\,.
\]
Recall that $\int_{\mathcal{B}}\varepsilon^{-d}d{\rm vol}_{\rho}=\rho(z)C_{0,0}+O(\varepsilon)$. Thus, we have
\begin{align*}
 \frac{1}{k_z}\sum_{m=1}^{k_z}Y_m^{(s,l)} 
=&\ \frac{\frac{1}{n}\sum_{i=1}^{n}Q_i^{(s,l)}\varepsilon^{-d}\chi_{ B_\varepsilon(z)}(x_i)}{\frac{1}{n}\sum_{i=1}^{n}\varepsilon^{-d}\chi_{ B_\varepsilon(z)}(x_i)}\\
=&\ \frac{1}{\rho(z)C_{0,0}+O(\varepsilon)}\int_M Q^{(s,l)} \varepsilon^{-d}\chi_{B_{\varepsilon}(z)}(x) d{\rm vol}_\rho+ O\left( \sqrt{\frac{\log(n)}{n\varepsilon^{d-4}}}\right)\\
=&\ (1+ O(\varepsilon))-\hspace{-3.8mm}\int_{B_{\varepsilon}(z)}  Q^{(s,l)}d{\rm vol}_\rho+ O\left( \sqrt{\frac{\log(n)}{n\varepsilon^{d-4}}}\right)
\end{align*}
and hence a control of ${\tt LD}_{AD}=O\left(\sqrt{\frac{\log(n)}{n\varepsilon^{d-4}}}\right)$ with probability higher than $1-Cn^{-3}$ for some $C>0$ when $n$ is sufficiently large. 

Besides $Y_m^{(s,l)}=({\bf q}_m)_s({\bf q}_m)_l$, we have to control terms like $Y_m^{(s)}:=({\bf q}_m)_s$, $Y_m^{(ss,l)}:=({\bf q}_m)_s^2({\bf q}_m)_l$, $Y_m^{(s,l,t)}:=({\bf q}_m)_s({\bf q}_m)_l({\bf q}_m)_t$, etc. The computation is similar and can be summarized by the following lemma. Since the proof of this lemma is the same as the above, we omit details.

\begin{Lma} \label{CorBern}
Suppose $X_1,\ldots,X_n$ are i.i.d. sampled from a random variable $X$.
Let $\varepsilon>0$ and $\lambda\in \mathbb{N}$. Suppose $E(X)=0$, $|X|\leq a=\Theta(\varepsilon^{-d+\lambda})$ and ${{\tt Var}}(X)=\Theta(\varepsilon^{-d+2\lambda})$. Then we have
$$
\mathbb{P}\left\{\left|\frac{1}{n}\sum_{i=1}^n X_i\right|> \sqrt{\frac{\log(n)}{n\varepsilon^{d-2\lambda}}}\right\}\leq Cn^{-3}
$$
for some $C>0$ when $n$ is sufficiently larger if $\varepsilon\to 0$  and $\frac{n\varepsilon^d}{\log(n)}\to \infty$ hold when $n\to \infty$.
\end{Lma}

Note that the control of ${\tt LD}_{AD}$ is a special case of Lemma \ref{CorBern}. Based on this Lemma, the other terms of ${\tt LD}$ can be computed similarly and we have the following lemma.

\begin{Lma} \label{LD control bound}
Denote $\omega:=\sqrt{\frac{\log(n)}{n\varepsilon^d}}$ to simplify the notation. Then
$${\tt LD}=\left[\begin{array}{cccc}
{\tt LD}_{AA}&{\tt LD}_{AB}& {\tt LD}_{AC}& {\tt LD}_{AD}\\
- &{\tt LD}_{BB}&{\tt LD}_{BC}& {\tt LD}_{BD}\\ 
- & - & {\tt LD}_{CC}&{\tt LD}_{CD}\\
- & - &- &{\tt LD}_{DD}
\end{array}\right]
=\left[\begin{array}{cccc}
0&O(\varepsilon\omega)& O(\varepsilon^{2}\omega)& O(\varepsilon^{2}\omega)\\
- &O(\varepsilon^{2}\omega)&O(\varepsilon^{3}\omega)& O(\varepsilon^{3}\omega)\\
- & - & O(\varepsilon^{4}\omega)&O(\varepsilon^{4}\omega)\\
- & - &- &O(\varepsilon^{4}\omega)
\end{array}\right]
$$
with probability higher than $1-Cn^{-3}$ for some $C>0$ when $n$ is sufficiently large if $\varepsilon\to 0$  and $\frac{n\varepsilon^d}{\log(n)}\to \infty$ hold when $n\to \infty$.
\end{Lma}

To simplify the heavy notation, from now on we follow the same notation relationship used in \eqref{Definition Ymsl}, \eqref{Definition Qmsl} and \eqref{Definition Xmsl} and consider notations $Y_m^{(s\ldots s)}:=({\bf q}_m)^\alpha_s$ if there are $\alpha$ copies of $s$ in the superscript of $Y_m$, $Y_m^{(s\ldots s,l\ldots l)}:=({\bf q}_m)^\alpha_s({\bf q}_m)_l^\beta$ if there are $\alpha$ copies of $s$ and $\beta$  copies of $l$ in the superscript of $Y_m$,  $Y_m^{(s\ldots s,l\ldots l,t\ldots t)}:=({\bf q}_m)^\alpha_s({\bf q}_m)^\beta_l({\bf q}_m)^\gamma_t$ if there are $\alpha$ copies of  $s$, $\beta$ copies of $l$ and $\gamma$ copies of $t$ in the superscript of $Y_m$, etc., and the associated $Q$ and $X$ notations with the proper superscripts.

\begin{proof}
It is easy to see that ${\tt LD}_{AA}=0$. For each of the other terms, we have to compute the crucial term ${\tt Var}(X)$ as follows. 
\begin{itemize}
\item $L_{AB}$: 
Let $Y_m^{(l)}:=({\bf q}_m)_l$ and $X_i^{(l)}:=Q_i^{(l)}\varepsilon^{-d}\chi_{ B_\varepsilon(z)}(X_i)-E(Q_i^{(l)}\varepsilon^{-d}\chi_{ B_\varepsilon(z)}(X_i))$. It follows from Lemma \ref{lemma of local poly integration}
that 
\begin{align*}
{\tt Var}(X^{(l)})&\,=\varepsilon^{-2d}\int_{\mathcal{B}}{\bf q}_l^2  d{\rm vol}_\rho -\varepsilon^{-2d}\left(\int_{\mathcal{B}} {\bf q}_l  d{\rm vol}_\rho\right)^2\\
&\,=\varepsilon^{-d}\rho(z)C_{0,0}\alpha_i+O(\varepsilon^{-d+3}),\ i=1\mbox{ or }2\\
&\,=\Theta(\varepsilon^{-d+2})
\end{align*}
and, by Lemma \ref{CorBern}, ${\tt LD}_{AB}=O(\varepsilon\omega)$.

\item $L_{AC}$: By a direct computation as that in Lemma \ref{lemma of local poly integration}, we consider $Y_m^{(ll)}:=({\bf q}_m)_l^2$ and derive  
\begin{align*}
{\tt Var}(X^{(ll)})=&\ \varepsilon^{-2d}\int_{\mathcal{B}}{\bf q}_l^4  d{\rm vol}_\rho-\varepsilon^{-2d}\left(\int_{\mathcal{B}} {\bf q}_l^2  d{\rm vol}_\rho\right)^2\\
=&\ \varepsilon^{-d}\rho(z)C_{0,0}\beta_i
+O(\varepsilon^{-d+5}),\ i=1\mbox{ or }3\\ 
=&\ \Theta(\varepsilon^{-d+4})
\end{align*}
and thus ${\tt LD}_{AC}=O(\varepsilon^2\omega)$.
\item $L_{AD}$: We had shown that ${\tt Var}(X^{(s,l)})=\Theta(\varepsilon^{-d+4})$ in the previous paragraph. 

\item $L_{BB}$: The computation of $L_{BB}$ involves exactly the same terms as $L_{AC}$ and $L_{AD}$. So the variance is also in $\Theta(\varepsilon^{-d+4})$.

\end{itemize}
Similarly, it is easy to check that ${\tt Var}(X^{(s,ll)})$ and ${\tt Var}(X^{(s,l,t)})$ are of $\Theta(\varepsilon^{-d+6})$ while ${\tt Var}(X^{(ss,ll)})$, ${\tt Var}(X^{(ss,l,t)})$, and ${\tt Var}(X^{(s,l,t,r)})$ are all of $\Theta(\varepsilon^{-d+8})$.
Therefore, using Lemma \ref{CorBern}, we obtain the desired large deviations. 
\end{proof}
With the relationship assumption of $n$ and $\varepsilon$, from Lemmas \ref{C.W bound of L and L0} and \ref{LD control bound}, we conclude that $\frac{1}{k_z} Z^T Z$ satisfies 
{\small
\begin{align}
&\ \frac{1}{k_z} Z^T Z=L^0+W+{\tt LD}\nonumber\\
=&\nonumber\ \left[\begin{array}{cccc}
L^0_{AA}&L^0_{AB}& L^0_{AC}& L^0_{AD}\\
- &L^0_{BB}&L^0_{BC}&L^0_{BD}\\ 
- & - & L^0_{CC}&L^0_{CD}\\
- & - &- &L^0_{DD}
\end{array}\right]
+
\left[\begin{array}{cccc}
W_{AA}&W_{AB}& W_{AC}& W_{AD}\\
- &W_{BB}&W_{BC}&W_{BD}\\ 
- & - & W_{CC}&W_{CD}\\
- & - &- &W_{DD}
\end{array}\right]
+
\left[\begin{array}{cccc}
{\tt LD}_{AA}&{\tt LD}_{AB}& {\tt LD}_{AC}& {\tt LD}_{AD}\\
- &{\tt LD}_{BB}&{\tt LD}_{BC}&{\tt LD}_{BD}\\ 
- & - & {\tt LD}_{CC}&{\tt LD}_{CD}\\
- & - &- &{\tt LD}_{DD}
\end{array}\right]\\
=&\ \left[\begin{array}{cccc}
1&L^0_{AB}& L^0_{AC}& {\bf 0}\\
- &L^0_{BB}&L^0_{BC}&L^0_{BD}\\ 
- & - & L^0_{CC}&{\bf 0}\\
- & - &- &L^0_{DD}
\end{array}\right]
+
\left[\begin{array}{cccc}
0 & O(\varepsilon^2) & O(\varepsilon^3) &O(\varepsilon^3) \\
- & O(\varepsilon^3) & O(\varepsilon^4) &O(\varepsilon^4) \\
    - & - & O(\varepsilon^5) &O(\varepsilon^5) \\
    - & - & - &O(\varepsilon^5)
\end{array}\right]
+
\left[\begin{array}{cccc}
0&O(\varepsilon\omega)& O(\varepsilon^{2}\omega)& O(\varepsilon^{2}\omega)\\
- &O(\varepsilon^{2}\omega)&O(\varepsilon^{3}\omega)& O(\varepsilon^{3}\omega)\\
- & - & O(\varepsilon^{4}\omega)&O(\varepsilon^{4}\omega)\\
- & - &- &O(\varepsilon^{4}\omega)
\end{array}\right]\nonumber
\end{align}
}
with probability higher than $1-Cn^{-3}$ for some $C>0$ when $n$ is sufficiently large. 
Hence the large
deviation of $\frac{1}{k_z} Z^T Z{\bf f}_{Taylor}$ is 
$\left[\begin{array}{c}
O(\varepsilon\omega)\\
O(\varepsilon\omega)\\
O(\varepsilon^2\omega)\\
O(\varepsilon^2\omega)
\end{array}\right]$.
On the other hand, 
we are going to show that the large deviation of
$\frac{1}{k_z}Z^T{\bf f}$ is
$\left[\begin{array}{c}
O(\omega)\\
O(\varepsilon\omega)\\
O(\varepsilon^2\omega)\\
O(\varepsilon^2\omega)
\end{array}\right]$
and thus 
$$V
=\left[
\begin{array}{c}
V_A  \\
V_B \\
V_C \\
V_D
\end{array}\right]
=\left[
\begin{array}{c}
O(\varepsilon^{2+\min(\tau,\kappa)})  \\
O(\varepsilon^{3+\min(\tau,\kappa)}) \\
O(\varepsilon^{4+\min(\tau,\kappa)}) \\
O(\varepsilon^{4+\min(\tau,\kappa)}) 
\end{array}\right]
+
\left[\begin{array}{c}
O(\omega)\\
O(\varepsilon\omega)\\
O(\varepsilon^2\omega)\\
O(\varepsilon^2\omega)
\end{array}\right].$$

\begin{Lma} \label{LD ZTf}
The large deviation of
$\frac{1}{k_z}Z^T{\bf f}=\overline{L}+\overline{\tt LD}$ is
$\overline{\tt LD}=\left[\begin{array}{c}
O(\omega)\\
O(\varepsilon\omega)\\
O(\varepsilon^2\omega)\\
O(\varepsilon^2\omega)
\end{array}\right]$.
\end{Lma}

\begin{proof}
The computation is very similar to what we have done in Lemma \ref{LD control bound}.
For the first block $\overline{\tt LD}_A$, let $\widehat X_i:=f(X_i)\varepsilon^{-d}\chi_{ B_\varepsilon(z)}(X_i)-E\left(f(X_i)\varepsilon^{-d}\chi_{ B_\varepsilon(z)}(X_i)\right)$. 
Then $|\widehat X_i|\leq a= \Theta(\varepsilon^{-d})$ and 
\begin{align*}
{\tt Var}(\widehat X_i)
=&\ \varepsilon^{-2d}\int_{\mathcal{B}} f^2 d{\rm vol}_{\rho} -\varepsilon^{-2d}\left(\int_{\mathcal{B}} f d{\rm vol}_{\rho} \right)^2\\
=&\ \rho(z)C_{0,0}f^2(z) \varepsilon^{-d}+O(\varepsilon^{-d+1})\\
=&\ \Theta(\varepsilon^{-d}).
\end{align*}
Using the technical Lemma \ref{CorBern} again, this implies that $\overline{\tt LD}_A=O(\omega)$. For the second block
$\overline{\tt LD}_B$,
let $\widehat X_i^{(l)}:=Q_i^{(l)}f(X_i)\varepsilon^{-d}\chi_{ B_\varepsilon(z)}(X_i)-E\left(Q_i^{(l)}f(X_i)\varepsilon^{-d}\chi_{ B_\varepsilon(z)}(X_i)\right)$.
Then $|\widehat X^{(l)}_i|\leq a= \Theta(\varepsilon^{-d+1})$ and by the proof of Lemma \ref{lemma of local poly integration}(ii) we can show that 
\begin{align*}
{\tt Var}(\widehat X_i^{(l)})
=&\ \varepsilon^{-2d}\int_{\mathcal{B}}\langle  x-z,e_l\rangle^2 f^2 d{\rm vol}_{\rho} -\varepsilon^{-2d}\left(\int_{\mathcal{B}} \langle  x-z,e_l\rangle f d{\rm vol}_{\rho} \right)^2\\
=&\ \varepsilon^{-d}\rho(z)C_{0,0}\alpha_j f^2(z)+O(\varepsilon^{{-d+3}}),\ \ j= 1 \mbox{ or } 2\\
=&\ \Theta(\varepsilon^{-d+2})\,.
\end{align*}
Thus, $\overline{\tt LD}_B= O(\varepsilon\omega)$.

For the third block $\overline{\tt LD}_C$,
let $\widehat X_i^{(ll)}:=Q_i^{(ll)}f(X_i)\varepsilon^{-d}\chi_{ B_\varepsilon(z)}(X_i)-E\left(Q_i^{(ll)}f(X_i)\varepsilon^{-d}\chi_{ B_\varepsilon(z)}(X_i)\right)$.
Then by similar calculation we have $|\widehat X^{(ll)}_i|\leq \Theta(\varepsilon^{-d+2})$ and
${\tt Var}(\widehat X_i^{(ll)})= \Theta(\varepsilon^{-d+4})$, which implies $\overline{\tt LD}_C=O\left(\varepsilon^2\omega\right)$.

For the fourth block, let $\widehat X_i^{(s,t)}:=Q_i^{(s,t)}f(X_i)\varepsilon^{-d}\chi_{ B_\varepsilon(z)}(X_i)-E\left(Q_i^{(s,t)}f(X_i)\varepsilon^{-d}\chi_{ B_\varepsilon(z)}(X_i)\right)$.
Then we have $|\widehat X^{(s,t)}_i|\leq \Theta(\varepsilon^{-d+2})$ and
${\tt Var}(\widehat X_i^{(s,t)})= \Theta(\varepsilon^{-d+4})$ and 
thus $\overline{\tt LD}_D=O\left(\varepsilon^2\omega\right)$.
\end{proof}

In the next section, we will combine all the estimates together to compute $\left(\frac{1}{k_z}Z^TZ\right)^{-1}V$ and finish the proof of our main theorem.

\subsection{Finish the proof of Theorem \ref{main}}

\begin{proof}
To compute the inverse of $\frac{1}{k_z}Z^TZ=L^0+W+{\tt LD}=L^0+E$, one can use the formula 
\begin{equation}\label{formula: inversion with approximation}
(L^0+E)^{-1}=(L^0)^{-1}-(L^0)^{-1}E(L^0)^{-1}+((L^0)^{-1}E)^2(L^0)^{-1}-((L^0)^{-1}E)^3(L^0)^{-1}\pm\dots.
\end{equation}
We will see the inverse exists almost surely when $\varepsilon\to 0$ and $\frac{n\varepsilon^d}{\log(n)}\to\infty$ as $n\to\infty$. 
Recall that $L^0$ is given in (\ref{L0}). One can compute its inverse directly and derive
\begin{align}\label{inverseL0}
(L^0)^{-1}
=\left[\begin{array}{c|cccc|ccccc|c}
l_1 & 0& \cdots& 0 &l_3  &l_5 &\cdots      &  &l_5   & l_{9}       &0\\ \hline
- &\beta_2l_2&0&\cdots&0 &       0  &\cdots&\cdots&\cdots&  0       &\\ 
- & - &\ddots &\ddots&\vdots &\vdots&\ddots&\ddots&\ddots&\vdots    &((L^0)^{-1})_{BD}\\ 
- & - &- &\beta_2l_2& 0&         0  &\cdots&\cdots&\cdots&  0       &\\ 
- & - &- &- &l_4&               l_6 &\cdots&\cdots&l_{6}& l_{10}    &\\ \hline 
- & - &- &- &-  &               l_7 & l_8 &\cdots   & l_8  & l_{11} &0\\
- & - &- &- &-  &                -  & l_7 &\ddots   &\vdots& l_{11} &0\\
- & - &- &- &-  &                -  &  -  &\ddots   & l_8  & \vdots &0\\
- & - &- &- &-  &                -  &  -  &  -      & l_7  & l_{11} &0\\
- & - &- &- &-  &                -  &  -  &  -      &  -   & l_{12} &0\\ \hline 
- & - &- &- &-  &                -  &  -  &  -      &  -   &        &((L^0)^{-1})_{DD}  
\end{array}\right],
\end{align}
with 
$$
((L^0)^{-1})_{BD}=
\left[
\begin{array}{cccc|cccc|c|ccc|cc|c}
0&\cdots&0&-\mu_2l_2 &   0&\cdots& 0& 0    &  \cdots  &  0&0&0&      0&0      &   0 \\   
0&\cdots&0& 0    &   0&\cdots& 0&-\mu_2l_2 &  \cdots  &  0&0&0&      0&0      &   0 \\   
 &      & &      &    &      &  &      &  \ddots  &   & & &       &       &    \\   
0&\cdots&0& 0    &   0&\cdots& 0& 0    &  \cdots  &  0&0&-\mu_2l_2&  0&0      &   0 \\
0&\cdots&0& 0    &   0&\cdots& 0& 0    &  \cdots  &  0&0&0&      0&-\mu_2l_2  &   0 \\
0&\cdots&0& 0    &   0&\cdots& 0& 0    &  \cdots  &  0&0&0&      0&0      &-\mu_2l_2\\
0&\cdots&0& 0    &   0&\cdots& 0& 0    &  \cdots  &  0&0&0&      0&0      &   0 \\   
\end{array}
\right]
$$ 
and 
$$
((L^0)^{-1})_{DD}
=\left[
\begin{array}{cccc|c|ccc|cc|c}
\beta_4^{-1} & &  &  &&    & & &     &    & \\
  & \ddots&  &  &&    & & &     &    & \\
  & & \beta_4^{-1}&  &&    & & &     &    & \\
  & &     &\alpha_2l_2 &&    & & &     &            & \\ \hline
  & &     &&\ddots&  & & &     &            & \\ \hline
  & &     &     &&  \beta_4^{-1}  & & &     &        & \\
  & &     &     &&    &\beta_4^{-1}& &     &          & \\
  & &     &     &&    &      &\alpha_2l_2&     &    & \\ \hline
  & &     &     &&    &      &       &\beta_4^{-1}&         &\\
  & &     &     &&    &      &       &       & \alpha_2l_2& \\ \hline
  & &     &     &&    &      &       &   &        & \alpha_2l_2
\end{array}
\right],
$$ 
where all $l_j$'s are constants depending on $\alpha,\beta,\gamma,\mu$ and the dimension $d$. The exact values of these $l_j$'s are complicated and not important. What matters is the order of them. Indeed, by examining these constants, one has $(L^0)^{-1}
=\left[\begin{array}{cccc}
O(1)&O(\varepsilon^{-1})&O(\varepsilon^{-2})&       0           \\
 -  &O(\varepsilon^{-2})&O(\varepsilon^{-3})&O(\varepsilon^{-3})\\
 -  &         -         &O(\varepsilon^{-4})&       0           \\
 -  &         -         &          -        &O(\varepsilon^{-4})
\end{array}\right]$. For curious readers, the most frequently occurring number, $l_2$, which is also the simplest one, is actually $(\alpha_2\beta_2-\mu_2^2)^{-1}$. So $l_2$ is of $O(\varepsilon^{-6})$ and thus the orders of $((L^0)^{-1})_{BD}$ and $((L^0)^{-1})_{DD}$ can be observed.

Hence, one can further derive 
\begin{align}
\begin{split}
(L^0)^{-1}W
=&\
\left[\begin{array}{cccc}
O(\varepsilon^1)    &O(\varepsilon^2) &O(\varepsilon^3) &O(\varepsilon^3)\\
O(1)                &O(\varepsilon^1) &O(\varepsilon^2) &O(\varepsilon^2)\\ 
O(\varepsilon^{-1}) &O(1)             &O(\varepsilon^1) &O(\varepsilon^1)\\
O(\varepsilon^{-1}) &O(1)             &O(\varepsilon^1) &O(\varepsilon^1)
\end{array}\right],\\
((L^0)^{-1}W)^2
=&\ 
\left[\begin{array}{cccc}
O(\varepsilon^2) &O(\varepsilon^3) &O(\varepsilon^4) &O(\varepsilon^4)\\
O(\varepsilon^1) &O(\varepsilon^2) &O(\varepsilon^3) &O(\varepsilon^3)\\ 
O(1)             &O(\varepsilon^1) &O(\varepsilon^2) &O(\varepsilon^2)\\
O(1)             &O(\varepsilon^1) &O(\varepsilon^2) &O(\varepsilon^2)
\end{array}\right],\mbox{ etc}.
\end{split}
\end{align}

Moreover, $$(L^0)^{-1}{\tt LD}
=\left[\begin{array}{cccc}
O(\omega)&O(\varepsilon\omega)& O(\varepsilon^{2}\omega)& O(\varepsilon^{2}\omega)\\
O(\varepsilon^{-1}\omega)&O(\omega)&O(\varepsilon\omega)& O(\varepsilon\omega)\\
O(\varepsilon^{-2}\omega) & O(\varepsilon^{-1}\omega) & O(\omega)&O(\omega)\\
O(\varepsilon^{-2}\omega) & O(\varepsilon^{-1}\omega) & O(\omega)&O(\omega)
\end{array}\right],
$$
$$(L^0)^{-1}V=\left[
\begin{array}{c}
O(\varepsilon^{2+\min(\tau,\kappa)})  \\
O(\varepsilon^{1+\min(\tau,\kappa)}) \\
O(\varepsilon^{\min(\tau,\kappa)}) \\
O(\varepsilon^{\min(\tau,\kappa)}) 
\end{array}\right]
+
\left[\begin{array}{c}
O(\omega)\\
O(\varepsilon^{-1}\omega)\\
O(\varepsilon^{-2}\omega)\\
O(\varepsilon^{-2}\omega)
\end{array}\right],$$ etc.
It is easy to see that the order of $(L^0)^{-1}E=(L^0)^{-1}W+(L^0)^{-1}{\tt LD}$ becomes higher and higher after self-multiplication. 
Combining all of these, we have 
$$\left(Z^TZ\right)^{-1}Z^T{\bf f}
= {\bf f}_{Taylor}-\left(\frac{1}{k_z}Z^TZ\right)^{-1}V
= {\bf f}_{Taylor}
+
\left[
\begin{array}{c}
O(\varepsilon^{2+\min(\tau,\kappa)})  \\
O(\varepsilon^{1+\min(\tau,\kappa)}) \\
O(\varepsilon^{\min(\tau,\kappa)}) \\
O(\varepsilon^{\min(\tau,\kappa)}) 
\end{array}\right]
+
\left[\begin{array}{c}
O(\omega)\\
O(\varepsilon^{-1}\omega)\\
O(\varepsilon^{-2}\omega)\\
O(\varepsilon^{-2}\omega)
\end{array}\right].$$
\end{proof}

\section{An alternative proof of the main theorem -- a direct expansion approach}\label{Appendix proof direct approach}

In this section, we showcase an attempt to prove the main theorem by a direct expansion; that is, evaluate the bias and variance terms between $\frac{1}{k_z} Z^T Z$ and its continuous counterpart, as well as those between $\frac{1}{k_z} Z^T{\bf f}$ and its continuous counterpart (see Lemma \ref{ZTf} for an example), and then multiply $\left(Z^T Z\right)^{-1}$ and $Z^T{\bf f}$. 
As we will demonstrate, this approach is possible, but when the point is near the boundary, this approach is limited.

To save readers' time, here is an outline of encountered challenge. We will show that the asymptotic orders of entries of $Z^T Z$ and $Z^T{\bf f}$ are heterogeneous. Upon inversion, heterogeneous asymptotic orders in $Z^T Z$ give rise to several new terms of heterogeneous asymptotic orders that complicate the multiplication of $\left(Z^T Z\right)^{-1}$ and $Z^T{\bf f}$. The critical challenge comes from the inversion, which reverses the order of the dominant entries. These entries with reversed orders (see \eqref{G} and \eqref{hot} for an example), when multiplied with entries of $Z^T{\bf f}$, lead to potentially nontrivial low order terms that could contaminate the main terms we have interest. To address all terms with heterogeneous asymptotic orders, we require higher-order local approximations of the function, density function, and curvature when evaluating entries of $Z^T{\bf f}$ compared with the proof shown in Section \ref{section: proof of main theorem heuristic approach}, further complicating the analysis.

To appreciate this challenge, see $\partial
_j\rho$, $\mathcal{U}_{jj}$ and $\mathcal{V}_{s,t}$ that appear in Corollary \ref{intfu} as an example. These terms appear in $W_{AC},W_{AD}$ and $W_{BD}$ and eventually make the entries of $\left(Z^T Z\right)^{-1}$ of order $O(\varepsilon^{-4}
)$. While we expect to balance these asymptotically blowing up terms by matching them with the higher order terms in $Z^T{\bf f}$ through the multiplying $\left(Z^T Z\right)^{-1}$ and $Z^T{\bf f}$, the terms $\partial_j\rho$, $\mathcal{U}_{jj}$ and $\mathcal{V}_{s,t}$ remain. Unfortunately, some of these terms persist in the leading order term, potentially contaminating the Hessian term of interest (cf. Lemma \ref{ZTf}). Thus, we must carefully track these terms and look for cancellation.
By computing the precise leading quantities, we will see how there terms cancel out, allowing us to obtain the desired result: 
$$
\left(\frac{1}{k_z}Z^TZ\right)^{-1}\frac{1}{k_z}Z^T{\bf f}={\bf f}_{Taylor}
+\left[
\begin{array}{c}
O(\varepsilon^2)  \\
O(\varepsilon^2) \\
O(\varepsilon) \\
O(\varepsilon) 
\end{array}\right]
+
\left[\begin{array}{c}
O(\omega)\\
O(\varepsilon^{-1}\omega)\\
O(\varepsilon^{-2}\omega)\\
O(\varepsilon^{-2}\omega)
\end{array}\right].
$$
This gives an independent proof to part of our Main Theorem when $z\in M\setminus M_{\sigma}$. 
In this Section, we focus on a fixed $z\in M\setminus M_\sigma$ that is away from the boundary. The proof is divided into several steps like that in Section \ref{section: proof of main theorem heuristic approach}.

\subsection{Step 1: The bias term of $\frac{1}{k_z}Z^TZ$}

In the following, we again use dashes to indicate the lower triangular portion of a symmetric matrix.
As in Section \ref{section: proof of main theorem heuristic approach}, we consider 
\[
\frac{1}{k_z}Z^TZ=L+{\tt LD}, 
\]
where
$$
L=\left[\begin{array}{cccc}
1&L_{AB}& L_{AC}& L_{AD}\\
- &L_{BB}&L_{BC}&L_{BD}\\ 
- & - & L_{CC}&L_{CD}\\
- & - &- &L_{DD}
\end{array}\right]\in \mathbb{R}^{(1+d+d+\frac{d(d-1)}{2})\times (1+d+d+\frac{d(d-1)}{2})}.
$$
In the following lemma, we show that it can be well approximated by the following simplified matrix
{\small \begin{align}\label{definition L0}
L^0
=\left[\begin{array}{c|cccc|cccc|c}
1 & \alpha\rho^{-1}\partial_1\rho & \cdots& \cdots &\alpha\rho^{-1}\partial_d\rho  &\alpha & \cdots&\cdots& \alpha & {\bf 0}\\ \hline
- &\alpha &0& \cdots & 0&3\beta\rho^{-1}\partial_1\rho&\beta\rho^{-1}\partial_1\rho&\cdots&\beta\rho^{-1}\partial_1\rho&{\bf 0}\\ 
- & - &\ddots &\ddots&\vdots &\beta\rho^{-1}\partial_2\rho&\ddots& \ddots&\vdots&{\bf 0}\\ 
- & - &- &\ddots& 0&\vdots& \ddots&\ddots& \beta\rho^{-1}\partial_{d-1}\rho&{\bf 0}\\ 
- & - &- &- &\alpha  &\beta\rho^{-1}\partial_d\rho&\cdots&\beta\rho^{-1}\partial_d\rho&3\beta\rho^{-1}\partial_d\rho&{\bf 0}\\ \hline 
- & - &-&-&-& 3\beta&\beta&\cdots&\beta&{\bf 0}\\
- & - &-&-&-& -    &\ddots&\ddots&\vdots&{\bf 0}\\
- & - &-&-&-& -    &-     &\ddots&\beta&{\bf 0}\\
- & - &-&-&-& -    &-     &-     &3\beta &{\bf 0}\\ \hline
- & - &-&-&-& -    &-     &-     &-     &\beta I_{\frac{d(d-1)}{2}}
\end{array}\right],
\end{align}}
with $\alpha=\frac{\varepsilon^2}{d+2}$, $\beta=\frac{\varepsilon^{4}}{(d+2)(d+4)}$, $\rho:=\rho(z)$ and
$\partial_j\rho$ denoting partial derivative $\nabla_j\rho|_{z}$, $j=1,\dots, d$. In this manner we say that $L$ approximates $L^0$, and the approximation bias will be quantified precisely. 

\begin{Lma} \label{D.W bound of L and L0}
Let $\mathcal{U}_{jj}$ and $\mathcal{V}_{s,t}$ be defined as in Corollary \ref{intfu}. Define $L^0$ as in above. Then
$$W:=L-L^0
=\left[\begin{array}{cccc}
W_{AA}&W_{AB}&W_{AC}&W_{AD}\\
-     &W_{BB}&W_{BC}&W_{BD}\\ 
-     & -    &W_{CC}&W_{CD}\\
-     & -    &-     &W_{DD}
\end{array}\right]
=\left[\begin{array}{cccc}
0 &O(\varepsilon^3) & W^0_{AC}+O(\varepsilon^5) & W^0_{AD}+O(\varepsilon^5)\\
- &O(\varepsilon^4) &O(\varepsilon^5)&W_{BD}^0+O(\varepsilon^5)\\ 
- & - & O(\varepsilon^6)&O(\varepsilon^6)\\
- & - & -    &O(\varepsilon^6) 
\end{array}\right],$$
where $[W^0_{AC}]_j=\beta\left(\rho^{-1}\nabla_j^2\rho-\frac{1}{d+2}\rho^{-1}\Delta\rho+ \mathcal{U}_{jj}\right)$,
$W^0_{AD}:=\left[[W^0_{AD}]_{11}\ [W^0_{AD}]_{12}\ \cdots \ [W^0_{AD}]_{(d-1)d}\right]$ 
is a $1\times\frac{d(d-1)}{2}$ matrix with $[W^0_{AD}]_{st}=\beta(\rho^{-1}\nabla_s\nabla_t\rho+\mathcal{V}_{s,t})$,
and 
$$\footnotesize  W_{BD}^0=
\beta\rho^{-1}\left[ 
\begin{array}{ccccc|cccc|ccc|c|cc|c}
\partial_2\rho&\partial_3\rho&\partial_4\rho&\cdots&\partial_d\rho
      &0 &  0&\cdots&   0
         &  0&\cdots&   0
             &\cdots&  0&0& 0 \\
\partial_1\rho&     0&\cdots&\cdots&0
      &\partial_3\rho&\partial_4\rho&\cdots&\partial_d \rho
             &   0  &\cdots&   0
                    &\cdots&  0&0& 0  \\
0&\ddots& & & \vdots
 &\partial_2\rho&   0  &\cdots&  0 
        &\partial_4\rho&\cdots&\partial_d\rho
               &\cdots&  0&0&0 \\
\vdots&  & \ddots &  &\vdots&
0&\ddots& &  \vdots 
 &\partial_3\rho&\cdots&  0
        &\cdots&\partial_{d-1}\rho&\partial_d\rho&0 \\
\vdots&    &     &\ddots&  0  
 & \vdots  &    &\ddots&  0  
      &  \vdots &\ddots&\vdots  
            &\cdots&\partial_{d-2}\rho&0&\partial_d\rho \\
0& \cdots & \cdots& 0 &\partial_1\rho
 & 0 & \cdots &    0  &\partial_2 \rho
     & 0 &\cdots&\partial_3\rho
         &\cdots&0&\partial_{d-2}\rho&\partial_{d-1} \rho
\end{array}\right].
$$
\end{Lma}

\begin{proof} 
For the simplicity of notations, we omit the symbol $\iota$ in the following proof. First, it is easy to see that
$\displaystyle W_{AA}=-\hspace{-3.8mm}\int_{\widetilde B} 1 d{\rm vol}_{\rho} - 1=0$. Now we compute $W_{AB}$. By Corollary \ref{intfu}(i) with $f=1$, we have
$\displaystyle -\hspace{-3.8mm}\int_{\widetilde B}\langle x-z,u_j \rangle d{\rm vol}_{\rho}=\alpha \rho^{-1}\partial_j\rho+ O(\varepsilon^3)$.Hence, $W_{AB}=L_{AB}-L_{AB}^0$ is a vector of $O(\varepsilon^3)$. 

As for $W_{AC}$, by using Corollary \ref{intfu}(ii) with $f=1$, we have 
$$-\hspace{-3.8mm}\int_{\widetilde B}\langle x-z,e_j \rangle^2 d{\rm vol}_{\rho}
=\alpha +\beta\left(\rho^{-1}\nabla_j^2\rho-\frac{1}{d+2}\rho^{-1}\Delta\rho+ \mathcal{U}_{jj}\right)+O(\varepsilon^5).$$
Similarly, the entries in the vector $W_{AC}=L_{AC}-L^0_{AC}$ are $$\displaystyle-\hspace{-3.8mm}\int_{\widetilde B}\langle x-z,u_j \rangle^2 d{\rm vol}_\rho-\alpha=\beta\left(\rho^{-1}\nabla_j^2\rho-\frac{1}{d+2}\rho^{-1}\Delta\rho+ \mathcal{U}_{jj}\right)+O(\varepsilon^5)=:[W^0_{AC}]_j+O(\varepsilon^5)$$
and, using Corollary \ref{intf}(iii), one can derive
\begin{align}\label{error epsilon5 for Hessian}
[W_{AD}]_{st}=[W^0_{AD}]_{st}+O(\varepsilon^5)
= \beta(\rho^{-1}\nabla_s\nabla_t\rho+\mathcal{V}_{s,t})+O(\varepsilon^5).
\end{align}

To find $W_{BB}, W_{BC}$, and $W_{BD}$, one has to compute
$\displaystyle-\hspace{-3.6mm}\int_{\widetilde B}\langle x-z,u_j \rangle\langle x-z,u_k \rangle d{\rm vol}_\rho$,
$\displaystyle-\hspace{-3.6mm}\int_{\widetilde B}\langle x-z,u_j \rangle\langle x-z,u_k \rangle^2 d{\rm vol}_\rho$,
and 
$\displaystyle-\hspace{-3.6mm}\int_{\widetilde B}\langle x-z,u_j \rangle\langle x-z,u_s \rangle\langle x-z,u_t \rangle d{\rm vol}_\rho$,
respectively. Indeed,
$W_{BB}=L_{BB}-L^0_{BB}
=\displaystyle-\hspace{-3.6mm}\int_{\widetilde B}\langle x-z,u_j \rangle\langle x-z,u_k \rangle d{\rm vol}_\rho-\alpha I=O(\varepsilon^4)$ can be easily derived from Corollary \ref{intfu}(ii)(iii). Similarly, $W_{BC}=O(\varepsilon^6)$
and the matrix $W_{BD}$ can be derived from Corollary \ref{intfu}(iv). Other terms such as $L_{CC}$ and $L_{DD}$ are similarly obtained by using Corollary \ref{intfu}(v).
\end{proof}

\subsection{Step 2: The large deviation estimate for $\frac{1}{k_z}Z^TZ$} 

The large deviation terms of $\frac{1}{k_z} Z^T Z$ can also be derived from Lemma \ref{LD control bound} because the lemma is valid no matter the point $z$ is near or away from the boundary.
With the relationship assumption of $n$ and $\varepsilon$, from Lemmas \ref{D.W bound of L and L0} and \ref{LD control bound}, we conclude that $\frac{1}{k_z} Z^T Z$ satisfies 
\begin{align}
&\ \frac{1}{k_z} Z^T Z=L^0+W+{\tt LD}\label{Expansion of Z^TZ for the direct approach}\\
=&\nonumber\ \left[\begin{array}{cccc}
L^0_{AA}&L^0_{AB}& L^0_{AC}& L^0_{AD}\\
- &L^0_{BB}&L^0_{BC}&L^0_{BD}\\ 
- & - & L^0_{CC}&L^0_{CD}\\
- & - &- &L^0_{DD}
\end{array}\right]
+
\left[\begin{array}{cccc}
W_{AA}&W_{AB}& W_{AC}& W_{AD}\\
- &W_{BB}&W_{BC}&W_{BD}\\ 
- & - & W_{CC}&W_{CD}\\
- & - &- &W_{DD}
\end{array}\right]
+
\left[\begin{array}{cccc}
{\tt LD}_{AA}&{\tt LD}_{AB}& {\tt LD}_{AC}& {\tt LD}_{AD}\\
- &{\tt LD}_{BB}&{\tt LD}_{BC}&{\tt LD}_{BD}\\ 
- & - & {\tt LD}_{CC}&{\tt LD}_{CD}\\
- & - &- &{\tt LD}_{DD}
\end{array}\right]\\
=&\ \left[\begin{array}{cccc}
1&L^0_{AB}& L^0_{AC}& {\bf 0}\\
- &L^0_{BB}&L^0_{BC}&{\bf 0}\\ 
- & - & L^0_{CC}&{\bf 0}\\
- & - &- &L^0_{DD}
\end{array}\right]
+
\left[\begin{array}{cccc}
0 & O(\varepsilon^3) & O(\varepsilon^4) & O(\varepsilon^4)\\
- & O(\varepsilon^4) & O(\varepsilon^5) & O(\varepsilon^4)\\ 
- & - & O(\varepsilon^6)&O(\varepsilon^6)\\
- & - &- &O(\varepsilon^6)
\end{array}\right]
+
\left[\begin{array}{cccc}
0&O(\varepsilon\omega)& O(\varepsilon^{2}\omega)& O(\varepsilon^{2}\omega)\\
- &O(\varepsilon^{2}\omega)&O(\varepsilon^{3}\omega)& O(\varepsilon^{3}\omega)\\
- & - & O(\varepsilon^{4}\omega)&O(\varepsilon^{4}\omega)\\
- & - &- &O(\varepsilon^{4}\omega)
\end{array}\right]\nonumber
\end{align}
with probability higher than $1-Cn^{-3}$ for some $C>0$ when $n$ is sufficiently large.
Recall that $L^0_{AB}$ is the vector $\alpha\rho^{-1}\nabla\rho$, $L^0_{AC}$ is the vector $\alpha{\bf 1}_d^T$, $L^0_{BB}=\alpha I_d$ where $I_d$ is the identity matrix, $L^0_{BC}=\beta{\bf XY}$, where ${\bf X} = \rho^{-1}\left[\begin{array}{ccc}
    \partial_1\rho &&\\
    &\ddots&\\
    &&\partial_d\rho 
\end{array}\right]$,
${\bf Y} = (2 I_d+{\bf 11}^T)$,
$L^0_{CC}=\beta {\bf Y}$, $L^0_{DD}=\beta I_{\frac{d(d-1)}{2}}$,
and $\omega=\sqrt{\frac{\log(n)}{n\varepsilon^d}}$. The notations ${\bf X}, {\bf Y}$ and the decomposition $L_{BC}^0=\beta{\bf XY}$ are adopted to facilitate the computation of inverse of $L^0$.

\begin{Rmk}\label{remark of comparison of two methods, 1}
It is worthwhile to mention the differences between the reduction trick in Section C and the direct approach in this section. In section C, the matrix $L^0$ has a nonzero $L_{BD}^0$ while the $L^0$ here does not include this term. Instead, we place the $O(\varepsilon^4)$-term in $W_{BD}$, making it easier to compute the inverse of $L^0$. 
This highlights an advantage of our reduction trick, as we can handle more complex $L^0$ without calculating the exact values of all $l_j$ terms in $(L^0)^{-1}$ as in (\ref{inverseL0}). In the direct approach, we must compute the exact values of $(L^0)^{-1}$ and $E$, ensuring many of them cancel out in subsequent steps (cf. Lemma \ref{simplify} and equations (\ref{cancellation1}), (\ref{cancellation2})).
Since $L^0$ here is a bit simpler than the $L^0$ in the reduction-trick based approach, $W$ is more complicated, and we will see the consequences of this in Step 3. 
\end{Rmk}

\subsection{Step 3: The inverse of $\frac{1}{k_z}Z^TZ$}\label{subsection inversion of Z^TZ in the direct}
Recall (\ref{approximation of ZtZ/t}) and denote 
\[
E:=W+{\tt LD}\,;
\]
that is, $\frac{1}{k_z}Z^TZ=L^0+E$.  
To alleviate the heavy notation burden, we denote $\rho:=\rho(z)$ and $\partial_j\rho:=\nabla_j\rho|_{z}$, $j=1,\dots, d$.
Recall $L^0$ defined in \eqref{definition L0}.
Let 
${\bf X} = \rho^{-1}\left[\begin{array}{ccc}
    \partial_1\rho &&\\
    &\ddots&\\
    &&\partial_d\rho 
\end{array}\right]$,
${\bf Y} = \left[\begin{array}{cccc}
    3 & 1 & \cdots & 1\\
    1 & 3 & \ddots& \vdots\\
    \vdots & \ddots& \ddots&1\\
    1 & \cdots & 1 & 3 
\end{array}\right]\in \mathbb{R}^{d\times d}$.
Then the inverse of $L^0$ can be computed directly as
\begin{align}\label{L0-1}
(L^0)^{-1}
=\left[\begin{array}{cccc}
\lambda & 0      & \nu \bf 1^{T} &0\\
   -    & K_{11} & K_{12}        &0\\
   -    & K_{21} & K_{22}        &0\\
   -    & - & - &\frac{1}{\beta}I_{\frac{d(d-1)}{2}}
\end{array}\right],
\end{align}
where $\lambda$ and $\nu$ satisfy $
\left\{\begin{array}{lc}
     \lambda+\nu\alpha d &=1  \\
     \lambda\alpha +\nu\beta(d+2)& =0 
\end{array}\right.,$
\begin{align*}
K_{11} &\,= \frac{1}{\alpha}\left(I+\frac{\beta}{\alpha }{\bf XYQYX}\right),\ \ 
K_{12} = K_{21}^T=-\frac{1}{\alpha }{\bf XYQ},\ \ 
K_{22} =\frac{1}{\beta}{\bf Q}+\frac{\nu^2}{\lambda}{\bf 11}^T,\nonumber\\
{\bf Q}&\,= {\bf Y}^{-1}-F^{-1}+\frac{F^{-1}{\bf 11}^TF^{-1}}{1+{\bf 1}^TF^{-1}{\bf 1}}={\bf Y}^{-1}+O(\varepsilon^2), 
\end{align*}
and $F^{-1}={\bf X}^2(2{\bf X}^2-\frac{\alpha}{\beta}\rho^2I)^{-1}$. By a direct expansion, we have 
\begin{align*}
F^{-1}=-\frac{\beta}{\alpha}\rho^{-2}{\bf X}^2+O(\varepsilon^{4})
\end{align*}
and 
$$
{\bf Y}^{-1} = \frac{1}{2(d+2)}\left[\begin{array}{cccc}
    d+1 & -1 & \cdots & -1\\
    -1 & d+1 & \ddots& \vdots\\
    \vdots & \ddots& \ddots&-1\\
    -1 & \cdots & -1 & d+1 
\end{array}\right]=\frac{1}{2}I-\frac{1}{2(d+2)}{\bf 11}^T.$$
To evaluate the inverse of $\frac{1}{k_z}Z^TZ=L^0+E$, we consider the Taylor's expansion 
\begin{align}\label{L0+E inverse Taylor}
(L^0+E)^{-1}= (L^0)^{-1}-(L^0)^{-1}E(L^0)^{-1}+[(L^0)^{-1}E]^2(L^0)^{-1}-\cdots\,;
\end{align}
but to achieve this expansion, we need the control of the norm of $(L^0)^{-1}E$.

\begin{Rmk}\label{remark of comparison of two methods, 2}
As mentioned at the end of the previous section, $W_{BD}^0$ causes a problem in the direct approach. 
When using the expansion formula \eqref{formula: inversion with approximation} to compute the inverse of $L^0+E$, $W_{BD}^0$ pollutes the terms in the last column of $(L^0)^{-1}E$ in (\ref{omega123}), causing the expansion to converge slowly.
\end{Rmk}

A direct computation gives 
{\scriptsize
\begin{align*}
&(L^0)^{-1}E\\
=&\,\left[\begin{array}{cccc}
\nu {\bf 1}^T E_{CA}
&\lambda E_{AB}+\nu {\bf 1}^TE_{CB}
& \lambda E_{AC}+\nu {\bf 1}^TE_{CC}
& \lambda E_{AD}+\nu {\bf 1}^TE_{CD}\\
K_{11}E_{BA}+K_{12}E_{CA}
&K_{11}E_{BB}+K_{12}E_{CB}
&K_{11}E_{BC}+K_{12}E_{CC}
&K_{11}E_{BD}+K_{12}E_{CD}\\
K_{21}E_{BA}+K_{22}E_{CA}
&\nu{\bf 1}E_{AB}+K_{21}E_{BB}+K_{22}E_{CB}
&\nu{\bf 1}E_{AC}+K_{21}E_{BC}+K_{22}E_{CC}
&\nu{\bf 1}E_{AD}+K_{21}E_{BD}+K_{22}E_{CD}\\
\frac{1}{\beta}E_{DA}
&\frac{1}{\beta}E_{DB}
&\frac{1}{\beta}E_{DC}
&\frac{1}{\beta}E_{DD}
\end{array}\right],
\end{align*}}\\
where we decompose $E$ into blocks in the same way as that in \eqref{decomposition block LD} and use the associated notation.
Recall that $$(L^0)^{-1}=\left[\begin{array}{cccc}
O(1) & 0 & O(\varepsilon^{-2}) &0\\
- &O(\varepsilon^{-2}) &O(\varepsilon^{-2})&0\\ 
- & - & O(\varepsilon^{-4})&0\\
- & - & -    &O(\varepsilon^{-4}) 
\end{array}\right] 
\mbox{ and } 
E=\left[\begin{array}{cccc}
0 &O(\varepsilon^3) & O(\varepsilon^4) & O(\varepsilon^4)\\
- &O(\varepsilon^4) &O(\varepsilon^5)&O(\varepsilon^4)\\ 
- & - & O(\varepsilon^6)&O(\varepsilon^6)\\
- & - & -    &O(\varepsilon^6) 
\end{array}\right]+\mathtt{LD}.$$
Therefore, we have 
$$(L^0)^{-1}E=
\left[\begin{array}{cccc}
O(\varepsilon^2) &O(\varepsilon^3) &O(\varepsilon^4) &O(\varepsilon^4)\\
O(\varepsilon)   &O(\varepsilon^2) &O(\varepsilon^3) &O(\varepsilon^2)\\ 
O(1) & O(\varepsilon) & O(\varepsilon^2)   &O(\varepsilon^2)\\
O(1) & O(1)           & O(\varepsilon^2)   &O(\varepsilon^2) 
\end{array}\right]+\mathtt{LD},$$
and hence
\begin{align}\label{omega123}
\begin{split}
\Omega_1:=&\ (L^0)^{-1}E(L^0)^{-1}=
\left[\begin{array}{cccc}
O(\varepsilon^2) &O(\varepsilon) & O(1) & O(1)\\
O(\varepsilon)   &O(1)           & O(\varepsilon^{-1}) &O(\varepsilon^{-2})\\ 
O(1) & O(\varepsilon^{-1}) & O(\varepsilon^{-2})&O(\varepsilon^{-2})\\
O(1) & O(\varepsilon^{-2}) & O(\varepsilon^{-2})   &O(\varepsilon^{-2}) 
\end{array}\right]+\mathtt{LD},\\
\Omega_2:=&\ [(L^0)^{-1}E]^2(L^0)^{-1}
=
\left[\begin{array}{cccc}
O(\varepsilon^4) & O(\varepsilon^2) & O(\varepsilon^2) & O(\varepsilon)\\
O(\varepsilon^2) & O(1) & O(1) & O(1)\\ 
O(\varepsilon^2) & O(1) & O(1) & O(\varepsilon^{-1})\\
O(\varepsilon)   & O(1) & O(\varepsilon^{-1}) & O(\varepsilon^{-2})
\end{array}\right]+\mathtt{LD}
,\\
\Omega_3:=&\ [(L^0)^{-1}E]^3(L^0)^{-1}=
\left[\begin{array}{cccc}
O(\varepsilon^5) &O(\varepsilon^3) & O(\varepsilon^3) & O(\varepsilon^2)\\
O(\varepsilon^3) &O(\varepsilon^2) &O(\varepsilon^1)&O(1)\\ 
O(\varepsilon^3) & O(\varepsilon) & O(\varepsilon)&O(1)\\
O(\varepsilon^2) & O(1) & O(1)   &O(1) 
\end{array}\right]+\mathtt{LD}, \mbox{ etc}.
\end{split}
\end{align}
Denote $$\Omega_m:=
\left[\begin{array}{cccc}
(\Omega_m)_{AA}&(\Omega_m)_{AB}& (\Omega_m)_{AC}& (\Omega_m)_{AD}\\
(\Omega_m)_{BA} &(\Omega_m)_{BB}&(\Omega_m)_{BC}&(\Omega_m)_{BD}\\ 
(\Omega_m)_{CA} & (\Omega_m)_{CB} & (\Omega_m)_{CC}&(\Omega_m)_{CD}\\
(\Omega_m)_{DA} & (\Omega_m)_{DB} & (\Omega_m)_{DC} &(\Omega_m)_{DD}
\end{array}\right].$$
Note that only
$(L^0)^{-1}$, $\Omega_1$, $(\Omega_2)_{BB}$, and $(\Omega_2)_{DD}$ contain the first and the second lowest order terms, and the remaining terms inside $\Omega_m:=[(L^0)^{-1}E]^m(L^0)^{-1}, m\geq 2,$ are of higher order compared to $(L^0)^{-1}-\Omega_1$. 
To be precise, if we denote 
\begin{align}\label{G}
G:=(L^0)^{-1}-\Omega_1
+\left[\begin{array}{cccc}
0  &  0            &  0  &  0\\
0  &(\Omega_2)_{BB}&  0  &  0\\ 
0  &  0            &  0  &  0\\
0  &  0            &  0  &  (\Omega_2)_{DD}
\end{array}\right],
\end{align}
then the higher order terms are
\begin{align}\label{hot}
(L^0+E)^{-1}-G=
\left[\begin{array}{cccc}
O(\varepsilon^4) &O(\varepsilon^2) & O(\varepsilon^2) & O(\varepsilon)\\
O(\varepsilon^2) &O(\varepsilon^2) &O(1)&O(1)\\ 
O(\varepsilon^2) & O(1) & O(1)&O(\varepsilon^{-1})\\
O(\varepsilon) & O(1) & O(\varepsilon^{-1})   &O(1) 
\end{array}\right].
\end{align}

\begin{Rmk}\label{remark of comparison of two methods, 3}
We shall mention that although non-zero entries in $(L^0)^{-1}$ dominate their corresponding entries in $\Omega_1$ and $\Omega_2$, we cannot discard terms in $\Omega_1$ and some entries in $\Omega_2$ because at the vary last step of the proof when we put everything together, some of those dominated entries will play a critical role via canceling each other. See the proof of Theorem \ref{ThmInterior}, particularly \eqref{higher order residue}, for details. 
\end{Rmk}

{\allowdisplaybreaks
In the following paragraph, $(\Omega_1)_{CA}$ and $(\Omega_1)_{CC}$ in $G$ are computed explicitly because they will be used later, while other terms in $\Omega_1$ and $\Omega_2$ are kept implicit.
Namely, the matrix in (\ref{G}) can be expressed as 
$$G
=\left[\begin{array}{cccc}
G_{AA}&G_{AB}&G_{AC}&G_{AD}\\
G_{BA}&G_{BB}&G_{BC}&G_{BD}\\ 
G_{CA}&G_{CB}&G_{CC}&G_{CD}\\
G_{DA}&G_{DB}&G_{DC}&G_{DD}
\end{array}\right],$$
where the first row is
\begin{align*}
G_{AA} 
=&\ \lambda
-(\Omega_1)_{AA},\\
G_{AB}
=&\ 
-(\Omega_1)_{AB},\\
G_{AC}
=&\ \nu {\bf 1}^T 
-(\Omega_1)_{AC},\\
G_{AD}
=&\ 
-(\Omega_1)_{AD};
\end{align*}
the second row is 
\begin{align*}
G_{BA}
=&\ 
-(\Omega_1)_{BA},\\
G_{BB}
=&\ K_{11} 
-(\Omega_1)_{BB}+(\Omega_2)_{BB},\\
G_{BC}
=&\ K_{12} 
-(\Omega_1)_{BC},\\
G_{BD}
=&\ 
-(\Omega_1)_{BD};
\end{align*}
the third row is
\begin{align*}
G_{CA}
=&\ \nu{\bf 1} -\left( 
\lambda K_{21}E_{BA}+\lambda K_{22}E_{CA}
+\nu{\bf 1}E_{AC}\nu{\bf 1}+K_{21}E_{BC}\nu{\bf 1}+K_{22}E_{CC}\nu{\bf 1}\right),\\
G_{CB}
=&\ K_{21}
-(\Omega_1)_{CB},\\
G_{CC}
=&\ K_{22} -\left( K_{21}E_{BA}\nu{\bf 1}^T+ K_{22}E_{CA}\nu{\bf 1}^T\right)\\
&\ \hspace{6.65mm}-\left[ \left(\nu{\bf 1}E_{AB}+K_{21}E_{BB}+ K_{22}E_{CB}\right)K_{12} + \left(\nu{\bf 1}E_{AC}+K_{21}E_{BC}+K_{22}E_{CC}\right)K_{22}\right],\\
G_{CD}
=&\ 
-(\Omega_1)_{CD};
\end{align*}
and the fourth row is
\begin{align*}
G_{DA}
=&\ -\frac{1}{\beta}\left(\lambda E_{DA}+E_{DC}\nu{\bf 1}\right),\\
G_{DB}
=&\ -\frac{1}{\beta}\left(E_{DB}K_{11}
+E_{DC}K_{21}\right),\\
G_{DC}
=&\ -\frac{1}{\beta}\left(E_{DA}\nu{\bf 1}^T
+E_{DB}K_{12}+E_{DC}K_{22}\right),\\
G_{DD}
=&\ \frac{1}{\beta}I
-(\Omega_1)_{DD}+(\Omega_2)_{DD}.
\end{align*}
}

In the following computation, we use $K_{22}^0$ to denote the leading order terms of $K_{22}=\frac{1}{\beta}{\bf Q}+\frac{\nu^2}{\lambda}{\bf 11}^T=\frac{1}{\beta}{\bf Y}^{-1}+\frac{\nu^2}{\lambda}{\bf 11}^T+O(\varepsilon^{-2})$, i.e., $K_{22}^0=\frac{1}{\beta}{\bf Y}^{-1}+\frac{\nu^2}{\lambda}{\bf 11}^T$. Similarly, $K_{11}^0=\frac{1}{\alpha}I$ and $K_{12}^0 = K_{21}^0=-\frac{1}{\alpha }{\bf X}$. 
For later use, we need the following formulas which can be derived by routine computations, and we omit the details.
\begin{Lma}\label{simplify}
Denote
$\lambda=\frac{(d+2)\beta}{(d+2)\beta-\alpha^2 d}$, 
$\nu=\frac{-\alpha}{(d+2)\beta-\alpha^2 d}$.
Then we have
\begin{enumerate}
\item [\rm (i)] $\displaystyle\frac{\nu}{\lambda}=\frac{-\alpha}{(d+2)\beta}=\varepsilon^{-2}\frac{d+4}{d+2}$.

\item [\rm (ii)] $\displaystyle \nu{\bf 1}+\alpha K_{22}^0{\bf 1}=0$.

\item [\rm (iii)] $\displaystyle 2\beta\frac{\nu^2}{\lambda}{\bf 1}+\alpha\nu{\bf 1}+\beta K_{22}^0{\bf 1}=\frac{1}{d+2}{\bf 1}$.

\item [\rm (iv)] $\displaystyle E_{DB}K_{11}^0\alpha[\nabla_j\eta]_{d\times 1}=\beta\rho^{-1}[\nabla_s\rho\nabla_t\eta+\nabla_t\rho\nabla_s\eta]_{\frac{d(d-1)}{2}\times 1}$, where $\eta$ might be $f$ or $\rho$.

\item [\rm (v)] $\displaystyle E_{DB}K_{12}^0\alpha{\bf 1}=-\beta\rho^{-2}[\nabla_s\rho\nabla_t\rho+\nabla_t\rho\nabla_s\rho]_{\frac{d(d-1)}{2}\times 1}$.
\end{enumerate}
\end{Lma}

\subsection{Step 4: The estimate for $(Z^TZ)^{-1}Z^T{\bf f}$}

In the following, we will compute 
$\frac{1}{k_z}Z^T{\bf f}$
and express it as 
\begin{align*}
\frac{1}{k_z}Z^T{\bf f}=\overline{L^0}+ \overline{W} +\overline{{\tt LD}}. 
\end{align*}
By using similar computation we have done in previous steps, we obtain the following estimation.

\begin{Lma}\label{ZTf} Let $\alpha=\frac{\varepsilon^2}{d+2}$, $\beta=\frac{\varepsilon^4}{(d+2)(d+4)}$
and $\omega=\sqrt{\frac{\log(n)}{n\varepsilon^d}}$. Then
{\allowdisplaybreaks
{\small
\begin{align}\label{Expansion of Z^Tf for the direct approach}
&\ \frac{1}{k_z}Z^T{\bf f}
= (\overline{L^0}+ \overline{W})+\overline{{\tt LD}}\\
=&\ 
\left[\begin{array}{c}
f(z)+\frac{1}{2}\alpha\left(\Delta f(z)+2\rho^{-1}\langle\nabla f,\nabla \rho\rangle(z)\right)+O(\varepsilon^4)\\
\\
\alpha \left(\nabla_1f|_{z} + \rho^{-1}(z)f(z)\nabla_1\rho|_{z}\right)+ O(\varepsilon^3)\\
\vdots\\
\alpha \left(\nabla_df|_{z}+ \rho^{-1}(z)f(z)\nabla_d\rho|_{z}\right)+O(\varepsilon^3)\\
\\
\alpha f(z)+\beta\rho^{-1}(z)\left(\nabla_1^2(f\rho)|_{z}+\frac{1}{2}\Delta (f\rho)(z)-\frac{d
+4}{2(d+2)}f(z)\Delta\rho(z)\right)+\beta\mathcal{U}_{11}f(z)+O(\varepsilon^5)\\
\vdots\\
\alpha f(z)+\beta\rho^{-1}(z)\left(\nabla_d^2(f\rho)|_{z}+\frac{1}{2}\Delta (f\rho)(z)-\frac{d
+4}{2(d+2)}f(z)\Delta\rho(z)\right)+\beta\mathcal{U}_{dd}f(z)+O(\varepsilon^5)\\
\\
\beta\rho^{-1}\nabla_1\nabla_2(f\rho)|_{z} +\beta\mathcal{V}_{1,2}f(z)+O(\varepsilon^5)\\
\vdots\\
\beta\rho^{-1}\nabla_{d-1}\nabla_d(f\rho)|_{z} +\beta\mathcal{V}_{d-1,d}f(z)+O(\varepsilon^5)\\
\end{array}\right]
+
\left[\begin{array}{c}
O(\omega)\\
\\
O(\varepsilon\omega)\\
\vdots\\
O(\varepsilon\omega)\\
\\
O(\varepsilon^{2}\omega)\\
\vdots\\
O(\varepsilon^{2}\omega)\\
\\
O(\varepsilon^{2}\omega)\\
\vdots\\
O(\varepsilon^{2}\omega)
\end{array}\right],\nonumber
\end{align}}}
where $\mathcal{U}_{jj}$ and 
$\mathcal{V}_{s,t}$ are defined in Corollary \ref{intfu}, with probability higher than $1-Cn^{-3}$ for some $C>0$ if $\varepsilon\to 0$  and $\frac{n\varepsilon^d}{\log(n)}\to \infty$ holds when $n\to \infty$.
\end{Lma}

\begin{proof}
The proof is similar to the above, so we sketch it without providing all details. Denote
$$
\frac{1}{k_z}Z^T{\bf f}
=\left[\begin{array}{c}
N_{A}\\
N_{B}\\
N_{C}\\
N_{D}
\end{array}\right]\,,
$$ 
where $
N_{A}=\left[\begin{array}{c}
\displaystyle\frac{1}{k_z}\sum_{j=1}^{k_z} f(x_{z,j})
\end{array}\right]$, $N_{B}=\left[\begin{array}{c}
\displaystyle\frac{1}{k_z}\sum_{j=1}^{k_z} ({\bf q}_j)_1f(x_{z,j})\\
\vdots\\ 
\displaystyle\frac{1}{k_z}\sum_{j=1}^{k_z}({\bf q}_j)_df(x_{z,j})
\end{array}\right]
$,\\
$N_{C}=\left[\begin{array}{c}
\displaystyle\frac{1}{k_z}\sum_{j=1}^{k_z}({\bf q}_j)_1^2f(x_{z,j})\\
\vdots\\
\displaystyle\frac{1}{k_z}\sum_{j=1}^{k_z} ({\bf q}_j)_d^2f(x_{z,j})\\
\end{array}\right]$ and
$N_{D}=\left[\begin{array}{c}
\displaystyle\frac{1}{k_z}\sum_{j=1}^{k_z} ({\bf q}_j)_1({\bf q}_j)_2f(x_{z,j})\\
\vdots\\
\displaystyle\frac{1}{k_z}\sum_{j=1}^{k_z}({\bf q}_j)_{d-1}({\bf q}_j)_df(x_{z,j})
\end{array}\right]$.
The computation and control of $\overline{L^0}$, $\overline{W}$ and $\overline{{\tt LD}}$ follow the same steps as the above.
We first approximate $\frac{1}{k_z}\sum_{j=1}^{k_z} ({\bf q}_j)_lf(x_{z,j})$ by their continuous  counterparts $\overline{L^0}+ \overline{W}$ via the large deviation. For instance, fix any $l=1,\dots,d$, thanks to Theorem \ref{intf} we have
\begin{align*}
& \frac{1}{k_z}\sum_{j=1}^{k_z} ({\bf q}_j)_lf(x_{z,j})\\
=&\ \frac{1}{k_z}\sum_{j=1}^{k_z} \left((\widetilde{\bf q}_j)_l+O(\varepsilon^2)\|x_{z,j}-z\|\right) f(x_{z,j})\\
=&\ -\hspace{-4mm} \int_{\widetilde{B}} \langle x-z, e_l \rangle f(x)d{\rm vol}_{\rho}+ 
-\hspace{-3.8mm} \int_{\widetilde{B}} O(\varepsilon^3)f(x)d{\rm vol}_{\rho}+\textrm{L.D.}\\
=&\ \alpha(\nabla_lf|_{z}+\rho^{-1}f\nabla_l\rho)+ O(\varepsilon^3)+\textrm{L.D.}\,,
\end{align*}
where $\textrm{L.D.}$ means the large deviation of $\frac{1}{k_z}\sum_{j=1}^{k_z} ({\bf q}_j)_lf(x_{z,j})$.
Other bias terms also follow from Corollary \ref{intf} immediately. 
At last, the large deviation terms are computed in Lemma \ref{LD ZTf}. 
\end{proof}

By putting the above calculation together, we have 
\[
(Z^TZ)^{-1}Z^T{\bf f}
=(G+h.o.t.)(\overline{L^0}+\overline{W}+\overline{{\tt LD}})\,,
\]
where $h.o.t.$ means the higher order terms.
We are finally ready to finish the proof of the desired quadratic fitting theorem, which is a version of Theorem \ref{main} when the manifold is boundary free. We note that the same procedure could be applied to the boundary case; however, the computation is significantly more intensive.

\begin{Thm}[cf. Theorem \ref{main}] \label{ThmInterior}
Assume Assumptions \ref{Assumption1}-\ref{Assumption4} hold for the sample data $\{x_i\}_{i=1}^n$. Moreover, assume that $\rho\in C^4(M)$ and $M$ has no boundary.
Let $f:M\to\mathbb{R}$ be a $C^4$-function. 
Consider $z\in M$ and denote the sample points in the $\varepsilon$-neighborhood $B_{\varepsilon}(z)\subset\mathbb{R}^p$ as $x_{z,j}$, where $j=1,\dots, k_z$. 
Denote the base matrix generated by local PCA as $ Z$ and $\omega=\sqrt{\frac{\log(n)}{n\varepsilon^d}}$. Then 
$$(Z^TZ)^{-1}Z^T{\bf f}
=\left[\begin{array}{c}
f(z)+O(\varepsilon^2)+O(\omega)\\
\\
\nabla_1 f|_{z}+O(\varepsilon^2)+O(\varepsilon^{-1}\omega)\\ 
\vdots\\ 
\nabla_d f|_{z}+O(\varepsilon^2)+O(\varepsilon^{-1}\omega)\\
\\
\frac{1}{2}\nabla_1\nabla_1 f|_{z}+O(\varepsilon)+O(\varepsilon^{-2}\omega)\\
\vdots\\
\frac{1}{2}\nabla_d\nabla_d f|_{z}+O(\varepsilon)+O(\varepsilon^{-2}\omega)\\
\\
\nabla_1\nabla_2 f|_{z}
+O(\varepsilon)+O(\varepsilon^{-2}\omega)\\
\vdots\\
\nabla_{d-1}\nabla_d f|_{z}
+O(\varepsilon)+O(\varepsilon^{-2}\omega)
\end{array}\right],
$$ 
with probability higher than $1-Cn^{-3}$ for some $C>0$ 
if $\varepsilon\to 0$ and $\frac{n\varepsilon^d}{\log(n)}\to \infty$ holds when $n\to \infty$.
\end{Thm}

\begin{proof}
The leading terms of $(Z^TZ)^{-1}Z^T{\bf f}$ can be computed separately as follows.

\smallskip

{\bf The constant part:}
The leading $O(1)$-term is
$\lambda f+\nu{\bf 1}^T\alpha f{\bf 1}=f$ because $\lambda+\nu\alpha d=1$. The second order terms are in $O(\varepsilon^2)$. Note that, by Theorem \ref{main}, we know the error term is in fact of $O(\varepsilon^3)$. Here we cannot confirm this unless more computation is involved. Interested readers may try it and verify that all $O(\varepsilon^2)$-terms will be canceled. This again shows the advantage of our approach by using ${\bf f}_{Taylor}$ and $V$ in the proof of Theorem \ref{main}. 
\smallskip

{\bf The  gradient part:}
We multiply the matrix in $\frac{1}{k_z}Z^T{\bf f}$ by $[G_{BA}\ G_{BB}\ G_{BC}\ G_{BD}]$. The leading term is  
$$K_{11}^0[\alpha(\nabla_jf + \rho^{-1}f\nabla_j\rho)]_{d\times 1}+K_{12}^0\alpha f{\bf 1},$$
which can be simplified and becomes the solely 
$[\nabla_jf]_{d\times 1}$.

\smallskip

{\bf The Hessian part:}
Now we multiply the matrix in $\frac{1}{k_z}Z^T{\bf f}$ by $[G_{CA}\ G_{CB}\ G_{CC}\ G_{CD}]$. The product contains many terms and only leading terms are presented here. Recall that $\nu\in O(\varepsilon^{-2})$, $\lambda\in O(1)$, $K_{11}^0,K_{12}^0,K_{21}^0\in O(\varepsilon^{-2})$, $K_{22}^0=\frac{1}{\beta}{\bf Y}^{-1}+\frac{\nu^2}{\lambda}{\bf 11}^T\in O(\varepsilon^{-4})$. Moreover, $K_{22}-K_{22}^0=\frac{1}{\alpha}\rho^{-2}{\bf X}^2+O(1)\in O(\varepsilon^{-2})$. We denote $K_{22}^1:=\frac{1}{\alpha}\rho^{-2}{\bf X}^2$.

\begin{itemize}
\item 
The $O(\varepsilon^{-2})$-term  
$\nu{\bf 1} f +K_{22}^0\alpha f{\bf 1}$  vanishes by Lemma \ref{simplify}(ii).

\item 
The $O(1)$-terms of $[G_{CA}\ G_{CB}\ G_{CC}\ G_{CD}]\cdot \frac{1}{k_z}Z^T{\bf f}$ can be gathered as
{\allowdisplaybreaks
\begin{align*}
&\nu{\bf 1}\alpha\frac{1}{2} \rho^{-1}\left[\Delta(f\rho)-f\Delta\rho\right] -\left( 
\lambda K_{22}^0E_{CA}
+\nu{\bf 1}E_{AC}\nu{\bf 1}+K_{22}^0E_{CC}\nu{\bf 1}\right)f\\
&\quad +K_{21}^0\alpha\rho^{-1}[\nabla_j(f\rho)]_{d\times 1}
 + K_{22}^0\beta\rho^{-1}\left[\nabla_j^2(f\rho)+\frac{1}{2}\Delta (f\rho) -\frac{d+4}{2(d+2)}f\Delta\rho\right]_{d\times 1}\\
&\quad + K_{22}^0\beta[\mathcal{U}_{jj}f]_{d\times 1} +K_{22}^1{\bf 1}\alpha f - \left(  K_{22}^0E_{CA}\nu{\bf 1}^T
+\nu{\bf 1}E_{AC}K_{22}^0+K_{22}^0E_{CC}K_{22}^0\right){\bf 1}\alpha f\\
=&\ \rho^{-1}\beta K_{22}^0\left[\nabla_j^2(f\rho)\right]_{d\times 1}+\frac{1}{2}\rho^{-1}\left( \alpha\nu{\bf 1}+\beta K_{22}^0{\bf 1}\right)\Delta(f\rho)-\frac{1}{2}\rho^{-1}\alpha\nu{\bf 1}f\Delta\rho\\
&\quad +\rho^{-1}\alpha K_{21}^0[\nabla_j(f\rho)]_{d\times 1} -\frac{1}{2}\rho^{-1}\frac{d+4}{d+2}\beta K_{22}^0{\bf 1}f\Delta\rho
+ f\beta K_{22}^0\left[ \mathcal{U}_{jj}\right]_{d\times 1}+\frac{1}{\alpha}\rho^{-2}X^2{\bf 1}\alpha f\\
&\quad - f\left( \lambda K_{22}^0E_{CA}+\nu{\bf 1}E_{AC}\nu{\bf 1}+K_{22}^0E_{CC}\nu{\bf 1}\right)\\
&\quad -f\left(K_{22}^0E_{CA}\nu{\bf 1}^T
+\nu{\bf 1}E_{AC}K_{22}^0+K_{22}^0E_{CC}K_{22}^0\right){\bf 1}\alpha\\
=&\ \rho^{-1} \beta\frac{1}{\beta}{\bf Y}^{-1}\left[\rho\nabla_j^2f+2\nabla_jf\nabla_j\rho+f\nabla_j^2\rho\right]_{d\times 1}
+\rho^{-1} \beta\frac{\nu^2}{\lambda}{\bf 11}^T\left[\nabla_j^2(f\rho)\right]_{d\times 1}\\ 
&\quad +\frac{1}{2}\rho^{-1} \left( \alpha\nu{\bf 1}+\beta K_{22}^0{\bf 1}\right)\Delta(f\rho)-\frac{1}{2}\rho^{-1}\left( \alpha\nu{\bf 1}+\alpha^2K_{22}^0{\bf 1}\right)f\Delta\rho\\
&\quad -\frac{1}{\rho}{\bf X}[\nabla_j(f\rho)]_{d\times 1}
+ f\beta K_{22}^0\left[ \mathcal{U}_{jj}\right]_{d\times 1}+f\rho^{-2}(\nabla_j\rho)^2{\bf 1}\\
&\quad - f\left( K_{22}^0E_{CA}(\lambda+\alpha\nu{\bf 1}^T{\bf 1}) 
+\nu{\bf 1}E_{AC}(\nu{\bf 1}+\alpha K_{22}^0{\bf 1})+K_{22}^0E_{CC}(\nu{\bf 1}+\alpha K_{22}^0{\bf 1})\right).
\end{align*}
Recall that ${\bf Y}^{-1}=\frac{1}{2}I-\frac{1}{2(d+2)}{\bf 11}^T$, $\lambda+\alpha\nu d=1$ and, by Lemma \ref{simplify}, $\nu{\bf 1}+\alpha K_{22}^0{\bf 1}=0$ and 
$$
-\frac{1}{2(d+2)}\Delta (f\rho){\bf 1}+\beta\frac{\nu^2}{\lambda}\Delta(f\rho){\bf 1}+\frac{1}{2} \left( \alpha\nu{\bf 1}+\beta K_{22}^0{\bf 1}\right)\Delta(f\rho)=0.
$$ 
Hence this $O(1)$-term can be simplified as
\begin{align}\label{cancellation1}
& [\frac{1}{2}\nabla_j^2f]_{d\times1} + [\frac{1}{\rho}\nabla_jf\nabla_j\rho]_{d\times1}
+ [\frac{1}{2\rho}f\nabla_j^2\rho]_{d\times1}
-\frac{1}{2(d+2)}\rho^{-1}{\bf 11}^T\left[\nabla_j^2(f\rho)\right]_{d\times 1}\\
&\quad +\rho^{-1}\beta\frac{\nu^2}{\lambda}{\bf 11}^T\left[\nabla_j^2(f\rho)\right]_{d\times 1}
 +\frac{1}{2}\rho^{-1} \left( \alpha\nu{\bf 1}+\beta K_{22}^0{\bf 1}\right)\Delta(f\rho)\nonumber\\
&\quad -\frac{1}{\rho^2}\left[\nabla_j\rho\nabla_j(f\rho)\right]_{d\times 1} 
+ f\beta K_{22}^0\left[ \mathcal{U}_{jj}\right]_{d\times 1} 
+f\rho^{-2}(\nabla_j\rho)^2{\bf 1}
- f K_{22}^0E_{CA}\nonumber\\
=&\ \frac{1}{2}[\nabla_j^2f]_{d\times1} + \frac{1}{\rho}[\nabla_jf\nabla_j\rho]_{d\times1}
+ \frac{1}{2\rho}f[\nabla_j^2\rho]_{d\times1}\nonumber\\
&\quad -\frac{1}{\rho^2}\left[\nabla_j\rho(\rho\nabla_jf+f\nabla_j\rho)\right]_{d\times 1}
+f\rho^{-2}(\nabla_j\rho)^2{\bf 1}
+ f\beta K_{22}^0\left[ \mathcal{U}_{jj}\right]_{d\times 1}
- f K_{22}^0E_{CA}\nonumber\\
=&\ \frac{1}{2}[\nabla_j^2f]_{d\times1} 
+ \frac{1}{2\rho}f[\nabla_j^2\rho]_{d\times1}
+ f\beta K_{22}^0\left[ \mathcal{U}_{jj}\right]_{d\times 1}
- f K_{22}^0E_{CA}\,.\nonumber
\end{align}
Finally, since the leading term of   $E_{CA}$ is $\beta\left[\rho^{-1}\nabla_j^2\rho-\frac{1}{d+2}\rho^{-1}\Delta\rho+\mathcal{U}_{jj}\right]_{d\times 1}$ and $\beta\frac{\nu^2}{\lambda}=\frac{d+4}{4(d+2)}$, this $O(1)$-term becomes 
\begin{align}\label{cancellation2}
&\ \frac{1}{2}[\nabla_j^2f]_{d\times1} 
+ \frac{f}{2\rho}[\nabla_j^2\rho]_{d\times1}
- \rho^{-1}f \left({\bf Y}^{-1}+\beta\frac{\nu^2}{\lambda}{\bf 11}^T\right)\left[\nabla_j^2\rho-\frac{1}{d+2}\Delta\rho\right]_{d\times 1}\\
=&\ \frac{1}{2}[\nabla_j^2f]_{d\times1} 
+ \frac{f}{2\rho}[\nabla_j^2\rho]_{d\times1}
- \rho^{-1}f \left(\frac{1}{2}I+\left(\beta\frac{\nu^2}{\lambda}-\frac{1}{2(d+2)}\right){\bf 11}^T\right)\left[\nabla_j^2\rho-\frac{1}{d+2}\Delta\rho\right]_{d\times 1}\nonumber\\
=&\ \frac{1}{2}[\nabla_j^2f]_{d\times1}
+\left[ 
-\frac{1}{2(d+2)}
+\left(\beta\frac{\nu^2}{\lambda}-\frac{1}{2(d+2)}\right)\left(\frac{2}{d+2}\right)
\right]\Delta\rho{\bf 1}\nonumber
\\
=&\ \frac{1}{2}[\nabla_j^2f]_{d\times1}.\nonumber
\end{align}
}
\end{itemize}

{\allowdisplaybreaks
Now we multiply the matrix in $\frac{1}{k_z}Z^T{\bf f}$ by $[G_{DA}\ G_{DB}\ G_{DC}\ G_{DD}]$. As before, we only demonstrate the leading terms, which are of $O(1)$. Recall that $E_{DA}=\beta[\rho^{-1}\nabla_s\nabla_t\rho+\mathcal{V}_{s,t}]_{\frac{d(d-1)}{2}\times 1}$.

\begin{align*}
&\ -\frac{1}{\beta}\left(\lambda E_{DA}+E_{DC}\nu{\bf 1}\right)f
-\frac{1}{\beta}E_{DB}K_{11}^0(\alpha(\nabla_jf+\rho^{-1}f\nabla_j\rho))\\
&\ -\frac{1}{\beta}\left(E_{DA}\nu{\bf 1}^T+E_{DB}K_{12}^0
+E_{DC}K_{22}^0\right)\alpha f{\bf 1} +[\rho^{-1}\nabla_s\nabla_t(f\rho)+f\mathcal{V}_{s,t}]_{\frac{d(d-1)}{2}\times 1}\\
=&\  -\frac{1}{\beta}f\left(\lambda E_{DA}+\alpha E_{DA}\nu{\bf 1}^T{\bf 1}
+E_{DC}\nu{\bf 1}\underline{-\alpha E_{DB}K_{11}^0\rho^{-1}\nabla_j\rho-\alpha E_{DB}K_{12}^0{\bf 1}}_{=0}+\alpha E_{DC}K_{22}^0{\bf 1}\right) \\
&\ +[\nabla_s\nabla_tf]_{\frac{d(d-1)}{2}\times 1}-\frac{1}{\beta}E_{DB}K_{11}^0\alpha\nabla_jf+\frac{1}{\rho}[\nabla_sf\nabla_t\rho+\nabla_s\rho\nabla_tf+f\nabla_s\nabla_t\rho]_{\frac{d(d-1)}{2}\times 1}+f[\mathcal{V}_{s,t}]_{\frac{d(d-1)}{2}\times 1}\\
=&\  -\frac{1}{\beta}f\left( E_{DA}(\lambda+d\alpha \nu)
+E_{DC}(\nu{\bf 1}+\alpha K_{22}^0{\bf 1})\right) 
+[\nabla_s\nabla_tf]_{\frac{d(d-1)}{2}\times 1}+\frac{1}{\rho}[f\nabla_s\nabla_t\rho]_{\frac{d(d-1)}{2}\times 1}+f[\mathcal{V}_{s,t}]_{\frac{d(d-1)}{2}\times 1}\\
=&\ [\nabla_s\nabla_tf]_{\frac{d(d-1)}{2}\times 1}
+\frac{1}{\rho}[f\nabla_s\nabla_t\rho]_{\frac{d(d-1)}{2}\times 1} 
+f[\mathcal{V}_{s,t}]_{\frac{d(d-1)}{2}\times 1}-\frac{1}{\beta}f E_{DA}\\
=&\ [\nabla_s\nabla_tf]_{\frac{d(d-1)}{2}\times 1}
+\frac{1}{\rho}[f\nabla_s\nabla_t\rho]_{\frac{d(d-1)}{2}\times 1} 
+f[\mathcal{V}_{s,t}]_{\frac{d(d-1)}{2}\times 1}
-f[\rho^{-1}\nabla_s\nabla_t\rho+\mathcal{V}_{s,t}]_{\frac{d(d-1)}{2}\times 1}\\
=&\ [\nabla_s\nabla_tf]_{\frac{d(d-1)}{2}\times 1}.
\end{align*}

To finish the proof, we check the residue terms (\ref{hot}). By a direct calculation, $[(L^0+E)^{-1}-G]\frac{1}{k_z}Z^T{\bf f}$ becomes
\begin{align}\label{higher order residue}
\left[\begin{array}{cccc}
O(\varepsilon^4) &O(\varepsilon^2) & O(\varepsilon^2) & O(\varepsilon)\\
O(\varepsilon^2) &O(\varepsilon^2) &O(1)&O(1)\\ 
O(\varepsilon^2) & O(1) & O(1)&O(\varepsilon^{-1)}\\
O(\varepsilon) & O(1) & O(\varepsilon^{-1})   &O(1) 
\end{array}\right]
\left[\begin{array}{c}
O(1)\\
O(\varepsilon^2)\\
O(\varepsilon^2)\\
O(\varepsilon^4)
\end{array}\right]
=
\left[\begin{array}{c}
O(\varepsilon^4)\\
O(\varepsilon^2)\\
O(\varepsilon^2)\\
O(\varepsilon)
\end{array}\right],
\end{align}
which shows that they are all of higher orders and can be ignored asymptotically.

At last, recalling $(L^0)^{-1}$ from (\ref{L0-1}), we can figure out the large deviation terms which come from 
 $(L^0)^{-1}\overline{{\tt LD}}
 =\left[\begin{array}{c}
      O(\omega) \\
      O(\varepsilon^{-1}\omega)\\
      O(\varepsilon^{-2}\omega)\\
      O(\varepsilon^{-2}\omega)
 \end{array}
 \right]$ and 
$$\left[(L^0)^{-1}{\tt LD}\right]\left[(L^0)^{-1}(\overline{L^0})\right]
=\left[\begin{array}{cccc}
      O(\omega)           & O(\varepsilon\omega) & O(\varepsilon^2\omega)& O(\varepsilon^2\omega)\\
      O(\varepsilon^{-1}\omega)&O(\omega)&O(\varepsilon\omega)&O(\varepsilon\omega)\\
      O(\varepsilon^{-2}\omega)&O(\varepsilon^{-1}\omega)&O(\omega)  &O(\omega)  \\
    O(\varepsilon^{-2}\omega)&O(\varepsilon^{-1}\omega)&O(\omega)  &O(\omega)  
 \end{array}
 \right]
 \left[\begin{array}{c}
      O(1) \\
      O(1)\\
      O(\varepsilon^{-2})\\
      O(1)
 \end{array}
 \right]
 =\left[\begin{array}{c}
      O(\omega) \\
      O(\varepsilon^{-1}\omega)\\
      O(\varepsilon^{-2}\omega)\\
      O(\varepsilon^{-2}\omega)
 \end{array}
 \right].$$
Therefore, the large deviation terms are 
$\left[\begin{array}{c}
      O(\omega) \\
      O(\varepsilon^{-1}\omega)\\
      O(\varepsilon^{-2}\omega)\\
      O(\varepsilon^{-2}\omega)
 \end{array}
 \right]$, i.e., the claimed term $\left[\begin{array}{c}
      O\left(\sqrt{\frac{\log(n)}{n\varepsilon^d}}\right) \\
      O\left(\sqrt{\frac{\log(n)}{n\varepsilon^{d+2}}}\right)\\
      O\left(\sqrt{\frac{\log(n)}{n\varepsilon^{d+4}}}\right)\\
      O\left(\sqrt{\frac{\log(n)}{n\varepsilon^{d+4}}}\right)
 \end{array}
 \right]$. We thus finish the proof.
 }
\end{proof}


\begin{thebibliography}{99999999}


\bibitem[AAL23]{AL19}
E. Aamari, C. Aaron, and C. Levrard, {\it Minimax boundary
estimation and estimation with boundary}, Bernoulli, {\bf 29}(4), pp. 3334–3368, 2023.

\bibitem[AF03]{AF03}
R. A. Adams and J. F. Fournier,
{\it Sobolev Spaces}, Second Edition,
Pure and Applied Mathematics, {\bf 140}, 2003. 

\bibitem[CW13]{CW13}
M.-Y. Cheng and H.-T. Wu, 
{\it Local linear regression on manifolds and its geometric interpretation}, 
J. Amer. Stat. Assoc. {\bf 108}(504),
pp. 1421--1434, 2013. 


\bibitem[CWY99a]{CWY99a}
S.-Y. A. Chang, L. Wang, and P. C. Yang,
{\it Regularity of harmonic maps}, 
Comm. Pure Appl. Math. {\bf 52}(9), 
pp. 1099--1111, 1999.


\bibitem[CWY99b]{CWY99b}
S.-Y. A. Chang, L. Wang, and P. C. Yang,
{\it A regularity theory of biharmonic maps}, 
Comm. Pure Appl. Math. {\bf 52}(9), 
pp. 1113--1137, 1999.

\bibitem[CY97]{CY97}
S.-Y. A. Chang and P. C. Yang,
{\it On uniqueness of solutions of n-th order differential equations in conformal geometry}, Math. Reas. Lett. {\bf 4}(1), 
pp. 91--102, 1997.



\bibitem[CZXC20]{CZXC20}
L. Cao, J. Zhao, J. Xu, S. Chen, G. Liu, S. Xin, Y. Zhou, Y. He, {\it Computing smooth quasi-geodesic distance field (QGDF) with quadratic programming}, Computer-Aided Design {\bf 127}:102879, 2020.

\bibitem[D91]{D91}
A. S. Dvinsky, 
{\it Adaptive grid generation from harmonic maps on Riemannian manifolds}, 
J. Comput. Phys. {\bf 95}, 
pp. 450--476, 1991.

\bibitem[DG03]{DG03}
D. L. Donoho and C. Grimes, 
{\it Hessian eigenmaps: Locally linear embedding techniques for high-dimensional data}, 
Proc. Natl. Acad. Sci. {\bf 100}(10), 
pp. 5591--5596, 2003.


\bibitem[DHL00]{DHL00}
Z. Djadli, E. Hebey, and M. Ledoux, 
{\it Paneitz-type operators and applications}, Duke Math. J. {\bf 104}(1), 
pp. 129--169, 2000.



\bibitem[FG96]{FG96}
J. Fan and I. Gijbels, 
{\it Local polynomial modeling and its applications}, 
Chapman \& Hall, 1996.

\bibitem[FN19]{FN19}
L. Fang and M. K. Ng, 
{\it Image colorization by using graph bi-Laplacian}, 
Advances in Computational Math. {\bf 45}, 
pp. 1521--1549, 2019.

\bibitem[FY03]{FY03}
J. Fan and Q. Yao, 
{\it Nonlinear Time Series: Nonparametric and Parametric Methods}, Springer Series in Statistics (SSS), Springer-Verlag New York, 2003.


\bibitem[G73]{G73}
A. Gray, 
{\it The volume of a small geodesic ball of a Riemannian manifold}, 
Michigan Math. J. {\bf 20}, 
pp. 329--344, 1973.

\bibitem[G04]{G04}
A. Gray
{\it Tubes},
Second Edition,
Springer, Basel, 2004.




\bibitem[GK90]{GK90}
C. P. Gupta and Y. C. Kwong, 
{\it Biharmonic eigenvalue problems and $L^p$ estimates}, 
Int. J. Math. Math. Sci. {\bf 13}(3), 
pp. 469--480, 1990.

\bibitem[GT01]{GT01}
D. Gilbarg and N. S. Trudinger,
{\it Elliptic Partial Differential Equations of Second Order}, Second Edition, 
Springer Berlin, Heidelberg, 2001.

\bibitem[H06]{H06}
A. Henrot, 
{\it Extremum problems for eigenvalues of elliptic operators},
Frontiers in Math., Birkhauser, Basel, 2006.

\bibitem[HWAG09]{HWAG09}
Q.-X. Huang, M. Wicke, B. Adams, L. Guibas,
{\it Shape decomposition using modal analysis},
Comput. Graph. Forum, {\bf 28} (2), 407--416, 2009.



\bibitem[ITZ22]{ITZ22}
M. Irakarama, M. Thierry-Coudon, M. Zakari, G. Caumon. {\it Finite element implicit 3D subsurface structural modeling}. Computer-Aided Design. {\bf 149}: 103267, 2022.


\bibitem[KM14]{KM14}
D.N. Kaslovsky, F.G. Meyer, {\it Non-asymptotic analysis of tangent space perturbation}. Information and Inference: a Journal of the IMA. {\bf 3}(2):134--187, 2014.

\bibitem[KSH09]{KSH09}
K. Kim, F. Steinke, and M. Hein, 
{\it Semi-supervised regression using hessian energy with an application to semi-supervised dimensionality reduction},
Advances in Neural Information Processing Systems {\bf 22}, 
pp. 979--987, 2009.

\bibitem[Lee12]{Lee12}
J. Lee,
{\it Introduction to Riemannian manifolds}, 2ed. 
GTM {\bf 176}, 
Springer, 2018.


\bibitem[Li12]{Li12}
P. Li,
{\it Geometric analysis}.
Cambridge Studies in Advanced Mathematics {\bf 134}, 
Cambridge, 2012.






\bibitem[LTZ01]{LTZ01}
R. Li, T. Tang, and P. Zhang,
{\it Moving mesh methods in multiple dimensions based on harmonic maps}, 
J. Comput. Phys. {\bf 170}, 
pp. 562--588, 2001.



\bibitem[LC21]{LC21}
L. Lin and C.-W. Chen,
{\it A new locally linear embedding scheme in light of Hessian eigenmap}, 
arXiv:2112.09086, 2021.




\bibitem[Mo05]{Mo05}
F. Morgan,
{\it Manifolds with density},
Notices of the AMS {\bf 52}(8),
pp. 853--858, 2005.


\bibitem[Mu94]{Mu94}
D. Mumford, 
{\it Elastica and computer vision},
In: Bajaj, C.L. (eds) Algebraic Geometry and its Applications. Springer, NY, 1994. 

\bibitem[Mu02]{Mu02}
D. Mumford,
{\it Pattern theory: the mathematics of perception},
ICM 2002 Vol. III\ 1–3.



\bibitem[MO06]{MO06}
S. Montaldo and C. Oniciuc, 
{\it A short survey of biharmonic maps between Riemannian manifolds}, 
Rev. Un. Mat. Argentina {\bf 47}(2), 
pp. 1--22, 2006. 


\bibitem[MWZ10]{MWZ10}
S. Mukherjee, Q. Wu, and D.-X. Zhou, 
{\it Learning gradients on manifolds}, 
Bernoulli {\bf 16}(1), 
pp. 181--207, 2010.

\bibitem[N11]{N11}
L. Nguyen, 
{\it Singular harmonic maps and applications to general relativity}, 
Comm. Math. Phys. {\bf 301}, 
pp. 411--441, 2011.


\bibitem[P16]{P16}
P. Petersen,
{\it Riemannian geometry}, 3ed.
GTM {\bf 171}, Springer, 2016.


 
\bibitem[PV20]{PV20}
S. Pigola and G. Veronelli,
{\it The smooth Riemannian extension problem},
Ann. Sc. Norm. Super. Pisa Cl. Sci., {\bf XX}(2020), 1507--1551.

\bibitem[S87]{S87}
N. S\`anchez, 
{\it Harmonic maps in general relativity and quantum field theory}, 
in {\it Harmonic Mappings, Twistors and Sigma - Models},
Adv. Series in Math. Phys., World Scientific Publ. Co., 
pp. 270--305, 1987.


\bibitem[S06a]{S06a}
K.-T. Sturm, 
{\it On the geometry of metric measure spaces. I},
Acta Math. {\bf 196}, 
pp. 65--131, 2006.

\bibitem[S06b]{S06b}
K.-T. Sturm, 
{\it On the geometry of metric measure spaces. II},
Acta Math. {\bf 196}, 
pp. 133--177, 2006.




\bibitem[SG64]{SG64} 
A. Savitzky and M. J. E. Golay,  
{\it Smoothing and differentiation of data by simplified least squares procedures}, 
Analytical Chemistry {\bf 36}(8), 
pp. 1627--1639, 1964.

\bibitem[SGWJ18]{SGWJ18}
O. Stein, E. Grinspun, M. Wardetzky, and A. Jacobson, 
{\it Natural boundary conditions for smoothing in geometry processing}, 
ACM Transactions on Graphics {\bf 37}(2), 
pp. 1--13, 2018.

\bibitem[SJWG20]{SJWG20}
O. Stein, A. Jacobson, M. Wardetzky, and E. Grinspun, 
{\it A smoothness energy without boundary distortion for curved surfaces}, 
ACM Transactions on Graphics {\bf 39}(3)
pp. 1--17, 2020.

\bibitem[SW12]{SW12}
A. Singer and H.-T. Wu, 
{\it Vector diffusion maps and the connection Laplacian}, 
Comm. Pure Appl. Math. {\bf 65}(8), 
pp. 1067--1144, 2012.


\bibitem[T05]{T05}
T. Lamm, 
{\it Biharmonic maps},
Doctoral thesis, 
University of Freiburg, 2005.


\bibitem[WGMM10]{WGMM10}
Q. Wu, J. Guinney, M. Maggioni, and S. Mukherjee, 
{\it Learning gradients: predictive models that infer geometry and statistical dependence}, 
J. Mach. Learn. Res. {\bf 11}, 
pp. 2175--2198, 2010.

\bibitem[WW23]{WW23}
H.-T. Wu, N. Wu, {\em When locally linear embedding hits boundary}, 
J. Mach. Learn. Res. {\bf 24}, 
pp. 1--80, 2023.

\end{thebibliography}
\end{document}